\documentclass[twoside,11pt]{article}

\usepackage{jmlr2e}

\usepackage[figuresright]{rotating}
\usepackage{algorithm}
\usepackage{algorithmic}

\usepackage{amsmath}

\DeclareMathOperator*{\argmin}{\arg \min}
\DeclareMathOperator*{\argmax}{\arg \max}
\newcommand{\dom}{\textnormal{dom}}

\newcommand{\layer}{L}
\newcommand{\loss}{\ell}
\newcommand{\bregmanloss}[1][\Psi]{B_{#1}}
\newcommand{\R}{\mathbb{R}}
\newcommand{\N}{\mathbb{N}}
\newcommand{\prox}[1][\Psi]{\textnormal{prox}_{#1}}

\newcommand{\Thetas}{\mathbf{\Theta}}
\newcommand{\Xs}{\mathbf{X}}
\newcommand{\Zs}{\mathbf{Z}}

\ShortHeadings{Lifted Bregman Learning}{Wang and Benning}
\firstpageno{1}

\begin{document}

\title{Lifted Bregman Training of Neural Networks}

\author{\name Xiaoyu Wang \email xw343@cam.ac.uk \\
       \addr Department of Applied Mathematics and Theoretical Physics\\
       University of Cambridge \\
       Cambridge, CB3 0WA, UK
       \AND
       \name Martin Benning \email m.benning@qmul.ac.uk \\
       \addr School of Mathematical Sciences\\
       Queen Mary University of London\\
       London, E1 4NS, UK}

\editor{}

\maketitle

\begin{abstract}%   <- trailing '%' for backward compatibility of .sty file
    We introduce a novel mathematical formulation for the training of feed-forward neural networks with (potentially non-smooth) proximal maps as activation functions. This formulation is based on Bregman distances and a key advantage is that its partial derivatives with respect to the network’s parameters do not require the computation of derivatives of the network's activation functions. Instead of estimating the parameters with a combination of first-order optimisation method and back-propagation (as is the state-of-the-art), we propose the use of non-smooth first-order optimisation methods that exploit the specific structure of the novel formulation. We present several numerical results that demonstrate that these training approaches can be equally well or even better suited for the training of neural network-based classifiers and (denoising) autoencoders with sparse coding compared to more conventional training frameworks.
    
\end{abstract}


\section{Introduction} 
Deep neural networks (DNNs) are extremely popular choices of model functions for a great variety of machine learning problems (cf. \cite{Goodfellow-et-al-2016}). The predominant strategy for training DNNs is to use a combination of gradient-based minimisation algorithm and the back-propagation algorithm \citep{rumelhart1986learning}. Thanks to modern automatic differentiation frameworks, this approach is easy to implement and yields satisfactory results for a great variety of machine learning applications. However, this approach does come with numerous drawbacks. 
First of all, many common activation functions in neural network architectures are not differentiable, which means that subgradient- instead of gradient-based algorithms are required. And even though there are mathematical justifications for such an approach (cf. \cite{bolte2020mathematical}), subgradient-based methods have the disadvantage of inferior convergence rates over other non-smooth optimisation methods like proximal gradient descent (cf. \cite{teboulle2018simplified}). Another drawback is that back-propagation suffers from vanishing gradient problems \citep{erhan2009difficulty}, where gradients of loss functions with respect to the parameters of early network layers become vanishingly small, which renders the use of gradient-based minimisation methods ineffective. While the use of non-saturating activation functions and the introduction of skip-connections to network architectures have been proposed to mitigate vanishing gradient problems, achieving state-of-the-art performance still requires meticulous hyper-parameter tuning and good initialisation strategies. 
Another issue of the back-propagation algorithm is that it is sequential in nature, which makes it difficult to distribute the computation of the network parameters across multiple workers efficiently. As a consequence of the aforementioned issues, numerous distributed optimisation approaches have been proposed as alternatives to gradient-based training in combination with back-propagation in recent years \citep{carreira2014distributed, taylor2016training, zhang2017convergence, askari2018lifted, li2019lifted, zach2019contrastive, gu2020fenchel, hoier2020lifted}. However, many of these approaches still suffer from limitations such as differentiating non-differentiable activation functions \citep{carreira2014distributed}, recovering affine-linear networks \citep{askari2018lifted}, or overly restrictive assumptions \citep{li2019lifted, gu2020fenchel}.

In this work, we propose a distributed optimisation framework for the training of parameters of feed-forward neural networks that does not require the differentiation of activation functions, that does train truly non-linear DNNs, that can be optimised with a variety of deterministic and stochastic first-order optimisation methods and that does come with extensive mathematical foundations. The proposed approach is a generalisation of the Method of Auxiliary Coordinates with Quadratic Penalty \citep{carreira2014distributed}; instead of quadratic penalty functions, we propose to use a novel penalty function based on a special type of Bregman distance, respectively Bregman divergence \citep{bregman1967relaxation}.

The main contributions of this paper are: 1) the proposal of a new loss function that, when differentiated with respect to its second argument, does not require the differentiation of an activation function, 2) its use as a penalty function within the method of auxiliary coordinates, 3) a detailed mathematical analysis of the new loss function, 4) the proposal of a variety of different deterministic and stochastic iterative minimisation algorithms for the empirical risk minimisation of DNNs, 5) the comparison of these iterative minimisation algorithms to more conventional approaches such as stochastic gradient descent and back-propagation, 6) the show-casing of the proposed approach for the training of (denoising) autoencoders with sparse codes, and 7) a demonstration that contrary to wide-spread belief, expressive feed-forward networks can overfit to training data without generalising well when non-standard optimisation methods are being used.

The paper is organised as follows. We describe the state-of-the-art approach of training neural networks with gradient-based algorithms and back-propagation, before we provide an overview over recent developments in distributed optimisation methods for the training of DNNs in Section \ref{sec:training-nns}. In Section \ref{sec:lifted-bregman-training}, we introduce our proposed lifted Bregman framework. We first establish mathematical foundations, before we define the Bregman loss function that replaces the quadratic penalty function in the method of auxiliary coordinates, and verify a whole range of mathematical properties of this loss function that guarantee the advantages over other approaches. In Section \ref{sec:numerical-realisation}, we provide an overview over a range of deterministic and stochastic optimisation approaches that are able to computationally solve the lifted Bregman training problem thanks to results provided in Section \ref{sec:lifted-bregman-training}. We also discuss suitable regularisation strategies for the regularisation of network parameters and outputs of hidden network layers. In Section \ref{sec:example-problems}, we discuss the example problems of classification, data compression and denoising. For the latter two, we introduce a regularised empirical risk minimisation approach that will produce autoencoders with sparse codes. In Section \ref{sec:numerical-results}, we provide numerical results for the example problems described in Section \ref{sec:example-problems} and extensive comparisons to other minimisation approaches. In Section \ref{sec:conclusions-outlook}, we conclude with a summary of the findings and an outlook of future developments.

%overview of neural networks and  about different strategies of training neural networks, such as more distributed optimisation for training of neural networks, as opposed to classical use of gradient based approaches back propagation. In section 2 and 3 we will begin with a discussion of the proposed Bregman loss function and the Bregman lifted learning approach we pursue, as oppose to more traditional lifted approaches. In section 4 and 5 we will discuss example applications and numerical validations.

%This approach is based on the representing a non-decreasing activation function as the argmin of an appropriate convex optimization problem. 

%Experiments indicate that the proposed models provide excellent … for weights for standard neural networks. 

\section{Training neural networks}\label{sec:training-nns}

Traditionally, the predominant strategy for training neural networks is the use of a gradient-based minimisation algorithm in combination with the back-propagation algorithm. Alternatively, distributed optimisation techniques for training neural networks have received growing attention in recent years (cf. \cite{carreira2014distributed}). In the following sections, we revisit the notion of feed-forward neural networks before we summarise the classical approach of gradient-based training and back-propagation. We then discuss distributed optimisation and lifted training strategies.

%begin with feed-forward neural network，then the learning of network parameters，then empirical risk minimisation this is done by back propagation with stochastic gradient descent or variants, but drawbacks: 

\subsection{Feed-forward neural networks}

%Extending the idea from a simple artificial neuron model to multi-layer perceptron architectures is straightforward and natural. 
A feed-forward neural network $\mathcal{N}:\R^n \times \mathcal{P} \rightarrow \R^m$ with $\layer$ layers can be defined as the composition of parametrised functions of the form
%\[f(x) =  W_{\layer}^\top(\sigma_{\layer-1}(W_{\layer-1}^\top \dots\sigma(W_{l}^\top \dots \sigma_{1}(W_{1}^\top x + b_{1}) \dots+ b_{l})  \dots + b_{\layer-1}) ) + b_{\layer} \] 
\begin{align}
    \mathcal{N}(x,\mathbf{\Theta}) = 
    \sigma_{\layer}(f(\sigma_{\layer-1}(f(\ldots \sigma_{1}(f(x_{1},\Theta_1))\ldots)) ,\Theta_{\layer})), \label{eq:feed-forward-nn}
\end{align}
for given input data $x \in \R^n$ and parameters $\mathbf{\Theta} \in \mathcal{P}$. Here $\{ \sigma_{l} \}_{l = 1}^\layer$ denotes the collection of nonlinear activation functions and $f$ denotes a generic function parametrised by parameters $\{ \Theta_l \}_{l = 1}^\layer$. To ease the notation, we use $\mathbf{\Theta}$ to denote $\{ \Theta_l \}_{l = 1}^\layer$. As an example, a rectified linear unit (ReLU) \citep{nair2010rectified} can be constructed via the composition of $\sigma_l(\cdot) := \text{max}(0, \cdot)$ and $f(x, \Theta_l) = W^\top_l x + b_l$ for $\Theta_l = (W_l,b_l)$, where the matrix $W_l$ is usually referred to as the \emph{weights} and the vector $b_l$ as the \emph{bias term} of layer $l$. We refer to the output of the neural network model as $\hat{y}$, i.e. $\hat{y} := \mathcal{N}(x,\mathbf{\Theta})$. 

%Note that feed-forward networks often include so-called \emph{bias terms} that are added to the output of each layer. For simplicity of notation, we have omitted these terms as they can also be incorporated into the weights if the input data and the outputs of each layer are augmented, i.e. by using $\tilde x = (1, x)^\top$ instead of $x$.

\subsection{Training feed-forward networks with first-order methods and back-propagation}
Estimating parameters $\mathbf{\Theta}$ of a feed-forward network in the supervised learning setting is usually framed as the minimisation of an empirical risk of the form
\begin{equation}
     \min_{\Thetas} \frac{1}{s} \sum_{i=1}^{s} \loss(y^i, \hat{y}^i) \, ,
\end{equation}\label{eq:mlp-objective}%
where $\loss$ is a chosen data error term pertaining to the learning task (usually referred to as \emph{loss}), while $\{(x^i, y^i)\}_{i=1}^s$ denotes the $s$ pairs of input and output samples that have to be provided a-priori. %and $R_l$ is the regularisation term pertaining to the $l$-th layer on the weight parameters, controlled by $\alpha_l$. 
%Consider the supervised learning setting where $(x^i,y^i)$ are input-output data pairs, we can write the learning problem down more compactly using the following matrix notation. With $X = [x1, \dots, x_s]$ denoting the input data matrix and the target output matrix by $Y = [y_1, \dots, y_s]$, we then have the problem stating as (omitting the bias terms):
%\begin{equation}
%    \argmin_{\{W\}_{l=1}^{\layer}} \loss(Y - \hat{Y}) + \sum_{l=1}^{\layer} \alpha_l R_l(W_l) \,,
%\end{equation}\label{eq:mlp-objective-matrix-notation}
%On the other hand, one could also unroll the network and receive the Solving for \eqref{eq:mlp-constrained} is equivalent to solving the original problem \eqref{eq:mlp-nested}, in the sense that minimisers to \eqref{eq:mlp-constrained} also minimise \eqref{eq:mlp-nested}. \xw{add ref, Carpinan,Wang}
A standard computational approach to solve \eqref{eq:mlp-objective} are (sub-)gradient-based algorithms such as (sub-)gradient descent and stochastic variants of it. Gradient descent for solving \eqref{eq:mlp-objective} reads
\begin{equation}
\mathbf{\Theta}^{k+1} = \mathbf{\Theta}^{k} - \tau_{\mathbf{\Theta}} \nabla_{\mathbf{\Theta}} \left( \frac{1}{s} \sum_{i = 1}^{s} \ell(y^i,\mathcal{N}(x^i,\Thetas^{k})) \right) \, ,\label{eq:gradient-descent}
\end{equation}
for $k \in \N$, a step-size parameter $\tau_{\mathbf{\Theta}} > 0$ and initial parameters $\mathbf{\Theta}^0$. The evaluation of the gradient $\nabla_{\mathbf{\Theta}} \left( \frac{1}{s} \sum_{i = 1}^{s} \ell(y^i,\mathcal{N}(x^i,\Thetas^{k})) \right)$ is computed with the back-propagation (cf. \cite{rumelhart1986learning}), which is summarised in Algorithm \ref{alg:backprop-algorithm} in Appendix \ref{sec:backprop}. As pointed out in \cite{lecun1988theoretical}, the back-propagation algorithm can also be deduced from a Lagrangian formulation of \eqref{eq:mlp-objective}, where solving the corresponding optimality system leads to the individual steps in the back-propagation algorithm. For each pair of data samples $(x^i,y^i)$ we can define
%\begin{equation*}
%    x_{l} = \sigma_{l} (W_l^\top x_{l-1}^{i}) \; \text{ for $l = 1, 2, \dots, \layer$} \, ,
%\end{equation*}
\begin{equation}
    x_{l}^{i} = \sigma_{l}(f(x_{l-1}^{i},\Theta_{l})) \qquad \text{for} \qquad l = 1, \ldots, \layer \, ,\label{eq:nonlinear-constraint}
\end{equation}
as the so-called \emph{activation variables} with initial value $x_0^{i} = x^{i}$. Further, we use $z_l^{i}$ to denote \emph{transformed input variables} by
%\begin{equation*}
%    z_{l} = W_l^\top x_{l-1}\; \text{ for $l = 1, 2, \dots, \layer$} \, ,
%\end{equation*}
\begin{equation*}
    z_{l}^{i} = f(x_{l-1}^{i},\Theta_{l})\qquad \text{for} \qquad l = 1, \ldots, \layer \, .
\end{equation*}%
To ease notation, we introduce $\Xs = \{ x_l^{i} \}_{l=1,\ldots, \layer}^{i=1,\ldots, s}$ and $\Zs = \{ z_l^{i} \}_{l=1,\ldots, \layer}^{i=1,\ldots, s}$ to denote the two groups of auxiliary variables.  
Having introduced $\Xs$ and $\Zs$, we can write the minimisation problem \eqref{eq:mlp-objective} in the following equivalent constrained form:
\begin{equation}
\begin{split}
    & \text{\;\;\;\;\; }\min_{\Thetas,\Xs,\Zs} \frac{1}{s} \sum_{i=1}^{s} \ell(y^i, x_{\layer}^{i})   \\
    & \text{\;\;\;\;\;\;\;s.t. } z_{l}^{i} = f(x_{l-1}^{i},\Theta_{l}) \,\text{  for $l=1, \ldots, \layer$} \,, \\
    & \text{\;\;\;\;\;\;\;s.t. } x_{l}^{i} = \sigma_{l}(f(x_{l-1}^{i},\Theta_{l})) \,\text{  for $l=1, \ldots, \layer$} \,.
\end{split}\label{eq:mlp-constrained}
\end{equation} 
A derivation of Algorithm \ref{alg:backprop-algorithm} from the Lagrangian formulation of \eqref{eq:mlp-constrained} is included in Appendix \ref{sec:backprop}.

Despite its popularity for training neural networks, the approach of using (sub-)gradient based methods such as \eqref{eq:gradient-descent} in combination with Algorithm \ref{alg:backprop-algorithm} nevertheless suffers from some some drawbacks. Two of the major problems are the vanishing and exploding gradient issues, where the gradients either decrease to vanishingly small values or increase to very large values, which causes major problems for the computation of \eqref{eq:gradient-descent} (cf. \cite{glorot2010understanding,bengio1994learning}). 

A third issue is that the back-propagation algorithm is sequential in nature and computation of individual partial derivatives cannot easily be distributed. %This is due to the difficulty to parallelise the computation of the gradients via back-propagation. Since the forward and backward pass require sequential updates of the parameters $\Xs$, $\Zs$ and $\delta$. 
A fourth drawback is that while it is easy to generalise gradient-based algorithms to proximal gradient methods that can handle non-smooth functions acting on the network parameters, it is less straight-forward to generalise these techniques so that they can handle non-smooth functions acting on the activation variables or the transformed input variables without using subgradients.

%Despite the popularity of the BP algorithm in combination with , it nevertheless suffers from the following drawbacks: 1) vanishing or exploding gradient problems during training often slow down convergence \cite{glorot2010understanding,bengio1994learning}, 2) and back-propagation does not easily parallelise \cite{taylor2016training,zach2019contrastive}.

These drawbacks and limitations of the combination of (sub-)gradient based-methods and back-propagation have motivated research to seek for alternative training methods for DNNs. In the next section, we will summarise some of these alternatives.  %therefore summarise distributed optimisation approaches of training neural networks. 

%alternatively, there's another line of work on lifting methods, lifting neural network, introduce equivalent constrained problem, formulate penalty formulation, how this is trained? advantage? parallel, also don't need to differentiate, projected methods, drawbacks: optimality system, essentially learned linear networks, a bit more stress on this part

%Similar to the idea of generalising perceptron learning, that intuition can be naturally extended to the learning of multi-layer neural network parameters as well. 

\subsection{Distributed optimisation approaches}\label{sec:distributed}

Learning network parameters via distributed optimisation has been given much attention in recent years \citep{carreira2014distributed, taylor2016training, zhang2017convergence, askari2018lifted, li2019lifted, zach2019contrastive, gu2020fenchel}. In these works, the overall Learning Problem \eqref{eq:mlp-objective} is reformulated to Problem \eqref{eq:mlp-constrained}. %The introduction of auxiliary variables $\Xs$ and $\Zs$ break down the nested dependency of the function mapping $\mathcal{N}$ into layer-wise equality constraints.

%Note that the combination of back-propagation and (sub-)gradient-based algorithm can also be viewed as solving \eqref{eq:mlp-constrained} implicitly. The forward-passes of the network ensure the equality constraints in \eqref{eq:mlp-constrained} are met, and in the mean time $x_{l}$ variables are stored to be used later in computing the gradients. The weight parameters are then updated via back-propagation by carrying out a (sub-) gradient descent step.
%$x_{\layer}$ variables are recursively eliminated via forward passes and thus satisfying the non-linear constraints. 

%activation variables $\Xs$ have been decoupled from the network, which allows for more flexibility to use alternative optimisation methods and opens up new paths to solve sub-problems in parallel. %and introduced into the optimisation problem as auxiliary variables.

This constrained optimisation Problem \eqref{eq:mlp-constrained} can further be relaxed by replacing the constraints for the activation variables with penalty terms in the objective function. Solving each layer-wise sub-problem can then be performed independently, which allows for distributed optimisation methods to be utilised. This attempt has mainly been approached from two directions: by using quadratic penalties \citep{carreira2014distributed, taylor2016training} or by framing activation functions as orthogonal projections on convex sets of constraints \citep{zhang2017convergence, askari2018lifted, li2019lifted, zach2019contrastive, gu2020fenchel}. 

\subsubsection{Method of Auxiliary Coordinates with Quadratic Penalty  (MAC-QP)}
The work of \cite{carreira2014distributed} is among the earliest to look into the quadratic penalty approach, where the authors propose the Method of Auxiliary Coordinates (MAC) that relaxes the Learning Problem \eqref{eq:mlp-objective} by replacing the constraints $x_l^i = \sigma_l(f(x_{l - 1}^i, \Theta_l))$ with a Quadratic Penalty (QP) objective where the corresponding minimisation problem reads
\begin{align}\label{eq:mlp-quadratic-penalty}
    & \min_{\Xs,\Thetas} \frac{1}{s} \sum_{i=1}^{s} \left[ \ell(y^i, f(x_{\layer - 1}^i, \Theta_{\layer})) + \sum_{l=1}^{\layer-1} \|x^i_{l} - \sigma_l(f(x_{l - 1}^i, \Theta_l)) \|_2^2 \right] \, .
\end{align}
%The QP objective can efficiently be minimised alternatingly for $\Xs$ and $\Thetas$ where \cite{carreira2014distributed} propose to solve each nonlinear sub-problem with the Gauss-Newton algorithm. 

%In addition, a post-processing step is also performed. After manually setting $x_l^i = \sigma_l(f(x_{l - 1}^i, \Theta_l))$ to eliminate auxiliary variables and satisfy the constraints, the last layer weight parameters are refit to the target data while all other layers' weight parameters being fixed to further reduce the nested error. 

If we consider the special case of affine-linear functions $f(x_{l-1},\Theta_{l}) = W_{l}^\top x_{l-1} + b_l$, then for $s = 1$ and a single training pair of samples $x$ and $y$ the optimality condition w.r.t the bias parameter $b_{l}$ reads
\begin{equation*}
0 = ( x_{l} - \sigma_{l}(W_{l}^\top x_{l-1}+b_{l})) \sigma_{l}^\prime (W_l^\top x_{l-1}+b_{l}) \,.
\end{equation*}
\noindent Whenever $\sigma_{l}^\prime (W_l^\top x_{l-1}+b_{l}) \neq 0$, this automatically implies 
\[x_{l} = \sigma_{l}(W_{l}^\top x_{l-1}+b_{l}) \, ,\]
ensuring that a critical point satisfies \eqref{eq:nonlinear-constraint}.
Following the MAC-QP approach, minimisation of network parameters can easily be distributed with the right choice of optimisation method. However, minimising the MAC-QP objective with gradient-based fist-order optimisation methods still requires the differentiation of the activation functions. This is different for lifted training approaches that we briefly discuss in the next section.

\subsubsection{Lifted training approach}\label{sec:lifted}

The second line of work to approach Problem \eqref{eq:mlp-constrained} is motivated by the observation that the ReLU activation function $\sigma(z, 0) = \max(z, 0)$ itself can be viewed as the solution to a constrained minimisation problem. More specifically, in the work of \cite{zhang2017convergence}, the authors characterise the activated output of the ReLU activation function $\sigma_l (W_l^\top x_{l-1} + b_{l})$ as% the solution to the following non-smooth convex optimisation problem, or respectively constrained convex optimisation problem with smooth objective:
\begin{align*}
    x_{l} &= \sigma_l (W_l^\top x_{l-1} + b_{l}) \nonumber \\
    &= \argmin_{x \geq 0} \|x - (W_l^\top x_{l-1} + b_{l})\|_2^2 \,\, .
\end{align*}%
This implies that $x_l$ is the vector closest to $W_l^\top x_{l-1} + b_{l}$ that lies in the non-negative orthant. %More general non-decreasing activation functions are shown capable to be represented via constructions of appropriate convex optimisation problems \cite{askari2018lifted}. 
%Hence, instead of writing out the dependency relation of activation variables $\Xs$ explicitly as shown in Equation \eqref{eq:mlp-constrained}, nonlinear constraints can be re-written into conditions similarly to the above as argmin maps. 

In the classical lifted approach \citep{askari2018lifted}, Problem \eqref{eq:mlp-constrained} is consequently approximated with% the following constrained formulation where restraints on activation variables are now directly imposed, i.e.
\begin{equation}\label{eq:mlp-classical-lifted-network-training}
    \min_{\Thetas,\{x_l^{i} \in \mathbb{R}_{+}^{n}\}_{l=1}^{\layer-1}} \frac{1}{s} \sum_{i=1}^{s} \left[ \ell(y^i, f(x^{i}_{\layer-1},\Theta_{\layer})) + \sum_{l=1}^{\layer-1} \|x^i_{l} - f(x^{i}_{l-1},\Theta_{l})\|_2^2 \right] \, ,
\end{equation}
where $\{x_l^{i} \in \mathbb{R}_{+}^{n} \}_{l=1}^{\layer-1}$ denotes that the activation variables across layers $l=1, \ldots, \layer-1$ for each sample $i$ lie in the non-negative orthant, i.e. $\mathbb{R}_{+}^{n} := \{v \in \mathbb{R}^n \,| v_j \geq 0, \; \forall \, j \in \{1 \ldots n\} \}$.  %If ReLU activation is used for example, then the objective is minimised over $\{x_l \geq 0 \}_{l=1}^{\layer-1}$.
%"The core idea of using convex energies to determine hidden unit activations is by observing that ReLU non-linearities (and also other ones such as hard-sigmoid and leaky ReLU) can be stated as proximal operators \xw{add ref}"
This approach is often referred to as the ``lifted" approach because the search space for network parameters is lifted to become a higher dimensional space that now also includes auxiliary variables. However, this ``lifting" is obviously also present in MAC-QP. 

%One arrives at the lifted formulation \eqref{eq:mlp-classical-lifted-network-training} by replacing ReLU activation functions in \eqref{eq:mlp-quadratic-penalty} with non-negativity constraints. 
%constraints on activation variables are relaxed and added to the objective function as penalties. 

%the individual constrained minimisation problems on activation variables are added directly to the learning objective. 
%As an example, training a neural network with affine-linear functions $f(x_{l-1},\Theta_{l}) = W_{l}^\top x_{l-1} + b_{l}$ can now be formulated as the following problem, in which activated variables are enforced to be non-negative.

Even though the overall learning problem \eqref{eq:mlp-classical-lifted-network-training} is not convex, each sub-problem is convex in each individual variable when keeping all other variables fixed. Individual updates are also relatively easy to compute using orthogonal projections onto the non-negative orthant. In terms of distributing the computation of parameters, this approach also enjoys the same benefits that the MAC-QP scheme enjoys.
However, one major limitation that this approach suffers from is that it does not solve the original Problem \eqref{eq:mlp-objective}, but essentially trains an affine-linear network (cf. \cite{zach2019contrastive}). This restriction can be shown by examining the optimality system. Consider the case with $s=1$ and one single training sample pair $(x,y)$, then the optimality condition of \eqref{eq:mlp-classical-lifted-network-training} with respect to parameters $b_l$ reads
%\begin{equation*}
%    (W_{l}^{\top} x_{l-1} + b_{l} - x_{l}) = 0 \, .
%\end{equation*}

%\begin{align*}
%    \nabla_{W_{\layer}} E(W_{\layer};\{x\}) = 0 &\Longleftrightarrow (W_{\layer}^{\top} x_{\layer-1} - x_{\layer})
%    x_{\layer-1}^{\top} \\
%    &\Longleftrightarrow \forall j,j': x_{\layer j} = 0 \vee (W_{\layer}^{\top} x_{\layer-1})_{j'} = x_{\layer j'} ,
%\end{align*}

%for non-vanishing activated inputs, i.e. $x_{l-1} \neq 0$, 
%This suggests that at optimality, the nonlinear activation layers of the learned network would behave linearly:
\begin{equation}\label{eq:optimality_condition_lifted_approach}
   (W_{l}^{\top} x_{l-1} + b_{l} - x_{l}) = 0 \qquad \Longrightarrow \qquad x_{l} = W_{l}^{\top} x_{l-1} + b_{l} \,\,.
\end{equation}%
Despite this major limitation this approach has been considered for computing initialisations of network parameters for other optimisation algorithms \citep{li2019lifted,zach2019contrastive}. 

%Lifted networks introduce explicit variables to represent the activations of network units, and these activations are determined implicitly as minimizers of an underlying optimisation problem, which we call the network energy. 

%Easily . 
%In this section, we examine the optimality system for weight parameters update and point out observations that further supports that our approach is capable to avoid some limitations of the classical lifting network approach.

\section{Lifted Bregman training framework}\label{sec:lifted-bregman-training}

Building upon the previous approaches, we introduce a lifted Bregman framework for the training of feed-forward networks. This framework will fix the aforementioned issues and 1) be capable of recovering a network with nonlinear activation functions, while 2) taking partial derivatives with respect to network parameters will not require computing derivatives of the activation functions. %This allows the unconditional application of gradient-based minimisation algorithms to train these networks even when the activation functions are not differentiable, and further circumvents issues such as vanishing and exploding gradients.

\begin{figure}[!ht]
    \centering
    \includegraphics[width=1\textwidth]{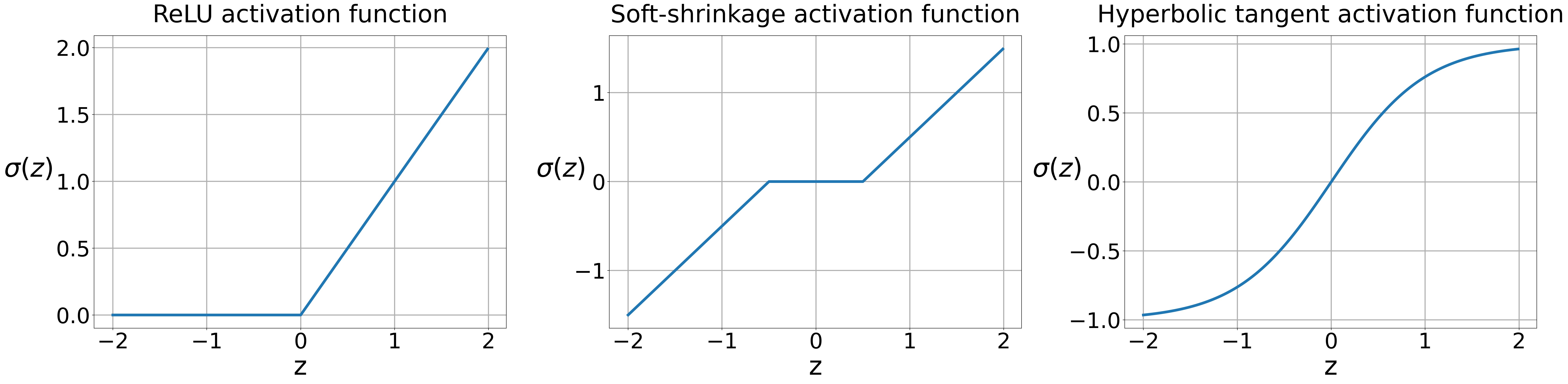}
    \includegraphics[width=0.7\textwidth]{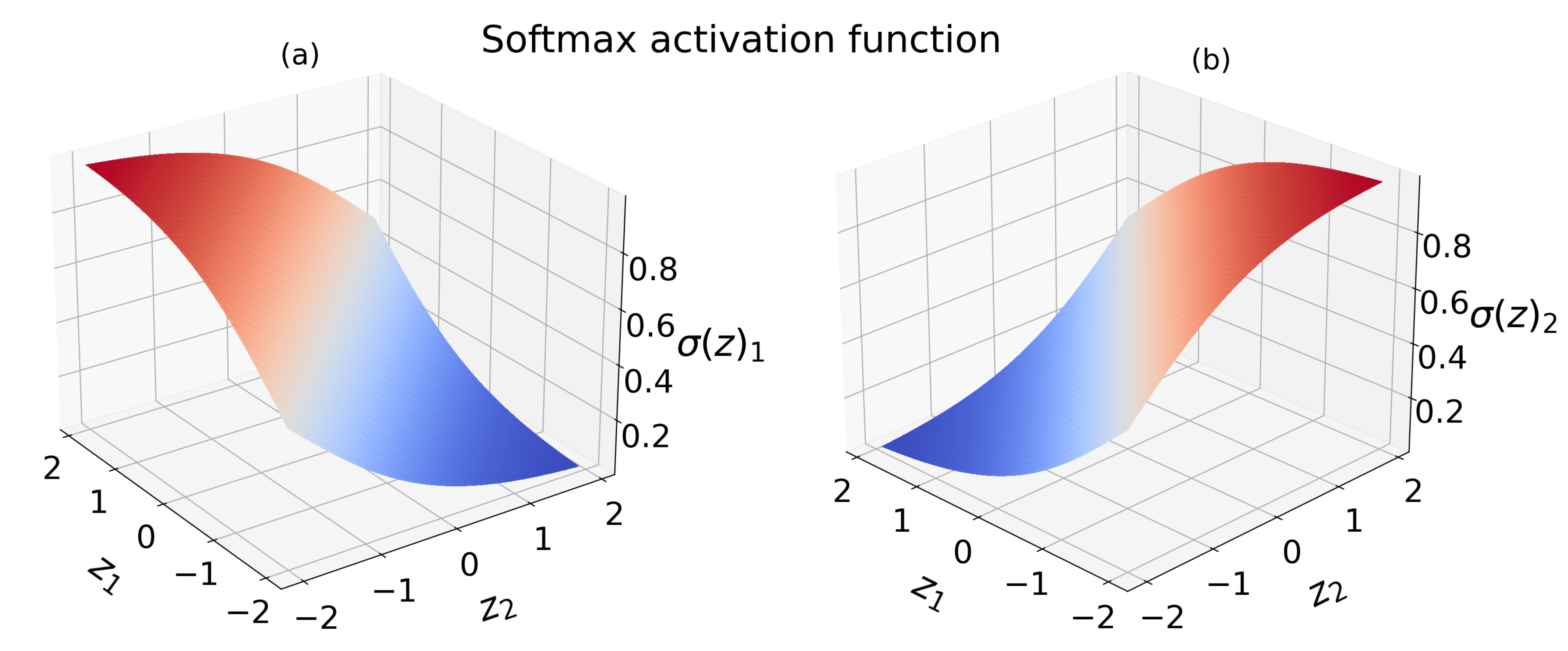}
    \caption{This figure shows example activation functions that are also proximal maps. From left to right on the top row are the ReLU activation function, the soft-thresholding activation function (threshold value $\alpha=0.5$) and the hyperbolic tangent activation function respectively. The bottom two figures visualise components of the softmax activation function for the case of $n=2$.}
    \label{fig:prox_activation_function}  
    \end{figure}

Our approach differs from Problem \eqref{eq:mlp-quadratic-penalty} and Problem \eqref{eq:mlp-classical-lifted-network-training} in the sense that instead of penalising the quadratic Euclidean norm of the nonlinear constraint $x_l^{i} = \sigma_{l} \left(f(x_{l - 1}^{i}, \Theta_l)\right)$ we propose an alternative penalisation. From now on we assume that the activation functions $\sigma_l$ that we consider come from the large class of proximal maps. \begin{definition}[Proximal map]\label{def:proximal-map}
The \emph{proximal map} $\sigma:\mathbb{R}^n \rightarrow \dom(\Psi) \subset \mathbb{R}^n$ of a proper, lower semi-continuous and convex function $\Psi:\mathbb{R}^n \rightarrow \mathbb{R} \cup \{ \infty \}$ is defined as
\begin{align}
    \sigma(z) = \prox(z) := \argmin_{u \in \mathbb{R}^n} \left\{ \frac12 \| u - z \|^2 + \Psi(u) \right\} \, .\label{eq:proximal-map}
\end{align}
\end{definition}
\begin{example}[Proximal maps]\label{ex:prox-maps}
There are numerous examples of commonly used activation functions in deep learning that are also proximal maps (cf. \cite{zhang2017convergence,li2019lifted,combettes2020deep,hasannasab2020parseval}). We give the following four examples. The first one is the \emph{rectifier} or \emph{ramp} function that is also known as Rectified Linear Unit (ReLU) (cf. \cite{nair2010rectified}). The ramp function is a special case of Definition \ref{def:proximal-map} for the characteristic function over the non-negative orthant, i.e.
\begin{align*}
    \Psi(u) := \begin{cases} 0 & u \in [0, \infty)^n \\ \infty & \text{otherwise} \end{cases} \qquad &\implies \qquad \sigma(z)_j = \max( 0, z_j ) \, , \; \forall j \in \{1, \ldots, n \} \, .
\end{align*}
The well-known \emph{soft-thresholding} function is a proximal map for a positive multiple of the one-norm, i.e.
\begin{align*}
    \Psi(u) := \alpha \| u \|_1 \qquad &\implies \qquad  \sigma(z)_j = \begin{cases} z_j - \alpha & z_j > \alpha \\ 0 & |z_j| \leq \alpha \\ z_j + \alpha & z_j < -\alpha \end{cases} \, , \; \forall j \in \{1, \ldots, n \} \, ,
\end{align*}
for $\alpha > 0$. Common smooth activation functions like the \emph{hyperbolic tangent} can also be framed as proximal maps. The hyperbolic tangent can be recovered by choosing the characteristic function
\begin{align*}
    \Psi(u) := \begin{cases} u \tanh^{-1}(u) + \frac12 \left( \log(1 - u^2) - u^2 \right) & |u| < 1 \\ \infty & \text{otherwise} \end{cases} \qquad \implies \qquad \sigma(z) = \tanh(z) \, .
\end{align*}
All previous examples were separable activation functions. Our final example is the non-separable \emph{softmax} activation function, which we can obtain as a proximal map for the characteristic function
\begin{align*}
    \Psi(u) &:= \begin{cases} \sum_{j = 1}^n \left[ u_j \log(u_j)  - \frac{1}{2} u_j^2 \right] & u_j \geq 0 \, , \, \sum_{j = 1}^n u_j = 1 \\ \infty & \text{otherwise} \end{cases} \\ \implies \qquad \sigma(z)_j &= \frac{\exp(z_j)}{\sum_{k = 1}^n \exp(z_k)} \, , \; \forall j \in \{1, \ldots, n \} \, .
\end{align*}%
In Figure \ref{fig:prox_activation_function}, we show visualisations of the above mentioned proximal activation functions.
\end{example}
For convex functions $\Psi$, the following remark gives an alternative characterisation of the proximal map.
\begin{remark}\label{rem:fermat}
Note that an equivalent characterisation of the proximal map $\sigma$ is
\begin{align*}
    z - \sigma(z) \in \partial \Psi(\sigma(z)) \, ,
\end{align*}    
as a consequence of Fermat's theorem. Here $\partial \Psi$ denotes the subdifferential of $\Psi$.
\end{remark}
With the help of Remark \ref{rem:fermat}, we discover that $x_l = \sigma(f(x_{l - 1}, \Theta_l))$ can also be written as
\begin{align}
    f(x_{l - 1}, \Theta_l) - x_l \in \partial \Psi(x_l) \, , \qquad \text{respectively, } \qquad f(x_{l - 1}, \Theta_l) \in \partial \left( \frac12 \| \cdot \|^2 + \Psi \right)(x_l) \, ,\label{eq:dl-inclusion}
\end{align}
where the second inclusion follows from \cite[Chapter 1, Section 5, Proposition 5.6]{ekeland1999convex}. We want to briefly recall the concept of Fenchel duality before we proceed. 
\begin{definition}[Fenchel conjugate]
For a proper, convex and semi-continuous function $F$ its \emph{convex-} or \emph{Fenchel-conjugate} $F^\ast$ is defined as
\begin{align*}
    F^\ast(y) := \sup_{x} \left\{ \langle x, y \rangle - F(x) \right\} \, .
\end{align*}
\end{definition}
The following neat equivalence for inclusion \eqref{eq:dl-inclusion} forms the basis for our new penalisation term. 
\begin{lemma}[Subdifferential characterisation]
Suppose the function $\Psi:\mathbb{R}^n \rightarrow \mathbb{R} \cup \{ \infty \}$ is a proper, lower semi-continuous and convex function. Then the inclusion $z \in \partial \left( \frac12 \| \cdot \|^2 + \Psi\right)(x)$ is equivalent to
\begin{align}
    \frac12 \| x \|^2 + \Psi(x) + \left( \frac12 \| \cdot \|^2 + \Psi\right)^\ast\left( z \right) = \langle x, z \rangle \, .\label{eq:subdifferential-equivalence}
\end{align}
Here, $\left( \frac12 \| \cdot \|^2 + \Psi\right)^\ast$ denotes the Fenchel conjugate of $\frac12 \| \cdot \|^2 + \Psi$.
\end{lemma}
A proof for this lemma can be found in \cite[Theorem 23.5]{rockafellar1970convex}, where a more general result is proven. Now, instead of enforcing \eqref{eq:subdifferential-equivalence} as a hard constraint, we define the penalisation function
\begin{align}
    \bregmanloss(x, z) := \frac12 \| x \|^2 + \Psi(x) + \left( \frac12 \| \cdot \|^2 + \Psi\right)^\ast\left( z \right) - \langle x, z \rangle \, ,\label{eq:bregman-loss}
\end{align}

%Our proposed approach relaxes the non-linear constrains in \eqref{eq:mlp-constrained} and includes the tailored Bregman distances $L_{\Psi}$ in the learning objective as penalties. In contrast to \eqref{eq:mlp-equivalent-activation-penalty}, this results in the following unconstrained optimisation problem:
and replace the squared Euclidean norm in \eqref{eq:mlp-classical-lifted-network-training} with it, which results in the unconstrained optimisation problem
\begin{equation}
    \min_{\Thetas,\Xs}  \sum_{i=1}^{s} \left[\ell\left(y^i, x_{\layer}^{i}\right) + \sum_{l=1}^{\layer-1} \,  \lambda_l \bregmanloss\left(x_{l}^{i}, f(x_{l - 1}^{i}, \Theta_l) \right)\right] \, .
    \label{eq:bregman-lifted}
\end{equation}%
for positive scalars $\{ \lambda_l \}_{l = 1}^{\layer - 1}$. 

%A key advantage of this approach is that at optimality, when computing gradients with respect to $W_{l}$,

We name this penalisation function $\bregmanloss(x, z)$ \emph{Bregman loss}, as we are going to show that it belongs to the family of generalised Bregman distances \citep{bregman1967relaxation,kiwiel1997proximal}. In the following sections, we provide a more detailed analysis of the Bregman loss and show how it can improve the classical lifted training approach and overcome its limitations.

\subsection{The Bregman loss function}

%Our proposed lifted Bregman framework relies on the core idea of the tailored loss function defined in \eqref{eq:bregman-loss}. Instead of imposing hard constraints on auxiliary activation variables as in Problem \eqref{eq:mlp-classical-lifted-network-training}, our approach relaxes the original constrained Problem \eqref{eq:mlp-constrained} via penalising the loss terms as described in \eqref{eq:bregman-lifted}. 
We want to verify that the proposed Bregman loss function \eqref{eq:bregman-loss} is a non-negative generalised Bregman distance \citep{bregman1967relaxation,kiwiel1997proximal} and satisfies some advantageous properties compared to the mean-squared error loss used in \eqref{eq:mlp-classical-lifted-network-training}.

The generalised Bregman distance induced by a proper, convex and lower semi-continuous function $\Phi$ is defined as follows.
\begin{definition}[Generalised Bregman distance]\label{def:bregman-distance}
The generalised Bregman distance of a proper, lower semi-continuous and convex function $\Phi$ is defined as
\begin{align*}
    D_\Phi^{q(v)}(u, v) &:= \Phi(u) - \Phi(v) - \langle q(v), u - v \rangle \, ,
\end{align*}%
Here, $q(v) \in \partial \Phi(v)$ is a subgradient of $\Phi$ at argument $v \in \mathbb{R}^n$.
\end{definition}
For examples of (generalised) Bregman distances we refer the reader to \cite{bregman1967relaxation,censor1981iterative, kiwiel1997proximal,burger2016bregman,benning2021bregman}. 

We now want to verify that the function \eqref{eq:bregman-loss} is a Bregman distance in the sense of Definition \ref{def:bregman-distance}. This directly follows from \eqref{eq:subdifferential-equivalence}, as we can replace $\left( \frac12 \| \cdot \|^2 + \Psi\right)^\ast\left( z \right)$ in \eqref{eq:bregman-loss} with $\langle \sigma(z), z \rangle - \frac12 \| \sigma(z) \|^2 - \Psi(\sigma(z))$ because of $z \in \partial (\frac12 \| \cdot \|^2 + \Psi)(\sigma(z))$. We then obtain
\begin{align*}
    \bregmanloss(x, z) &= \frac12 \| x \|^2 + \Psi(x) + \left( \frac12 \| \cdot \|^2 + \Psi\right)^\ast\left( z \right) - \langle x, z \rangle \, ,\\
    &= \frac12 \| x \|^2 + \Psi(x) - \frac12 \| \sigma(z) \|^2 - \Psi(\sigma(z)) + \langle \sigma(z), z \rangle - \langle x, z \rangle \, ,\\
    &=  \frac12 \| x \|^2 + \Psi(x) - \left( \frac12 \| \sigma(z) \|^2 + \Psi(\sigma(z)) \right) - \langle z, x - \sigma(z) \rangle \, \\
    &= D_{\frac12 \| \cdot \|^2 + \Psi}^{z}(x, \sigma(z)) \, .
\end{align*}
The last equality is correct because of $z \in \partial \left( \frac12 \| \cdot \|^2 + \Psi\right)(\sigma(z))$, which is equivalent to $z - \sigma(z) \in \partial \Psi(\sigma(z))$, respectively $\sigma(z) = \prox(z)$ as pointed out in Remark \ref{rem:fermat}. Since this is true by definition of $\sigma$, we have established that $\bregmanloss(x, z)$ is a Bregman distance. Since Bregman distances are non-negative, this automatically implies non-negativity of $\bregmanloss(x, z)$ for all arguments $x, z$, but we easily get a better lower bound if we split the Bregman distance, i.e.
\begin{align*}
    \bregmanloss(x, z) &= D_{\frac12 \| \cdot \|^2 + \Psi}^{z}(x, \sigma(z)) = D_{\frac12 \| \cdot \|^2}^{\sigma(z)}(x, \sigma(z)) + D_\Psi^{z - \sigma(z)}(x, \sigma(z)) \, ,\\
    &= \frac12 \| x - \sigma(z) \|^2 + D_\Psi^{z - \sigma(z)}(x, \sigma(z)) \geq \frac12 \| x - \sigma(z) \|^2 \, ,
\end{align*}
which we capture with the following corollary.
\begin{corollary}
The loss function $\bregmanloss$ as defined in \eqref{eq:bregman-loss} is a Bregman distance and bounded from below by $\frac12 \| x - \sigma(z) \|^2$.
\begin{proof}
This follows from the fact that $D_\Psi^{z - \sigma(z)}(x, \sigma(z))$ is a valid generalised Bregman distance, because of $z - \sigma(z) \in \partial \Psi(\sigma(z))$ by definition of $\sigma$. Hence, $D_\Psi^{z - \sigma(z)}(x, \sigma(z)) \geq 0$, which concludes the proof.
\end{proof}
\end{corollary}
Hence, we can characterise the exact discrepancy between $\bregmanloss(x, z)$ and $\frac12 \| x - \sigma(z) \|^2$ in terms of $D_\Psi^{z - \sigma(z)}(x, \sigma(z))$, which means we establish
\begin{align}
    \bregmanloss(x, z) &= \frac12 \| x - \sigma(z) \|^2 + \Psi(x) - \Psi(\sigma(z)) - \langle z - \sigma(z), x - \sigma(z) \rangle \, , \nonumber\\
    &= \frac12 \| x - \sigma(z) \|^2 + \Psi(x) + \Psi^\ast(z - \sigma(z)) - \langle x, z - \sigma(z) \rangle \, ,  \label{eq:intermediate-bregman-loss}
\end{align}
%\begin{corollary}
%Suppose $F$ is a proper, convex and semi-continuous function with Fenchel conjugate $F^\ast$. Then we know %that $y \in \partial F(x)$ is equivalent to $x \in \partial F^\ast(y)$. 
%\end{corollary}
We can further simplify this term with the help of the Moreau identity.
\begin{theorem}[Moreau identity \citep{moreau1962fonctions}]\label{thm:moreau}
Suppose $\Psi :\R^n \rightarrow \R \cup \{ \infty \}$ is a proper, convex and lower-semicontinuous function with Fenchel conjugate $\Psi^\ast$, and we have $\sigma(x) := \prox(x)$ and $\sigma^\ast(x) := \prox[\Psi](x)$. Then 
\begin{align*}
    x = \sigma(x) + \sigma^\ast(x)
\end{align*}
holds true for all $x \in \R^n$.
\end{theorem}
Hence, with Theorem \ref{thm:moreau}, we can rewrite \eqref{eq:intermediate-bregman-loss} to
\begin{align*}
    \bregmanloss(x, z) = \frac12 \| x - \sigma(z) \|^2 + \Psi(x) + \Psi^\ast(\sigma^\ast(z)) - \langle x, \sigma^\ast(z) \rangle \, .
\end{align*}
If we define $E_z(x) := \frac12 \| x - z \|^2 + \Psi(x)$ as the function for which we have $\sigma(z) = \prox(z) = \argmin_{x} E_z(x)$ and consequently $0 \in \partial E_z(\sigma(z))$, we can further conclude
\begin{align*}
    \bregmanloss(x, z) &= \frac12 \| x - \sigma(z) \|^2 + \Psi(x) - \Psi(\sigma(z)) - \langle z - \sigma(z), x - \sigma(z) \rangle \, , \\
    &= \frac12 \| x \|^2 + \Psi(x) - \langle x, z \rangle + \frac12 \| \sigma(z) \|^2 - \Psi(\sigma(z)) + \langle z - \sigma(z), \sigma(z) \rangle \, , \\
    &= \frac12 \| x - z \|^2 + \Psi(x) - \frac12 \| z \|^2 - \frac12 \| \sigma(z) \|^2 - \Psi(\sigma(z)) + \langle z, \sigma(z) \rangle\, , \\
    &= \frac12 \| x - z \|^2 + \Psi(x) - \frac12 \| \sigma(z) - z \|^2 - \Psi(\sigma(z)) \, , \\
    &= E_z(x) - E_z(\sigma(x)) = D_{E_z}^0(x, \sigma(x)) \, ,
\end{align*}
which is an alternative characterisation that is a valid Bregman distance because $0 \in \partial E_z(\sigma(z))$ is a valid subgradient of $E_z$ at $\sigma(z)$. This characterisation nicely links to recent work in parametric majorisation for data-driven energy minimisation methods \citep{geiping2019parametric}.

What makes the loss function \eqref{eq:bregman-loss} truly special, however, is its gradient with respect to its second argument. 
\begin{lemma}[Gradient of \eqref{eq:bregman-loss} with respect to second argument]\label{lemma:gradient-bregman-loss}
Suppose $\Psi :\R^n \rightarrow \R \cup \{ \infty \}$ is a proper, convex and lower semi-continuous function with proximal map $\sigma := \prox$. Then the gradient of $\bregmanloss(x, z)$ as defined in \eqref{eq:bregman-loss} with respect to its second argument, i.e. $\nabla_2 \bregmanloss(x, z)$, is
\begin{align*}
    \nabla_2 \bregmanloss(x, z) = \sigma(z) - x \, ,
\end{align*}
for all $x \in \dom(\Psi)$.
\begin{proof}
Note that the term $\frac12 \| x \|^2 + \Psi(x)$ in \eqref{eq:bregman-loss} is constant with respect to the argument $z$ and bounded because of $x \in \dom(\Psi)$, which implies
\begin{align*}
    \nabla_2 \bregmanloss(x, z) = \nabla_z \left[ \left( \frac12 \| \cdot \|^2 + \Psi \right)^\ast(z) - \langle x, z \rangle \right]
\end{align*}
for fixed $x$. The gradient of $\langle x, z \rangle$ with respect to the second argument $z$ is simply $x$; hence, what remains to be shown is that 
\begin{align}
    \nabla \left( \frac12 \| \cdot \|^2 + \Psi \right)^\ast(z) = \sigma(z) \, .\label{eq:fenchel-conjugate-gradient}
\end{align}
Since $\Psi$ is proper, convex and lower semi-continuous, the Fenchel conjugate $(\frac12 \| \cdot \|^2 + \Psi)^\ast(z) = \sup_x \langle x, z \rangle - \frac12 \| x \|^2 - \Psi(x)$ of $\frac12 \| \cdot \|^2 + \Psi$ is well-defined and the supremum over $\langle z, x \rangle - \frac12 \| x \|^2 - \Psi(x)$ is attained. More importantly, we have 
\begin{align*}
    \left(\frac12 \| \cdot \|^2 + \Psi\right)^\ast(z) &= \sup_x \langle x, z \rangle - \frac12 \| x \|^2 - \Psi(x) \, , \\
    &= \frac12 \| z \|^2 + \sup_x - \frac12 \| x - z \|^2 - \Psi(x) \, , \\
    &= \frac12 \| z \|^2 - \inf_x \left\{ \frac12 \| x - z \|^2 + \Psi(x) \right\} \, .
\end{align*}
We know that $\sigma(z) = \prox(z) = \argmin_{x} \left\{ \frac12 \| x - z \|^2 + \Psi(x) \right\}$, which implies $x = \sigma(z)$ and
\begin{align*}
    \left(\frac12 \| \cdot \|^2 + \Psi\right)^\ast(z) &= \frac12 \| z \|^2 - \frac12 \| \sigma(z) - z \|^2 - \Psi(\sigma(z)) \, .
\end{align*}
The quantity $M(z) := \frac12 \| \sigma(z) - z \|^2 + \Psi(\sigma(z))$ is the Moreau-Yosida regularisation of $\Psi$, see \citep{moreau1965proximite,yosida1964functional}, for which we know that its gradient is $\nabla M(z) = z - \prox(z) = z - \sigma(z)$, cf. \cite[Proposition 12.30]{bauschke2011convex}. The gradient of $\frac12 \| z \|^2$ with respect to $z$ is simply $z$, which is why we have verified  \eqref{eq:fenchel-conjugate-gradient}.
\end{proof}
\end{lemma}
\begin{remark}
Note that $\nabla_2 \bregmanloss(x, z)$ may be unbounded from below, which implies that $\nabla_2 \bregmanloss(x, z)$ may never be equal to zero for some choices of $x$. More precisely,  $\nabla_2 \bregmanloss(x, z) = 0$ is equivalent to $z - x \in \partial \Psi(x)$, which requires $\partial \Psi(x) \neq \emptyset$. Hence, if $x$ is chosen such that $\partial \Psi(x) = \emptyset$, we can never guarantee $\nabla_2 \bregmanloss(x, z) = 0$. A simple example is $\Psi(x) = \chi_{\geq 0}(x)$ with negative argument $x$. In that case, there obviously cannot exist an element $z$ such that $\sigma(z) = \max(z, 0) = x$ holds true because $\max(z, 0) \geq 0$ while $x < 0$.
\end{remark}

\noindent We summarise all previous findings in the following theorem.

\begin{theorem}[Bregman loss function]\label{thm:bregman-loss}
The Bregman loss function as defined in \eqref{eq:bregman-loss} satisfies the following properties.
\begin{enumerate}
    \item $\bregmanloss(x, z) = E_z(x) - E_z(\sigma(z))$.
    \item $\bregmanloss(x, z) = D_{\frac12 \| \cdot \|^2 + \Psi}^{z}(x, \sigma(z))$.
    \item $\bregmanloss(x, z) = \frac12 \| x - \sigma(z) \|^2 + D_\Psi^{z - \sigma(z)}(x, \sigma(z)) \geq \frac12 \| x - \sigma(z) \|^2$.
    \item $\bregmanloss(x, z) = \frac12 \| x - \sigma(z) \|^2 + \Psi(x) - \Psi(\sigma(z)) - \langle z - \sigma(z), x - \sigma(z) \rangle$.
    \item For fixed first argument $x$, the function $\bregmanloss$ is continuously differentiable w.r.t the second argument if the first argument satisfies $x \in \dom(\Psi)$.
    \item $\nabla_z \, \bregmanloss(x, z) = \sigma(z) - x$.
    \item The global minimum of $\bregmanloss$ is $x = \sigma(z)$ if $\partial \Psi(x) \neq \emptyset$.
    \item\label{item:bi-convex} $\bregmanloss$ is a bi-convex function, i.e. it is convex w.r.t. $x$ for fixed $z$ and convex w.r.t. $z$ for fixed $x$.
    \item\label{item:monotone} The operator $G_x(z) := \nabla_z \, \bregmanloss(x, z) = \sigma(z) - x$ is monotone.
    \item\label{item:lipschitz} The operator $G_x$ is Lipschitz-continuous with Lipschitz constant one, i.e.
    \begin{align*}
        \| G_x(z_1) - G_x(z_2) \| \leq \| z_1 - z_2 \| \, .
    \end{align*}
    \item\label{item:non-expansive} The operator $G_y$ is firmly non-expansive, i.e.
    \begin{align*}
        \langle G_x(z_1) - G_x(z_2), z_1 - z_2 \rangle \geq \| z_1 - z_2 \|^2 \, .
    \end{align*}
\end{enumerate}
\begin{proof}
The only properties left to prove are Items \ref{item:bi-convex}--\ref{item:non-expansive}. Item \ref{item:bi-convex} follows directly from the definition of $\bregmanloss(x, z)$ in \eqref{eq:bregman-loss}, since both $\frac12 \| \cdot \|^2 + \Psi - \langle \cdot, z \rangle$ and $\left( \frac12 \| \cdot \|^2 + \Psi\right)^\ast - \langle x, \cdot \rangle$ are convex functions for fixed $x$ and $z$, respectively. Item \ref{item:lipschitz} follows directly from the $1$-Lipschitz continuity of the proximal map $\sigma$ (cf. \cite[Theorem 6.42 (b)]{beck2017first}) via $\| G_x(z_1) - G_x(z_2) \| = \| \sigma(z_1) - x - \left( \sigma(z_2) - x \right) \| = \| \sigma(z_1) - \sigma(z_2) \| \leq \| z_1 - z_2 \|$. In similar fashion, Item \ref{item:non-expansive} follows directly from the firm non-expansiveness of the proximal map $\sigma$ (cf. \cite[Theorem 6.42 (a)]{beck2017first}) via $\langle G_x(z_1) - G_x(z_2), z_1 - z_2 \rangle = \langle \sigma(z_1) - x - \left( \sigma(z_2) - x \right), z_1 - z_2 \rangle = \langle \sigma(z_1) - \sigma(z_2), z_1 - z_2\rangle \geq \| z_1 - z_2 \|^2$. The firm non-expansiveness than automatically implies Item \ref{item:monotone}, since $\| z_1 - z_2 \|^2 \geq 0$ for all $z_1, z_2$.
\end{proof}
\end{theorem}

\begin{figure}[t]   
    \centering
    \includegraphics[width=0.83\textwidth]{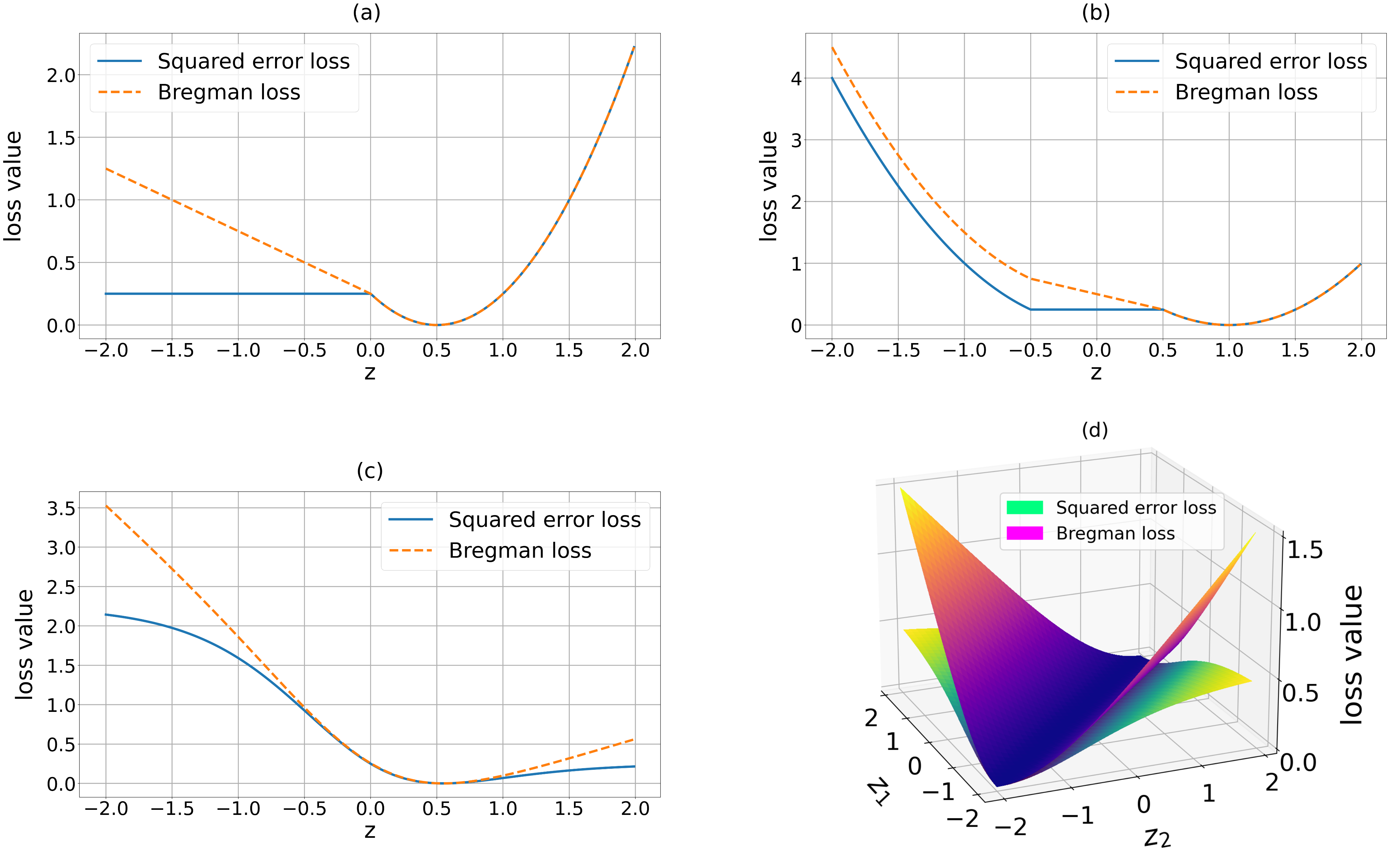}
    \caption{This figure plots the comparison between the Bregman loss function $\bregmanloss(x,z)$ and the squared Euclidean loss $\frac 1 2 \|x-\sigma(z)\|^2$ for $ -2\leq z \leq 2$. We show four cases where each $\sigma$ is: \textbf{a}) the ReLU activation function, \textbf{b}) the soft-thresholding activation function, \textbf{c}) the hyperbolic tangent activation function and (\textbf{d}) the softmax activation function for $n=2$ respectively. The target value $x$ is set to 0.5 when $n=1$ and $x=(0.5,0.5)$ when $n=2$.}
    \label{fig:bregmanloss_vs_squaredloss}
    \end{figure}

\begin{example}
Using the same examples as in Example \ref{ex:prox-maps}, we obtain the following Bregman loss associated with ReLU activation function:
\begin{align*}
    \bregmanloss(x, z) &= \begin{cases} \frac12 \| x - \max(z, 0) \|^2 + \langle x, \max(-z, 0) \rangle & x \in [0, \infty)^n \\ \infty & \text{otherwise} \end{cases} \, .
\end{align*}
%for $\Psi(u) = \begin{cases} 0 & x \in [0, \infty)^n \\ \infty & \text{otherwise} \end{cases}$. %\xw{Add the discussion for the rectifier in terms of the extra penalisation $\langle \max(0, -z), y \rangle$. Also add the functions for the other examples here (you can use Wolframalpha if you like, but try to simplify as much as possible); it might also be nice to add three small individual plots at the top of the page.}
%for $\Psi(u) = \alpha \| u \|_1 $.

Similarly, the Bregman loss associated with the \emph{soft-thresholding} activation function for the scalar case $x, z \in \R$ can be derived as:
\begin{align*}
    \bregmanloss(x, z) = \begin{cases} \frac{1}{2}z^2 -  z\,(\alpha + x)  & z > \alpha \\ -x z  & |z| \leq \alpha \\ \frac{1}{2}z^2 + z\,(\alpha - x) & z < -\alpha \end{cases} \, .
\end{align*}
for $\alpha > 0$.
For the \emph{hyperbolic tangent} activation function, its induced Bregman loss function is computed as
\begin{align*}
    \bregmanloss(x, z) = \begin{cases}  \frac 1 2 \log\left(\frac{1-x^2}{1-\tanh^{2}{z}}\right) + x\,( \tanh^{-1}(x)-z)  & |x| < 1 \\ \infty & \text{otherwise} \end{cases} \, .
\end{align*}

%\begin{align*}
%    \bregmanloss(x, z) = \begin{cases} \log(\cosh{z})  - xz  & |x| < 1 \\ \infty & \text{otherwise} \end{cases} \, .
%\end{align*}
Finally, when considering the non-separable \emph{softmax} activation function, its associated Bregman loss function is
\begin{align*}
    \bregmanloss(x, z) = \begin{cases}  \sum_{j=1}^{n}x_j (\log(x_j) - z_j) + \log( \sum_{j=1}^{n} \exp(z_j)) & x_j \geq 0 \,, \sum_{j = 1}^n x_j = 1 \\ \infty & \text{otherwise} \end{cases} \, ,
\end{align*}
which is basically a shifted version of the multinomial logistic regression loss function.
%\begin{align*}
%    \bregmanloss(x, z) = \begin{cases}  \sum_{j=1}^{n}x_j (\log(x_j) - z_j) - \sum_{j=1}^{n} \frac{\exp(z_j)}{\sum_{k = 1}^n \exp(z_k)}( \log(\frac{\exp(z_j)}{\sum_{k = 1}^n \exp(z_k)}) - z_j) & x_j \geq 0 \,, \\ & \sum_{j = 1}^n x_j = 1 \\ \infty & \text{otherwise} \end{cases} \, .
%\end{align*}
%\begin{align*}
%    \bregmanloss(y, z) = \langle y, \log y \rangle  -  \langle \sigma(z), \log(\sigma(z)) \rangle - \langle y, \sigma(z) \rangle - \langle z - \sigma(z), y - \sigma(z) \rangle \,,
%\end{align*}
%for $u$ such that $ u_j \geq 0 \,, \sum_{j = 1}^n u_j = 1 $ and $\bregmanloss(y, z) = \infty$ otherwise.
In Figure \ref{fig:bregmanloss_vs_squaredloss}, we visualise how each example Bregman loss compares to the squared Euclidean loss.
\end{example}

%\mb{add relevant definitions and properties here}

\subsection{Lifted Bregman training}

By replacing the squared Euclidean penalisations in \eqref{eq:mlp-quadratic-penalty} with the Bregman loss function as defined in \eqref{eq:bregman-loss}, we have generalised \eqref{eq:mlp-quadratic-penalty} to \eqref{eq:bregman-lifted}. To demonstrate the advantage of \eqref{eq:bregman-lifted}, let us consider the case with $s=1$ and one single training sample pair $(x,y)$ again. Using Lemma \ref{lemma:gradient-bregman-loss}, the gradient of \eqref{eq:bregman-lifted} with respect to $\Theta_l$ reads
\begin{align*}
    \nabla_{\Theta_{l}} \bregmanloss(x_{l}, f(x_{l - 1}, \Theta_l)) =  (\sigma(f(x_{l - 1}, \Theta_l)) - x_{l} ) J_{f}^{\Theta}(x_{l - 1}, \Theta_l)\, , %\frac{\partial f}{\partial \Theta_{l}}(x_{l - 1}, \Theta_l)\, ,
\end{align*}
where $J_{f}^{\Theta}$ denotes the Jacobian of $f$ with respect to argument $\Theta$, and for the specific choice $f(x_{l-1},\Theta_{l}) = W_{l}^\top x_{l-1} + b_{l}$ we observe
\begin{align*}
    \sigma(W_{l}^{\top} x_{l-1}^{i} + b_{l}) - x_{l}^{i} = 0 \, ,
\end{align*}
as the optimality condition of \eqref{eq:bregman-lifted} w.r.t. the bias term $b_l$.

This shows that we guarantee network consistency, i.e. $x_{l}^{i} = \sigma(W_{l}^{\top} x_{l-1}^{i} + b_{l})$, for all layers $l$. In other words, the network is truly nonlinear in contrast to the lifted networks described in Section \ref{sec:lifted}, and computing the optimality condition with respect to network parameters does not require differentiation of the activation function as in MAC-QP. %This shows that the Bregman lifted training approach is able to retain the true action of non-linear activations across the neural network, thus layers of a trained Bregman lifted neural network would still behave non-linearly. 

%In the classical lifted network formulation \citep{askari2018lifted}, weights parameters are updated such that the pre-activation $W_{l}^\top x_{l-1}^{i}$ and post-activation $x_{l}^{i}$ would need to match \eqref{eq:optimality_condition_lifted_approach}. Take ReLU activation function as an example, $x_{l}^{i}$ is subject to being non-negative as a hard constraint when solving the minimisation Problem \eqref{eq:mlp-classical-lifted-network-training}. This in turn enforces that $W_{l}^\top x_{l-1}^{i} + b_{l}$ will have to be constructed in a way such that it avoids cutting off negative values. 

%This suggests that the classical lifted approach learns by construction, a linear affine transformation such that the updates of post-activation $x_l^{i}$ is pushed towards the direction of being non-negative. Therefore, searching for an affine transformation under this condition will be restricted, as it has to ensure that the non-negativity of $x_l^{i}$ is achieved by forcing away from the negatives. This indicates that weight parameters $W_l$ in the Bregman lifted approach are updated with larger degree of freedom, thus circumvent some major limitations of the classical lifted approach. 

%This observation holds true in general where similar arguments follow for other proximal map type of activation functions such as soft-threshold activation function or hyperbolic tangent activation function. 

In Appendix \ref{sec:related-works} we discuss other works related to the lifted Bregman training framework and briefly address similarities and differences.

\section{Numerical realisation}\label{sec:numerical-realisation} %with coordinate Bregman methods}
In this section we discuss different computational strategies to solve Problem \eqref{eq:bregman-lifted} numerically. We divide the discussion into deterministic and stochastic strategies.

\subsection{Deterministic approaches}\label{sec:deterministic}
For simplicity of notation, we consider only one pair of samples $(x, y)$ in Problem \eqref{eq:bregman-lifted} and without loss of generality we assume that the loss function $\ell$ is also a Bregman loss function. In order to efficiently solve \eqref{eq:bregman-lifted} we split the objective in \eqref{eq:bregman-lifted} into differentiable and non-differentiable parts as follows: %Using the short-hand notations $\mathbf{\Theta} = \{ \Theta_l \}_{l = 1}^\layer$ and $\mathbf{X} = \{ x_l \}_{l = 1}^{\layer - 1}$, 
\begin{align*}
    E\left(\mathbf{\Theta}, \mathbf{X}\right) {} = {} &\sum_{l=1}^{\layer}  \bregmanloss(x_{l}, f(x_{l - 1}, \Theta_l)) \\
    {} = {} &\underbrace{\sum_{l = 1}^{\layer} \frac12 \| x_l \|^2 + \Psi(x_l)}_{=: G(\mathbf{X})} + \underbrace{\sum_{l = 1}^{\layer} \left[ \left( \frac12 \| \cdot \|^2 + \Psi \right)^\ast\left(f\left(x_{l - 1}, \theta_l\right)\right) - \left\langle x_l, f\left(x_{l - 1},\theta_l\right) \right\rangle \right]}_{=: H(\mathbf{\Theta}, \mathbf{X})} \\
    {} = {} &G(\mathbf{X}) + H(\mathbf{\Theta}, \mathbf{X}) \, .
\end{align*}
Here $H$ is a smooth-function in both arguments, while $G$ is potentially non-smooth but with closed-form proximal map. 

\begin{remark}
We want to emphasise that \eqref{eq:bregman-lifted} can also be split such that the squared Euclidean norm $\frac12 \| x_l \|^2$ of each variable $x_l$ is part of $H$, and not of $G$. Both splittings are perfectly reasonable, and we arbitrarily chose the one that includes $\frac12 \| x_l \|^2$ in $G$ because computing the proximal map of $G$ with added squared Euclidean norm requires only a trivial modification of the proximal map (cf. \cite[Theorem 6.13]{beck2017first}).
\end{remark}

\subsubsection{Proximal gradient descent}\label{sec:prox-grad}

A straight-forward approach for the computational minimisation of $E$ is proximal gradient descent (cf. \cite{teboulle2018simplified}), also known as forward-backward splitting \citep{lions1979splitting}. This method is a special instance of the following Bregman proximal algorithm \citep{censor1992proximal,teboulle1992entropic,chen1993convergence,eckstein1993nonlinear} that aims at finding minimisers of $E$ via the iteration
\begin{align}
    \left( \begin{matrix} \mathbf{\Theta}^{k + 1} \\ \mathbf{X}^{k + 1} \end{matrix} \right) = \argmin_{\mathbf{\Theta}, \mathbf{X}} \left\{ E\left(\mathbf{\Theta}, \mathbf{X}\right) + D_{J}\left( \left( \begin{matrix} \mathbf{\Theta} \\ \mathbf{X} \end{matrix} \right), \left( \begin{matrix} \mathbf{\Theta}^k \\ \mathbf{X}^k \end{matrix} \right) \right) \right\} \, .\label{eq:bregman-proximal}
\end{align}
Here, $D_J$ denotes the Bregman distance with respect to a function $J$, which in case of forward-backward splitting for $E(\mathbf{\Theta}, \mathbf{X}) = G(\mathbf{X}) + H(\mathbf{\Theta}, \mathbf{X})$ reads $J(\mathbf{\Theta}, \mathbf{X}) = \frac{1}{2\tau_{\mathbf{\Theta}}} \sum_{l = 1}^\layer \| \Theta_l \|^2 + \frac{1}{2\tau_{\mathbf{X}}} \sum_{l = 1}^{\layer - 1} \| x_l \|^2 - H(\mathbf{\Theta}, \mathbf{X})$. Here, $\tau_{\mathbf{\Theta}}$ and $\tau_{\mathbf{X}}$ are positive step-size parameters. Then, the optimisation problem $\eqref{eq:bregman-proximal}$ simplifies to
\begin{align*}
    \mathbf{\Theta}^{k + 1} &= \mathbf{\Theta}^k - \tau_{\mathbf{\Theta}} \nabla_{\mathbf{\Theta}} H(\mathbf{\Theta}^k, \mathbf{X}^k) \, ,\\
    \mathbf{X}^{k + 1} &= \argmin_{\mathbf{X}} \left\{ \frac12 \left\| \mathbf{X} - \left( \mathbf{X}^k - \tau_{\mathbf{X}} \nabla_{\mathbf{X}} H(\mathbf{\Theta}^k, \mathbf{X}^k) \right) \right\|^2 + \tau_{\mathbf{X}} G(\mathbf{X}) \right\} \, ,\\
    &= \text{prox}_{\tau_{\mathbf{X}} G}\left( \mathbf{X}^k - \tau_{\mathbf{X}} \nabla_{\mathbf{X}} H(\mathbf{\Theta}^k, \mathbf{X}^k) \right)
\end{align*}
for initial values $\mathbf{\Theta}^0$ and $\mathbf{X}^0$. More precisely, the updates with the previous definition of $G$ and $H$ read 
\begin{subequations}
\begin{align}
    \Theta_l^{k + 1} {} = {} &\Theta_l^k - \tau_{\mathbf{\Theta}} \left( \text{prox}_{\Psi}\left( f(x_{l - 1}^k, \Theta_l^k) \right) - x_l^k \right) J_{f}^{\mathbf{\Theta}}(x_{l - 1}^k, \Theta_l^k) \, , \\
    x_j^{k + 1} {} = {} &\text{prox}_{\frac{\tau_{\mathbf{X}}}{1 + \tau_{\mathbf{X}}} \Psi}\left( \frac{1}{1 + \tau_{\mathbf{X}}} \left( x_j^k - \tau_{\mathbf{X}} \left( \left( \text{prox}_\Psi\left( f(x_j^k, \Theta_{j + 1}^k) \right) - x_{j + 1}^k \right) J_f^x(x_j^k, \Theta_{j + 1}^k) \right. \right. \right. \\
    &\qquad \qquad \qquad \left. \left. \left. - f(x_{j - 1}^k, \Theta_j^k) \vphantom{\frac{1}{1 + \tau_{\mathbf{X}}}} \right) \right) \right) \, ,\nonumber
\end{align}\label{eq:bregman-proximal-network-training-algorithm}%
\end{subequations}
for $l \in \{1, \ldots, \layer\}$ and $j \in \{1, \ldots, \layer - 1\}$, where $J_{f}^{\Theta}$ and $J_{f}^{x}$ denote the Jacobians of $f$ with respect to $\Theta$ and $x$, respectively.

\begin{example}\label{exm:proximal-gd}
Suppose we design a feed-forward network architecture with $\Psi = \chi_{\geq 0}$, which implies $\text{prox}_\Psi(z) = \max(z, 0)$, and $f(x, \Theta_l) = W^\top_l x + b_l$, for $\Theta_l = (W_l, b_l)$. Then \eqref{eq:bregman-proximal-network-training-algorithm} reads
\begin{align*}
    W_l^{k + 1} &= W_l^k - \tau_{W}  \, x_{l - 1}^k \left( \max\left( (x_{l - 1}^k)^\top W_l^k + (b_l^k)^\top, 0\right) - (x_l^k)^\top \right) \, ,\\
    b_l^{k + 1} &= b_l^k - \tau_b \, \left( \max\left( (W_l^k)^\top x_{l - 1}^k + b_l^k \right) - x_l^k, 0 \right) \, ,\\
    x_j^{k + 1} &= \max\left( \frac{ x_j^k - \tau_x \left( W_{j + 1}^k \max\left( (W_{j + 1}^k)^\top x_j^k + b_{j + 1}^k, 0 \right) - x_{j + 1}^k\right) - (W_j^k)^\top  x_{j - 1}^k - b_j^k}{1 + \tau_x}, 0\right) \, ,
\end{align*}
for $l \in \{1, \ldots, \layer\}$, $j \in \{ 1, \ldots, \layer - 1 \}$ (with input $x_0 = x$ and output $x_\layer = y$), $k \in \mathbb{N}$ and the step-size parameters $\tau_{\mathbf{\Theta}} = (\tau_W, \tau_b)$ and $\tau_{\mathbf{X}} = \tau_x$.
\end{example}
Note that many modifications of proximal gradient descent can be applied, such as proximal gradient descent with Nesterov or Heavy-ball acceleration \citep{nesterov1983method,huang2013accelerated,teboulle2018simplified,mukkamala2020convex}. In Appendix \ref{app:numerical-realisation}, we also describe how $E$ can be minimised with alternating minimisation approaches such as coordinate descent and the alternating direction method of multipliers that better exploit the structure of the objective for distributed optimisation.

%\mb{Next coordinate descent, ADMM, linearised Bregman iteration; mention that we can obviously have more, but we have to process the whole batch}

\subsubsection{Regularisation of network parameters}
In empirical risk minimisation, it is common to design regularisation methods that substitute the empirical risk minimisation process, in order to combat ill-conditioning of the empirical risk minimisation and to improve the validation error. In the context of risk minimisation, regularisations are approximate inverses of the model function w.r.t. the model parameters. For more information about regularisation, we refer the interested reader to \cite{engl1996regularization,benning2018modern} for an overview of deterministic regularisation, to \cite{stuart2010inverse} for a Bayesian perspective on regularisation and \cite{de2005learning,de2021regularization} for regularisation in the context of machine learning. We discuss the two common approaches of variational and iterative regularisation and how to incorporate them into the lifted Bregman framework.\\

\noindent \textbf{Variational regularisation.} In variational regularisation, a positive multiple of a regularisation function is added to the empirical risk formulation. Denoting this regularisation function by $R$, we can simply modify $E$ to
\begin{align*}
    E(\mathbf{\Theta}, \mathbf{X}) = \sum_{l = 1}^\layer \bregmanloss\left(x_l, f(x_{l - 1}, \Theta_l)\right) + R(\mathbf{\Theta}) \, .
\end{align*}
However, there is a catch with this modification that is probably not obvious at first glance. Suppose $R$ is differentiable for now. If we compute the optimality condition w.r.t. a particular $\Theta_l$, we observe
\begin{align*}
    0 &= \nabla_{\Theta_l} \left( \bregmanloss(x_l^\ast, f(x_{l - 1}^\ast, \Theta_l^\ast)) + R(\Theta_1^\ast, \ldots, \Theta_l^\ast, \ldots, \Theta^\ast_\layer) \right) \, , \\
    &= \left( \sigma(f(x_{l - 1}^\ast, \Theta_l^\ast)) - x_l^\ast \right) J_{f}^{\Theta}(x_{l - 1}^\ast, \Theta_l^\ast) + \nabla_{\Theta_l} R(\mathbf{\Theta}^\ast) \, .
\end{align*}
Without $\nabla_{\Theta_l} R(\mathbf{\Theta}^\ast)$, the condition $0 = \left( \sigma(f(x_{l - 1}^\ast, \Theta_l^\ast)) - x_l^\ast \right) J_{f}^{\Theta}(x_{l - 1}^\ast, \Theta_l^\ast)$ guarantees network consistency $\sigma(f(x_{l - 1}^\ast, \Theta_l^\ast)) = x_l^\ast$ up to an element in the nullspace of $J_{f}^{\Theta}(x_{l - 1}^\ast, \Theta_l^\ast)$. Unless $\nabla_{\Theta_l} R(\mathbf{\Theta}^\ast) = 0$, this network consistency will be violated when adding a regularisation term $R$ that acts on the network parameters. Because of this shortcoming, we consider iterative regularisation strategies as an alternative in the next subsection.\\

\noindent \textbf{Iterative regularisation.} When we iteratively update the parameters $\mathbf{\Theta}$ and $\mathbf{X}$ via approaches like \eqref{eq:bregman-proximal} or \eqref{eq:bregman-coordinate-descent}, we can convert them into iterative regularisations known as (linearised) Bregman iterations \citep{osher2005iterative,cai2009linearized,benning2018modern}. We achieve this simply by including a regularisation function in the Bregman function. For example, in \eqref{eq:bregman-proximal} we can choose $J(\mathbf{\Theta}, \mathbf{X}) = \frac{1}{\tau_{\mathbf{\Theta}}} \left( R(\mathbf{\Theta}) + \frac{1}{2} \sum_{l = 1}^\layer \| \Theta_l \|^2 \right) + \frac{1}{2\tau_{\mathbf{X}}} \sum_{l = 1}^{\layer - 1} \| x_l \|^2 - H(\mathbf{\Theta}, \mathbf{X})$ instead of $J(\mathbf{\Theta}, \mathbf{X}) = \frac{1}{2\tau_{\mathbf{\Theta}}} \sum_{l = 1}^\layer \| \Theta_l \|^2 + \frac{1}{2\tau_{\mathbf{X}}} \sum_{l = 1}^{\layer - 1} \| x_l \|^2 - H(\mathbf{\Theta}, \mathbf{X})$. Note that the choice of $J$ doesn't affect the optimality conditions of $E$, but it allows to control the regularity of the model parameters $\mathbf{\Theta}$. If the objective $E$ has multiple or even infinitely many minimisers, a different choice of $R$ enables convergence towards network parameters with desired properties such as sparsity of the parameters.

In the context of neural networks, such a strategy has first been applied in \cite{benning2021choose} to train neural networks and also control the rank of the network parameters. We can achieve the same by choosing $R$ to be a positive multiple of the nuclear norm. Recently, in \cite{bungert2021bregman}, the idea of (linearised) Bregman iterations has been extended to stochastic first-order optimisation in order to effectively train neural networks with sparse parameters. Networks with sparse network parameters can be obtained by setting $R$ to a positive multiple of the one-norm. 

\subsubsection{Regularisation of auxiliary variables}\label{sec:regularisation-of-auxiliary-variables}
One of the great advantages of a lifted network approach is that it is straight-forward to also impose regularisation on the auxiliary variables $\{ x_l \}_{l = 1}^{\layer - 1}$. This can be extremely useful in different contexts. Suppose, for example, that we want to train an autoencoder neural network with sparse encoding. This would require us to impose regularity on the output of the encoder network, which is equal to one of the auxiliary variables in the lifted network formulation. We will discuss such a sparse autoencoder approach in greater detail in Section \ref{sec:example-problems-sparse-autoencoder}. In the following, we want to discuss how to adapt the two regularisation strategies discussed in the previous section.\\

\noindent \textbf{Variational regularisation.} In contrast to variational regularisation of network parameters, adding a regularisation function $R$ that acts on $\mathbf{X}$ to the objective $E$ does not impact the optimality system of $E$ w.r.t. the network parameters $\mathbf{\Theta}$. Hence, network consistency will not be violated and regularity can be imposed this way.\\

\noindent \textbf{Iterative regularisation.} In identical fashion to the previous section, we can incorporate regularity by modifying Bregman functions, for example in \eqref{eq:bregman-proximal} via $J(\mathbf{\Theta}, \mathbf{X}) = \frac{1}{2\tau_{\mathbf{\Theta}}} \sum_{l = 1}^\layer \| \Theta_l \|^2 + \frac{1}{\tau_{\mathbf{X}}} \left( R(\mathbf{X}) + \frac{1}{2} \sum_{l = 1}^{\layer - 1} \| x_l \|^2 \right) - H(\mathbf{\Theta}, \mathbf{X})$ instead of $J(\mathbf{\Theta}, \mathbf{X}) = \frac{1}{2\tau_{\mathbf{\Theta}}} \sum_{l = 1}^\layer \| \Theta_l \|^2 + \frac{1}{2\tau_{\mathbf{X}}} \sum_{l = 1}^{\layer - 1} \| x_l \|^2 - H(\mathbf{\Theta}, \mathbf{X})$. %This concludes this discussion on deterministic optimisation strategies.

\subsection{Stochastic approaches}
Having discussed suitable deterministic approaches for the minimisation of $E$, we want to describe how to adopt such approaches in stochastic settings. In particular, we consider the objective $E$ from the previous section, but for $s$ pairs of samples $\{ (x_i, y_i) \}_{i = 1}^s$:
\begin{align}
    E(\mathbf{\Theta}, \mathbf{X}) = \frac1s \sum_{i = 1}^s \sum_{l = 1}^\layer \bregmanloss\left(x_l^i, f(x_{l - 1}^i, \Theta_l) \right) \, .\label{eq:many-data-points-objective}
\end{align}
Here, $\mathbf{\Theta} = \{ \Theta_l \}_{l = 1}^\layer$ is again the short-hand notation for all parameters (for all layers), but $\Xs = \{ x_l^{i} \}_{l=1,\ldots, \layer}^{i=1,\ldots, s}$ is the collection of auxiliary variables that now also depend on the input and output samples $\{ (x_i, y_i) \}_{i = 1}^s$. We want to investigate which stochastic minimisation methods we can formulate that only use a (random) subset of the indices at every iteration.

The most straight-forward approach is to use out-of-the-box, state-of-the-art first order methods like stochastic gradient descent \citep{robbins1951stochastic,kiefer1952stochastic} or variants of it. Since the function $E$ is additively composed of a smooth ($H$) and a non-smooth ($G$) function, the overall function $E$ is non-smooth. Hence, gradient-based first-order methods applied directly to $E$ become subgradient-based first-order methods, with potentially slower convergence and artificial critical points, even though these usually do not pose serious issues with high probability  (cf. \cite{bolte2020mathematical}). However, fully explicit optimisation methods whose convergence speed depends on Lipschitz constants of gradients can easily suffer from an explosion of these constants when applied directly to non-smooth functions. 

Another disadvantage of applying methods such as stochastic gradient descent and evaluating the gradient or subgradient via backpropagation is that it is not straight-forward to easily distribute the computation of parameters as it is the case for distributed optimisation approaches such as the ones described in Section \ref{sec:distributed}, respectively \citep{carreira2014distributed,askari2018lifted,gu2020fenchel}.

And not to forget, we lose the advantage that the differential of the Bregman loss function $\bregmanloss$ does not require the differentiation of the activation function if we try to differentiate the non-smooth part $G$. We therefore discuss three alternative optimisation strategies; we begin with a discussion of stochastic proximal gradient descent, then we continue to focus on data-parallel (implicit) optimisation, before we conclude with implicit stochastic gradient descent.

\subsubsection{Stochastic proximal gradient descent}
Given the structure of $E$, it seems natural to consider stochastic proximal gradient descent \citep{duchi2009efficient,rosasco2019convergence}. However, most approaches, such as the ones discussed in \cite{rosasco2019convergence}, assume a structure of the form
\begin{align*}
    E(\mathbf{\Theta}) = G(\mathbf{\Theta}) + \frac1s \sum_{i = 1}^s H_i(\mathbf{\Theta}) \, ,
\end{align*}
where every $H_i$ is differentiable and $G$ is proximable. In comparison, our setting is of the form
\begin{align}
    E(\mathbf{\Theta}, \mathbf{X}) = \frac1s \sum_{i = 1}^s \left( G(\mathbf{X}_i) + H(\mathbf{\Theta}, \mathbf{X}_i) \right) \, ,\label{eq:all-sample-objective}
\end{align}
where $\mathbf{X}_i$ is the collection $\{ x_1^i, \ldots, x_{\layer - 1}^i \}$ for each index $i \in \{1, \ldots, s\}$. Here, the function $G$ only depends on the samples $\mathbf{X}$ and optimising with respect to $\mathbf{X}$ remains deterministic, while we can optimise $\mathbf{\Theta}$ by only using (random) subsets of $\{1, \ldots, s\}$ at each iteration. A possible approach is to compute
\begin{align*}
    \mathbf{\Theta}^{k + 1} &= \mathbf{\Theta}^k - \tau_{\mathbf{\Theta}}^k \nabla_{\mathbf{\Theta}} H(\mathbf{\Theta}^k, \mathbf{X}_i^k) \, , \\
    \mathbf{X}_i^{k + 1} &=  \text{prox}_{\tau_{\mathbf{X}_i} G}\left( \mathbf{X}_i^k - \tau_{\mathbf{X}_i} \nabla_{\mathbf{X}_i} H(\mathbf{\Theta}^k, \mathbf{X}_i^k) \right) \, ,
\end{align*}
or alternating variants. It should be straight-forward to show convergence for such an algorithm with standard arguments in the simpler but unrealistic setting in which both $H$ and $G$ are jointly-convex in both $\mathbf{X}$ and $\mathbf{\Theta}$. However, proving such a result is beyond the scope of this work.

%This means that most standard approaches are not directly applicable or the theoretical guarantees do not transfer to our setting without a modification of the underlying proofs. Due to other limitations described in the previous section, we do not focus on this approach but rather consider implicit stochastic optimisation approaches.  

\subsubsection{Data-parallel optimisation}\label{sec:data-parallel}
Instead of minimising \eqref{eq:all-sample-objective} for all samples at once, one can split the indices into $m$ randomised batches $B_p$ with $\bigcup_{p = 1}^m B_p = \{1, \ldots, s\}$ and $\bigcap_{p = 1}^m B_p = \emptyset$. We can then solve the optimisation problems for each batch individually and subsequently average all results, which is also widely known as model averaging \citep{zinkevich2010parallelized,mcdonald2010distributed}. The optimisation problem for each batch can be solved in parallel with any of the methods described in Section \ref{sec:deterministic}. In the next section, we focus on an alternative implicit optimisation technique that can be performed sequentially for all batches or in parallel, which in its sequential form is known as implicit stochastic gradient descent.

\subsubsection{Implicit stochastic gradient descent}
Implicit stochastic gradient descent \citep{toulis2017asymptotic} is a straight-forward modification of stochastic gradient descent where each update is implicit. Note that mini-batch stochastic gradient descent for an objective $E(\mathbf{\Theta}) = \frac1s \sum_{i = 1}^s f_i(\mathbf{\Theta})$, which in its usual form reads
\begin{align}
    \mathbf{\Theta}^{k + 1} = \mathbf{\Theta}^k -  \frac{\tau^k}{|B_{p}|} \sum_{i \in B_{p}} \nabla f_i(\mathbf{\Theta}^k) \, ,\label{eq:stochastic-gradient-descent}
\end{align}
can be formulated as
\begin{align}
    \mathbf{\Theta}^{k + 1} = \argmin_{\mathbf{\Theta}} \left\{ \frac{1}{|B_{p}|} \sum_{i \in B_{p}} f_i(\mathbf{\Theta}) + D_{J_k}(\mathbf{\Theta}, \mathbf{\Theta}^k) \right\} \, ,\label{eq:sgd}
\end{align}
see for instance \cite{benning2021bregman}, which shows that it is also a special-case of stochastic mirror descent \citep{nemirovski2009robust}, and where $J_k = \frac{1}{2\tau^k} \| \mathbf{\Theta} \|^2 - \frac{1}{|B_{p}|} \sum_{i \in B_{p}} f_i(\mathbf{\Theta})$, with $D_{J_k}$ being the corresponding Bregman distance. Here, the choice of function $J_k$ ensures the explicitness of mini-batch stochastic gradient descent (cf. \cite[Equation (6)]{benning2021bregman}). We can obviously replace $J_k$ simply with $\frac{1}{2\tau^k} \| \mathbf{\Theta} \|^2$, so that \eqref{eq:sgd} changes to
\begin{align}
    \mathbf{\Theta}^{k + 1} = \argmin_{\mathbf{\Theta}} \left\{ \frac{1}{|B_{p}|} \sum_{i \in B_{p}} f_i(\mathbf{\Theta}) + \frac{1}{2\tau^k} \| \mathbf{\Theta} - \mathbf{\Theta}^k \|^2 \right\} \, .\label{eq:implicit-sgd}
\end{align}
This method, proposed and studied in \cite{toulis2017asymptotic} is similar to classical stochastic gradient descent, but the update can no longer be computed explicitly and requires the implicit solution of a deterministic optimisation problem at each iteration. If we want to connect this to Section \ref{sec:data-parallel}, we can modify \eqref{eq:implicit-sgd} to not use the previous argument $\mathbf{\Theta}$ as the second argument in the Bregman distance, but to use the average of the previous epoch, i.e.
\begin{align}
    \mathbf{\Theta}^{k + 1}_p = \argmin_{\mathbf{\Theta}} \left\{ \frac{1}{|B_p|} \sum_{i \in B_p} f_i(\mathbf{\Theta}) + \frac{1}{2\tau^k_p} \left\| \mathbf{\Theta} - \frac{1}{m} \sum_{q = 1}^m \mathbf{\Theta}^k_q \right\|^2 \right\} \, ,\label{eq:implicit-parallel-sgd}
\end{align}
for $p \in \{1, \ldots, m\}$. The advantage of \eqref{eq:implicit-parallel-sgd} over \eqref{eq:implicit-sgd} is that each batch can be processed in parallel, at the cost of potentially inferior convergence speed. Applying \eqref{eq:implicit-sgd} to objective \eqref{eq:many-data-points-objective} yields the iteration
\begin{align}
    (\mathbf{\Theta}^{k + 1}, \mathbf{X}^{k + 1}) = \argmin_{\mathbf{\Theta}, \mathbf{X}} \left\{ \frac{1}{|B_{p}|} \sum_{i \in B_{p}} \sum_{l = 1}^\layer \bregmanloss(x_l^i, f(x_{l - 1}^i, \Theta_l))) + \frac{1}{2\tau^k} \| \mathbf{\Theta} - \mathbf{\Theta}^k \|^2 \right\} \, .\label{eq:implicit-sgd-for-many-data-points-objective}
\end{align}%
It is important to point out that in \eqref{eq:implicit-sgd-for-many-data-points-objective} not all $x$-variables are updated at once, but only the variables for the current batch $B_{p}$. Hence, an entire epoch is required to update all $x$-variables once. For simplicity, we introduce a short-hand notation $X_{l}^{p} = \{x_l^{i} \}^{i \in B_{p}}$ to denote the collection of $x_l$-variables associated with the current batch $B_{p}$.

As an example, consider the problem of training a feed-forward network by minimising \eqref{eq:bregman-lifted} with additional variational regularisation acting on the auxiliary variable as described in Section \ref{sec:regularisation-of-auxiliary-variables}, for $f(x_{l-1}, \Theta_l) = W^\top_l x_{l-1} + b_l$. Using \eqref{eq:implicit-sgd-for-many-data-points-objective} and \eqref{eq:bregman-proximal-network-training-algorithm} for each individual minimisation problem per mini-batch, the overall algorithm is summarised in Algorithm \ref{alg:implicit-stochastic-lifted-bregman-algorithm}. Here, $K$ refers to the number of epochs, and $N$ refers to the number of iterations of the inner algorithm.

%For each mini-batch we optimise w.r.t the Bregman losses via the proximal gradient descent approach following  \eqref{eq:bregman-proximal-network-training-algorithm}. 

%the following Implicit Batched Lifted Bregman learning algorithm:
%a pseudo Algorithm is summarised in Algorithm \ref{alg:bregman-parallel-algorithm}.

\begin{algorithm}[!t]
\caption{Implicit Stochastic Lifted Bregman Learning}   
\begin{algorithmic}
\label{alg:implicit-stochastic-lifted-bregman-algorithm}
    \STATE Initialise $\Theta_l^{0}$ for each $l \in \{1,2,\ldots,\layer\}$ \;
        \FOR{$k \in\{ 1, 2, \ldots, K\}$}
            \STATE Choose $B_{p} \subset \{1,2,\dots s\}$ either at random or deterministically.\\
            \STATE Initialise $x^{i}_{j} = \sigma_{j} (f(x_{j-1}^{i}, \Theta_{j}^{k}))\text{ for each } i \in B_{p} \text{ and } j \in \{1,2,\ldots,\layer-1\}.  $ \\
            \FOR{$l \in \{\layer, \ldots 1\}$ and $j \in \{\layer-1,\ldots 1\}$ }
            \FOR{$n \in\{0, 1, \ldots, N-1 \}$}
                \STATE Compute step size parameters $\tau_{\Theta}$ and $\tau_{x}$.
                \STATE {\centering $ \displaystyle
                    \begin{aligned}
                    \Theta_l^{n + 1} \leftarrow  \Theta_l^n - \tau_{\Theta} \left( \frac{1}{|B_{p}|} \sum_{i \in B_{p}} \left( \text{prox}_{\Psi} \right. \right. & \left. \left( f(x_{j - 1}^{i,n}, \Theta_l^n) \right) - x_{j}^{i,n}  J_{f}^{\Theta}(x_{j - 1}^{i,n}, \Theta_l^n) \right) \\
                    &\left. \phantom{\sum_{i \in B_{p}}} +\tau^{k} (\Theta_{l}^{n} - \Theta_{l}^{k-1}) \right) \, \end{aligned} $ \par}
                    \FOR{$i \in B_{p}$}
                \STATE \begin{align*}
                x_{j}^{i,n+1} \leftarrow  &\text{prox}_{\frac{\tau_{x}}{|B_{p}| + \tau_{x}} \left(\Psi+R \right) }\left( \frac{|B_{p}|}{|B_{p}| + \tau_{x}} \left( x_j^{i,n} + \frac{\tau_{x}}{|B_{p}|} \left( f(x_{n - 1}^{i,n}, \Theta_j^k) \right. \right. \right. \\
    &\qquad \left. \left. \left. - \left( \text{prox}_\Psi\left( f(x_j^{i,n}, \Theta_{j + 1}^n) \right) - x_{j + 1}^{i,n} \right) J_f^x(x_j^{i,n}, \Theta_{j + 1}^n) \right) \right) \right) \, 
                \end{align*} 
                \ENDFOR
            \ENDFOR
        \STATE $ \Theta_{l}^{k} \leftarrow \Theta_{l}^{N}$
        \ENDFOR
    \ENDFOR
\end{algorithmic}
\end{algorithm}
%\vphantom{\frac{1}{1 + \tau_{\mathbf{X}}}}
Note that similar to the examples described in Section \ref{sec:deterministic}, one can replace the squared-two-norm term in \eqref{eq:implicit-sgd-for-many-data-points-objective} with a Bregman distance term to design explicit-implicit variants of \eqref{eq:implicit-sgd-for-many-data-points-objective}. For example, if we replace $\frac{1}{2\tau^k} \| \mathbf{\Theta} - \mathbf{\Theta} \|^2$ with the Bregman distance $D_{J_k}(\mathbf{\Theta}, \mathbf{\Theta}^k)$ for the choice $J_k(\mathbf{\Theta}) = \frac{1}{\tau^k} \| \mathbf{\Theta} \|^2 - \frac{1}{|B_{p}|} \sum_{i \in B_{p}} \sum_{l = 1}^\layer \bregmanloss(x_l^i, f(x_{l - 1}^i, \Theta_l)))$, we can make \eqref{eq:implicit-sgd-for-many-data-points-objective} explicit with respect to the parameters $\mathbf{\Theta}$, with the potential drawback of more complicated implicit optimisation problems for $\mathbf{X}$.

\section{Example problems}\label{sec:example-problems}
In this section, we present some example problems and demonstrate how the lifted Bregman approach can be utilised to train neural networks. We discuss both supervised and unsupervised learning examples. We consider the supervised learning task of image classification, and the unsupervised learning task of training a sparse autoencoder and a sparse denosing autoencoder.

%We want to emphasis that none of the applications shown are really large scale applications. The idea of this section is rather to demonstrate that the algorithms are applicable to a wide range of different problems, offering the potential to enhance actual large-scale problems. 

%Recall from earlier, ref section and eq, the total training objective will be defined as the following:
%\begin{equation*}
%    \min_{\Thetas,\Xs}  \sum_{i=1}^{s} \left[\ell\left(y^i, x_{\layer}^{i}\right) + \sum_{l=1}^{L-1} \lambda_l \,  \bregmanloss\left(x_{l}^{i}, f(x_{l - 1}^{i}, \Theta_l) \right)+ \alpha_l R_l(\Theta_l) \right] \,.
%\end{equation*} \label{eq:bregman-lifted-regularised}

%where the last layer loss function is chosen depending on the learning task. 

%(last layer use general loss function, then in specific problem, we use bregman loss.

\subsection{Classification}\label{sec:example-problems-classification}
The first example problem that we consider is image classification. In the fully supervised learning setting, pairwise data samples $\{x_i,y_i\}_{i=1}^{s}$ are provided, with each $x_i$ being an input image while $y_i$ represents the corresponding class label. The task is to train a classifier that correctly categorises images by assigning the correct labels. Training an $\layer$-layer network via the lifted Bregman approach can be formulated as the following minimisation problem:
\begin{equation}
    \min_{\Thetas,\Xs}  \frac{1}{s} \sum_{i=1}^{s} \left[\bregmanloss[\Psi_{\layer}]\left(y^i, x_{\layer}^{i}\right) + \sum_{l=1}^{\layer-1} \,  \bregmanloss[\Psi_{l}]\left(x_{l}^{i}, f(x_{l - 1}^{i}, \Theta_l) \right)\right] \,.
\end{equation} \label{eq:mlp-objective_classification}%
Note that we allow for different functions $\Psi_l$ in order to allow for different activation functions for different layers. The loss function $\ell$ in the Learning Problem \eqref{eq:mlp-objective} can also be chosen to be a Bregman loss function $\bregmanloss[\Psi_{\layer}]\left(y^i, x_{\layer}^{i}\right)$, which here measures the discrepancy between the predicted output and the target labels. Note that the special case for $\layer=1$ has already been covered in \cite{wang2020generalised}.

%Here, each function $\Psi_{l}$ for $l=1,\ldots, \layer-1$ determines the choice of activation function $\sigma_{l}$ of the neural network. 
%$\ell$ denotes the Bregman loss function with respect to a stable activation function \xw{add ref proximal network} as defined earlier. For this example, we do not include any additional regularisation further by setting $R=0$. The prediction rule of the trained network is via a forward pass: $ \hat{Y} = W_{L} X_{L-1} + b_{L}$ where $X_{L-1}$ is defined recursively through:\[ X_{\layer} = \sigma_{\layer} (W_{\layer} X_{\layer-1}) \; \text{ for $\layer = 1, 2, \dots, L-1$.}\]

\subsection{Sparse autoencoder}\label{sec:example-problems-sparse-autoencoder}
%In this section, we demonstrate how the lifted Bregman neural network can be used in unsupervised learning scenarios where only unlabelled data is provided. 
%For the next two application examples we consider sparse autoencoders (cf. \cite{lee2007sparse,ranzato2007sparse,ranzato2007unified,glorot2011deep}). 
For the next two application examples we consider sparse autoencoders. Sparse autoencoders aim at transforming signals into sparse latent representations, effectively compressing the original signals. In this section, we formulate an unsupervised, regularised empirical risk minimisation for sparse autoencoders in the spirit of Section \ref{sec:regularisation-of-auxiliary-variables}.

An autoencoder is a composition of two mathematical operators: an encoder and a decoder. The encoder maps input data onto latent variables in the latent space. The latent variables are usually referred to as code. The decoder aims to recover the input data from the code. Unlike regular autoencoders that reduce the dimension of the latent space, the latent space of a sparse autoencoder can have the same or an even larger dimension than the input space. The compression of the input data is achieved by ensuring that only relatively few coefficients of the code are non-zero. The advantage over conventional autoencoders is that the position of the non-zero coefficients can vary for individual signals.

Sparsity on latent nodes can be achieved by explicitly regularising the latent representations \citep{ng2011sparse,xie2012image}. One approach is to penalise the Kullback Leibler divergence between the hidden node sparsity rate and the target sparsity level during training \citep{ng2011sparse, xie2012image}. In this work, we explore an alternative approach using $\ell_1$ regularisation to promote sparsity of the codes in the spirit of compressed sensing (cf. \cite{candes2006robust,donoho2006compressed}), and proceed as discussed in Section \ref{sec:regularisation-of-auxiliary-variables}. We formulate the sparse autoencoder training problem as the minimisation of the energy
\begin{equation}
    \min_{\Thetas,\Xs}  \frac{1}{s} \sum_{i=1}^{s} \left[\bregmanloss[\Psi_{\layer}]\left(x_{0}^{i}, x_{\layer}^{i}\right) + \sum_{l=1}^{\layer-1} \lambda_l \,  \bregmanloss[\Psi_l]\left(x_{l}^{i}, f(x_{l - 1}^{i}, \Theta_l) \right)+ \alpha \|x_{\layer/2}^{i}\|_1  \right] \, ,
\end{equation}\label{eq:mlp-objective_sparse_autoencoder}%
with $\{ x_{0}^{i} \}_{i = 1}^s$ denoting the provided, unlabelled data. Here $\{x_{\layer/2}^{i} \}_{i=1}^s$ are the activation variables that correspond to the code, i.e. the output of the encoder, which we arbitrarily chose at the middle layer $\frac \layer 2$. The regularisation parameter $\alpha \geq 0$ is the hyper-parameter controlling the sparsity of the codes. %Mean Squared Error is used to measure the reconstruction quality.

As discussed in Section \ref{sec:regularisation-of-auxiliary-variables}, one of the advantages of the lifted Bregman framework is that additional regularisation terms acting on the activation variables can easily be incorporated into various algorithmic frameworks. %This helps avoiding to have to differentiate with respect to possible non-smooth regularisation functions, as is usually the case for back-propagation method when computing (sub-)gradients.

\subsection{Sparse denoising autoencoder}\label{sec:example-problems-denoising-sparse-autoencoder}
We present the third example that aims at learning sparse autoencoders for denoising. During training, a denoising autoencoder receives corrupted versions $\tilde{x}$ of the input $x$ and aims to reconstruct the clean signal $x$. By changing the reconstruction objective, the denoising autoencoder is compelled to learn a robust mapping against small random perturbations and must go beyond finding an approximation of the identity function \citep{bengio2013representation, alain2014regularized}. In analogy to \eqref{eq:mlp-objective_sparse_autoencoder}, the proposed learning problem reads
%Sparse denosing autoencoders further combine the advantages of both the sparse coding approach \citep{lee2006efficient, olshausen1996emergence} and ideas of denoising autoencoders \cite{vincent2010stacked} and have been shown effective for denosing noisy images \citep{cho2013simple,xie2012image}. The overall learning problem writes:
\begin{equation}
    \min_{\Thetas,\Xs} \frac{1}{s} \sum_{i=1}^{s} \left[\bregmanloss[\Psi_{\layer}]\left(x_{0}^{i}, \tilde{x}_{\layer}^{i}\right) + \sum_{l=1}^{\layer-1} \lambda_l \,  \bregmanloss[\Psi_l]\left(\tilde{x}_{l}^{i}, f(\tilde{x}_{l - 1}^{i}, \Theta_l) \right)+ \alpha \|\tilde{x}_{\layer/2}^i\|_1  \right] \,,
\end{equation}\label{eq:mlp-objective_sparse_denoising_autoencoder}%
where $\tilde{x}_{0}^{i}$ and $x_{0}^{i}$ are the corrupted and clear input data respectively. Similar to training sparse autoencoders, we regularise the $\ell_1$-norm of the hidden code, i.e. the middle layer activation variable $x_{\layer/2}$ to enforce activation sparsity.

\section{Numerical results}\label{sec:numerical-results}
In this section, we present numerical results for the example applications described in Section \ref{sec:example-problems}. We implement the lifted Bregman framework with Algorithm \ref{alg:implicit-stochastic-lifted-bregman-algorithm} and compare this method to three first-order methods: 1) the classical (sub-)gradient descent method described in \eqref{eq:gradient-descent}, 2) the stochastic (sub-)gradient descent method which follows \eqref{eq:stochastic-gradient-descent}, and 3) the implicit stochastic (sub-)gradient descent that follows \eqref{eq:implicit-sgd}.%, where (sub-)gradients are computed with back-propagation Algorithm \ref{alg:backprop-algorithm}. 

All results have been computed using PyTorch 3.7 on an Intel Xeon CPU E5-2630 v4. \textbf{Code related to this publication will be made available through the University of Cambridge data repository at \url{https://doi.org/10.17863/CAM.86729}}. 

We use the MNIST dataset \citep{lecun1998gradient} and the Fashion-MNIST dataset \citep{xiao2017fashion} for all numerical experiments. For both datasets, we pre-process the data by centring the images and by converting the labels to one-hot vector encodings.

\subsection{Classification}
For the classification task, we follow the work of \cite{zach2019contrastive} and consider a fully connected network with $\layer = 4$ layers with ReLU activation functions. More specifically, we use $f(x_{l-1},\Theta_{l}) = W_{l}^\top x_{l-1} + b_{l}$ with $W_l \in \R^{m_{l-1} \times m_{l}}$ and $b_{l} \in \R^{m_{l} \times 1}$ where $m_{1} = 784$, $m_{2} = m_{3} = 64$ and $m_{4} = 10$. In solving the classification problem described in Section \ref{sec:example-problems-classification}, we do not include any additional regularisation term and set $R\equiv0$ in the learning objective \eqref{eq:mlp-objective_classification}. We use the Mean Squared Error (MSE) between the prediction and the one-hot encoded target values as our classification loss, which also is a Bregman loss function \eqref{eq:bregman-loss} for $\Psi_\layer \equiv 0$. %In the lifted Bregman training scheme, setting $\Psi_{\layer}=0$ naturally corresponds to using an identity map as last layer's activation function and hence also recovers the MSE loss function \citep{benning2021bregman, wang2020generalised}. 
%\eqref{eq:implicit-sgd-for-many-data-points-objective}, 
%follow the lifted Bregman learning scheme to 

We apply Algorithm \ref{alg:implicit-stochastic-lifted-bregman-algorithm} to train the classification network with the lifted Bregman framework. In solving each mini-batch sub-problem \eqref{eq:implicit-sgd-for-many-data-points-objective}, we run a maximum of $N = 15$ iterations with $\tau^{k} = 100$. The other step-size parameters are set to $\tau_{W_l} = 1.99/\left( \|X_{l}^{p}\|_2^2 + \tau_k/2 \right)$ and $\tau_{x} = 1.99/\|W_{l}\|_2^2$ to ensure convexity of $J$ in \eqref{eq:bregman-proximal}. For comparison, the same network architecture is trained via the stochastic (sub-)gradient method \eqref{eq:stochastic-gradient-descent} in combination with the back-propagation Algorithm \ref{alg:backprop-algorithm}, which we will refer to as the SGD-BP approach. Out of the learning rates $\{1 \times 10^{-5}, 5 \times 10^{-5}, 1 \times 10^{-4}, 5 \times 10^{-4}, 1 \times 10^{-3}, 5 \times 10^{-3}, 1 \times 10^{-2},5 \times 10^{-2}, 1 \times 10^{-1}, 5 \times 10^{-1} \}$, we found that $1 \times 10^{-3}$ works best empirically in terms of receiving lowest training objective values. %$\tau_b = 1.99$. The same batch size $|B_k| = 100$ is used. 

For both the MNIST and Fashion-MNIST datasets, we use 60,000 images for training and 10,000 images for validation. Network parameters in all experiments are identically initialised following \cite{glorot2010understanding}. We choose batch size $|B_k| = 100$ and train the network for 100 epochs. 

Table \ref{classification_accuracy_table} summarises the achieved training and validation classification accuracy after training. For the classical lifted training approach, we quote results from \cite{zach2019contrastive}. The lifted Bregman scheme achieves comparable classification accuracy to the SGD-BP approach and shows substantial improvements over the classical lifted approach that produces affine-linear networks.

%\begin{table}

%Standard Lifted Network & 85.6\% & 86.3\%   \\
%Bregman Lifted Network &  99.3\% & 96.8\% \\
%SGD with back-prop &  99.6\% & 96.9\% \\

%Standard Lifted Network & 81.4\% & 80.0\%   \\
%Bregman Lifted Network &  93.5\% & 85.7\% \\
%SGD with Back-prop &  95.8\% & 87.9\% \\

\begin{table}[!ht]
\begin{center}
\begin{tabular}{ | c | c | c | c | c |}
\hline

& \multicolumn{2}{|c|}{MNIST} & \multicolumn{2}{|c|}{Fashion-MNIST} \\  
\hline
Model & Train & Test & Train & Test \\  
\hline
Standard Lifted Network & 85.6\% & 86.3\%  & 81.4\% & 80.0\%  \\
\hline
Lifted Bregman Network &  99.3\% & 96.8\%  &  93.5\% & 85.7\%  \\
\hline
SGD-BP &  99.6\% & 96.9\% &  95.8\% & 87.9\% \\
\hline
\end{tabular}
\end{center}
\vspace{1.5mm}
\caption{This table summarises training and test classification accuracies for the MNIST and Fashion-MNIST datasets, where the network is trained with the lifted Bregman approach, the classical lifted training approach and the stochastic (sub-)gradient method with back-propagation.}
\label{classification_accuracy_table}
\end{table}

In Table \ref{linear_activation_percentage_table} we present the evaluated percentage of linear activations of the network under different training strategies. This rate computes the percentage of number of nodes in each hidden layer that perform linearly, i.e. ($z^i = \max(z^i,0)$). As described in Section \ref{sec:lifted}, classical lifted network training produces affine-linear networks. We empirically verify that the lifted Bregman framework overcomes this limitation.

\begin{table}[!ht]
\begin{center}
\begin{tabular}{ |c | c | c | c | c | c | c |}
\hline
& \multicolumn{3}{|c|}{MNIST} & \multicolumn{3}{|c|}{Fashion-MNIST} \\ 
\hline
Model & Layer 1 & Layer 2 & Layer 3  & Layer 1 & Layer 2 & Layer 3 \\
\hline
Standard Lifted Network & 99.9\%  &99.9\%  & 99.9\%  & 99.9\%  &99.9\%  & 99.9\% \\
\hline
Lifted Bregman Network &  47.4\% & 82.2\% & 70.1\%  &  59.2\% & 85.7\% & 64.7\%  \\
\hline
SGD-BP &  40.1\% & 29.7\% & 73.6\% &  32.2\% & 28.1\% & 78.4\%   \\
\hline
\end{tabular}
\end{center}
\vspace{1.5mm}
\caption{In this table, we show the percentage of nodes in the network's hidden layers that act linearly. The network is trained with the lifted Bregman approach, the classical lifted approach and the the stochastic (sub-)gradient method with back-propagation respectively, on both the MNIST and Fashion-MNIST datasets.}
\label{linear_activation_percentage_table}
\end{table}

%Standard Lifted Network & 99.9\%  &99.9\%  & 99.9\%  \\
%Bregman Lifted Network &  47.4\% & 82.2\% & 70.1\%  \\
%SGD with back-prop &  40.1\% & 29.7\% & 73.6\%  \\
%Standard Lifted Network & 99.9\%  &99.9\%  & 99.9\%  \\
%Bregman Lifted Network &  59.2\% & 85.7\% & 64.7\%  \\
%SGD with back-prop &  32.2\% & 28.1\% & 78.4\%  \\

\subsection{Sparse Autoencoder}\label{sec:numerical-results-sparse-autoencoder}

As discussed in Section \ref{sec:example-problems-sparse-autoencoder}, we train a fully-connected autoencoder with $\layer = 4$ layers with $f(x_{l-1},\Theta_{l}) = W_{l}^\top x_{l-1} + b_{l}$. The hidden dimension is set to $m_{l} = 784$ for $l=1, \ldots,4$, i.e. there is no explicit reduction in dimension. We use the MSE reconstruction loss and the regularisation parameter $\alpha$ is chosen at $\alpha = 0.09$ and does not ensure optimal validation errors, but a sparsity rate of the code of approximately 90\% to guarantee an implicit reduction in dimension instead. 
%$\alpha=0.09$ 
%As discussed in \ref{sec:regularisation-of-auxiliary-variables, sec:example-problems-sparse-autoencoder}. 

Note that the additional regularity on the code can easily be implemented in the lifted Bregman framework provided we choose a suitable activation function for the activation variable that corresponds to the code. Specifically, when the activation variable $x_{2}$ is activated by the soft-thresholding function, the additional non-smooth $\ell_1$-norm regularisation term in \eqref{eq:mlp-objective_sparse_autoencoder} can easily be incorporated by modifying the $\Psi_2$ function in the update step of the $x_2$ variable in Algorithm \ref{alg:implicit-stochastic-lifted-bregman-algorithm}. In this example, we apply the soft-thresholding activation function after the second affine-linear transformation and adopt ReLU activation functions for all other layers. 

%$\alpha = 0.05$
%the additional non-smooth $\ell_1$-norm regularisation term easily handled by a slight modification of the proximal step update of the hidden $X$-variable in Algorithm \ref{alg:bregman-parallel-algorithm}. 

%Specifically, when the latent activation variable $X_{\layer/2}$ is activated by the soft-thresholding function, the additional non-smooth $\ell_1$-norm regularisation term can be easily handled by a slight modification of the proximal step update of the hidden $X$-variable in Algorithm \ref{alg:bregman-parallel-algorithm}. 
\textbf{MNIST-1K} For this application example, we consider a slightly more challenging learning scenario (referred to as MNIST-1K). We limit the amount of training data, and compare how well the trained network generalises on a larger validation set. We choose $s = 1,000$ images from the MNIST dataset at random and use it as our training dataset and use all $10,000$ images from the validation dataset for validation.
%and $\tau_b = 1.9$. 

%an additional $\ell_1$-norm regularisation term. 
%\label{eq:mlp-objective_sparse_denoising_autoencoder}
%he autoencoder parameters are approximated by iteratively 

We experiment with both deterministic and stochastic implementations of the lifted Bregman framework. For the stochastic implementation (referred to as LBN-S), network parameters are trained via Algorithm  \ref{alg:implicit-stochastic-lifted-bregman-algorithm}. We choose batch size $|B_k|=20$ and step-size $\tau_k = 1$. In each batch sub-problem, the minimisation problems are iterated for $N = 15$ iterations, with step-sizes computed via $\tau_{w_l} = 1.99s/(\|X^{p}_l\|_2^2+\tau_k/2)$ and $\tau_b = 1.99$. The deterministic implementation (LBN-D) follows the same update steps described in Algorithm \ref{alg:implicit-stochastic-lifted-bregman-algorithm} with $\tau^{k}=0$ and $|B_k| = s$, i.e. we use all data at every epoch. The step size $\tau_{w_l}$ is computed as $ 1.99s/\|X_l\|_2^2$. 

For comparison, we train the autoencoder also with two first-order methods: (sub-)gradient descent (GD-BP), which follows \eqref{eq:gradient-descent}, and the stochastic (sub-)gradient descent approach (SGD-BP) described in \eqref{eq:stochastic-gradient-descent}, both in combination with the back-propagation Algorithm \eqref{alg:backprop-algorithm} for the computation of (sub-)gradients. 

As learning objective we choose the MSE loss plus $\alpha$ times the $\ell_1$-norm regularisation. Various learning rates in the range $[10^{-4},10^{-1}]$ have been tested. In terms of minimising the training objective values, we find that $4 \times 10^{-2}$ and $3 \times 10^{-2}$ work optimal for GD-BP and SGD-BP approaches, respectively.
%\begin{equation*}
%\mathbf{\Theta}^{k+1} = \frac{1}{s} \sum_{i = 1}^{s} \left[ \ell(y^i,\mathcal{N}(x^i,\Thetas)) + \alpha \|X_{\layer/2}\|_1  \right]\,.
%\end{equation*}%

\begin{figure}[!t]
    \centering
    %original scale was scale=0.285
    \includegraphics[width=0.4\textwidth]{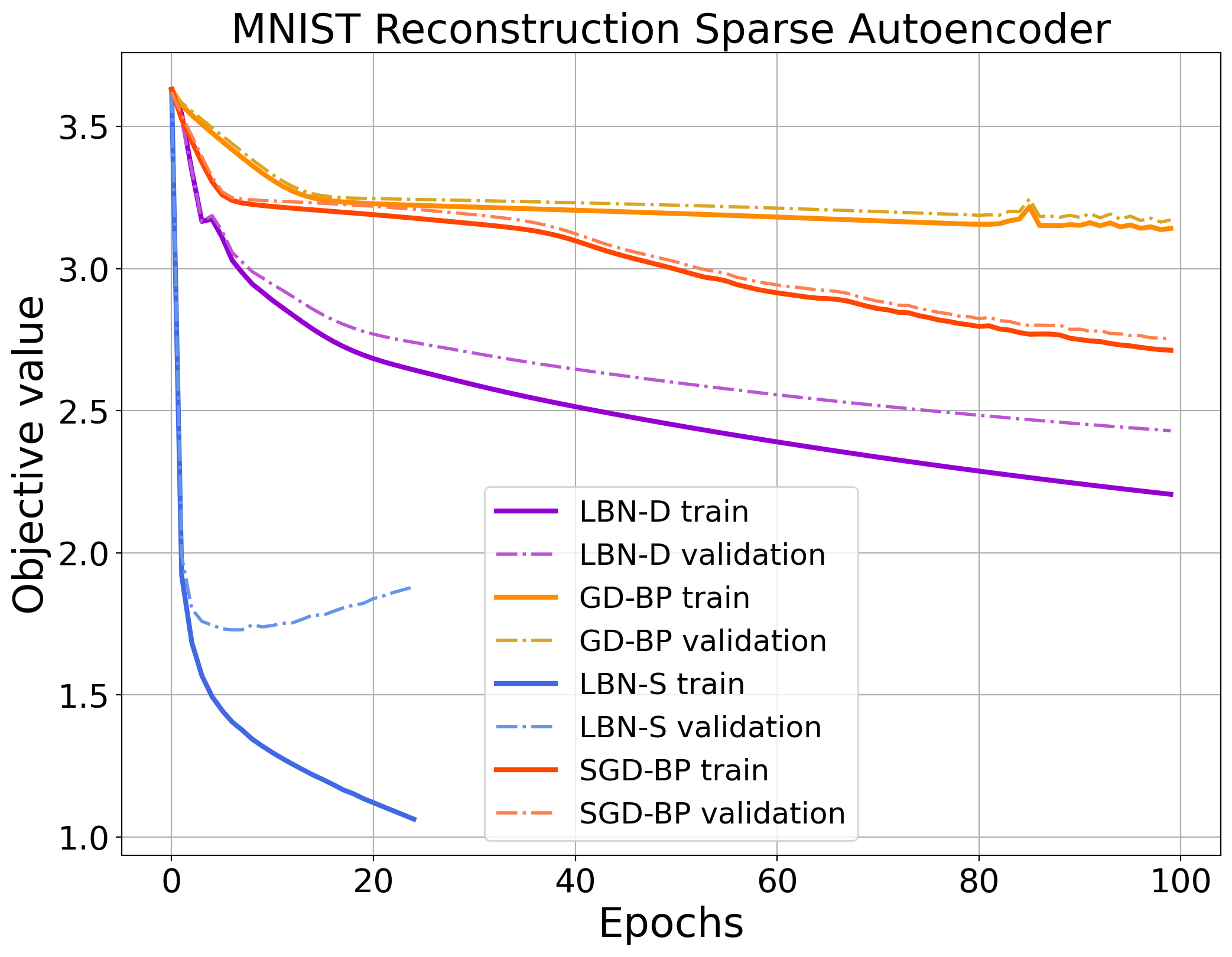}\hspace{0.5cm}
    \includegraphics[width=0.4\textwidth]{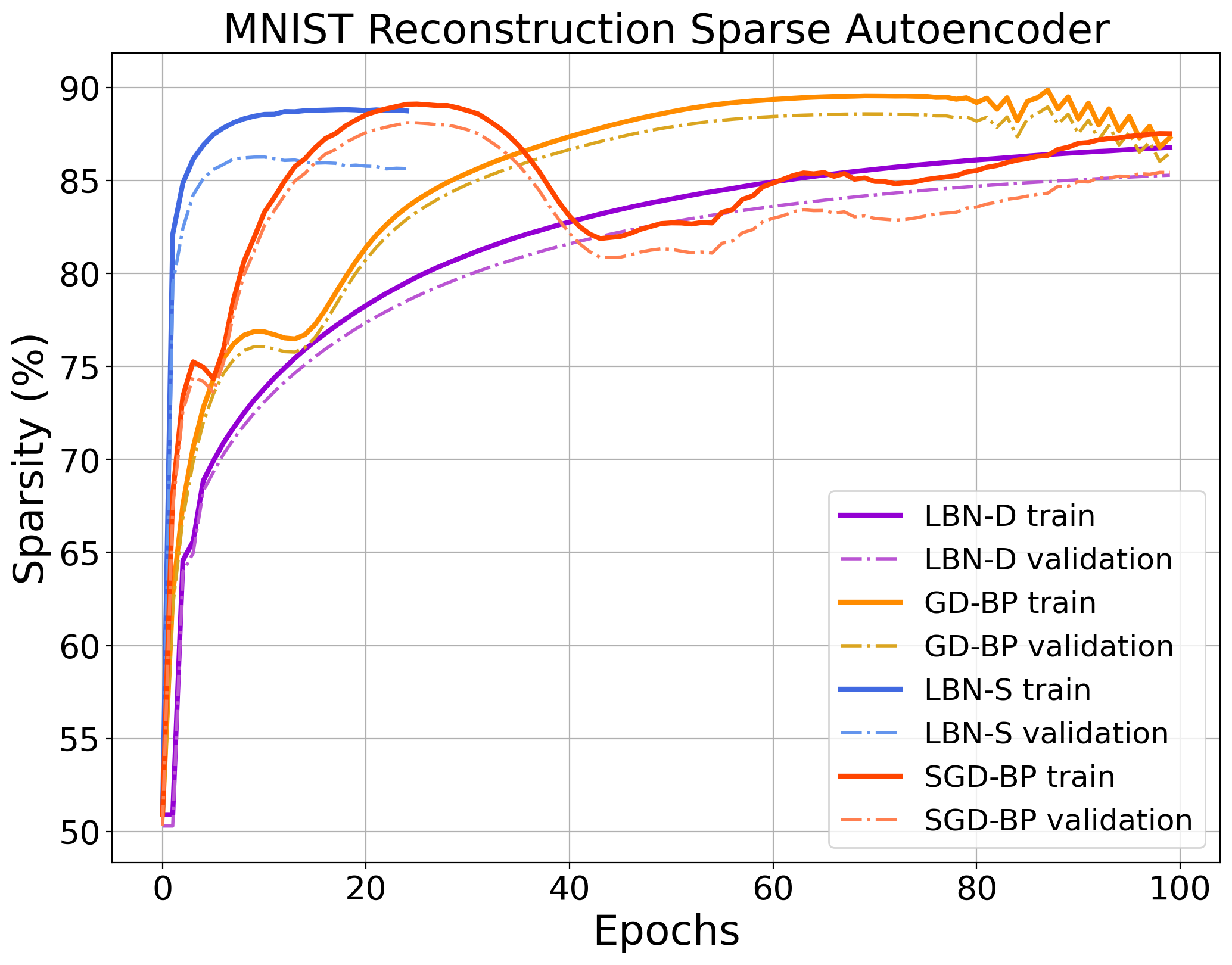}
    \caption{Sparse Autoencoder trained on MNIST-1K dataset. \textbf{Left:} Objective values per epoch for each learning scheme as explained in Section \ref{sec:numerical-results-sparse-autoencoder}. The training objectives (solid lines) record MSE loss plus $\ell_1$-norm regularisation and the validation objectives (dashed lines) report the MSE loss values. \textbf{Right:} Sparsity level of the code per training epoch for each learning scheme.}
    \label{figure:sparse-autoencoder-mnist-objective-sparsity}
\end{figure}

%Network parameters are initialised with the same He initialisation \cite{he2015delving} and 
We train the autoencoder network for 100 epochs and plot the decay of the objective values over all epochs in Figure  \ref{figure:sparse-autoencoder-mnist-objective-sparsity} (Left). While the training objective value tracks the sum of reconstruction loss and $\alpha$ times $\ell_1$-norm regularisation, the validation objective only records the MSE loss. The sparsity rate per epoch is computed as the percentage of zero-valued entries of the code $X_2$, which we visualise in Figure \ref{figure:sparse-autoencoder-mnist-objective-sparsity} (Right). 
Overall, the lifted Bregman implementations outperform the (sub-)gradient-based first-order methods, although it should be emphasised that LBN-D and LBN-S are computationally more expensive than GD-BP and SGD-BP (and we will present a fairer comparison in the denoising autoencoder section). We also observe that the stochastic variants of both learning schemes can speed up training and provide faster convergence compared to their deterministic counterparts. The effect on LBN-D is more noticeable than on GD-BP as we observe significant drops in objective value in early training epochs.

The validation loss curves recovered with GD-BP and SGD-BP follow their training curves quite closely. LBN-S on the other hand shows a drop of the validation loss before it increases again after around 7 epochs. A similar increase could be observed for LBN-D if we were to train the parameters for more epochs. This demonstrates that we can overfit this autoencoder with 3136 parameters to the 1000 training data samples of dimension 784 with the LBN algorithms, but we do not observe the same phenomenon with the GD-BP and SGD-BP algorithms. In light of \cite{zhang2021understanding}, this raises the question whether good generalisation properties of neural networks are more down to the choice of optimisation technique than architectural design choices. 
%This suggests that an overfitting phenomenon could occur when fitting to small training dataset with a lifting training scheme. One possible explanation could be that the lifting scheme increases the search parameter space dimensions which allows for larger model complexity. Lifting scheme is hence capable of solving the parameter estimation problem w.r.t training data more precisely therefore causing an overfit. Nevertheless this potential issue can be effectively mitigated using early-stopping during training.

When looking at a sparsity rate comparison per epoch, we notice the oscillatory behaviours of both (sub-)gradient based methods. As their sparsity levels climb higher, their corresponding objective value curves remain stagnant. The decreases in objective values for both GD-BP and SGD-BP roughly coincide with when their sparsity rate drops. In contrast, we can see that the lifted Bregman implementations handle the same level of regularisation more smoothly. They achieve similar sparsity levels in the end but in a more stable and controlled manner.

\begin{figure}[!t]
    \centering
    \begin{minipage}{0.49\textwidth}
        \centering
        \includegraphics[width=\textwidth]{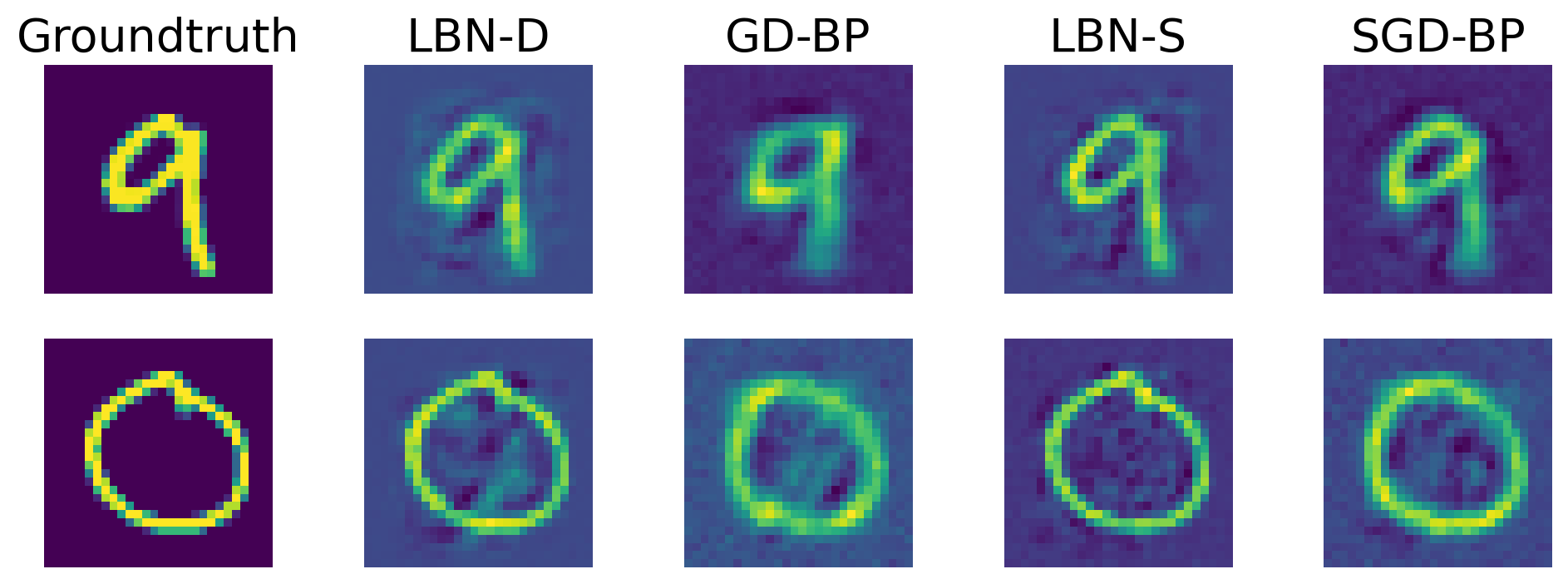}
        \caption{Reconstructed images from the MNIST training dataset for the LBN-D, LBN-S, GD-BP and SGD-BP training strategies respectively, along with ground truth images.}
        \label{figure:sparse-autoencoder-mnist-train-reconstruction}
    \end{minipage}\hfill
    \begin{minipage}{0.49\textwidth}
        \centering
        \includegraphics[width=\textwidth]{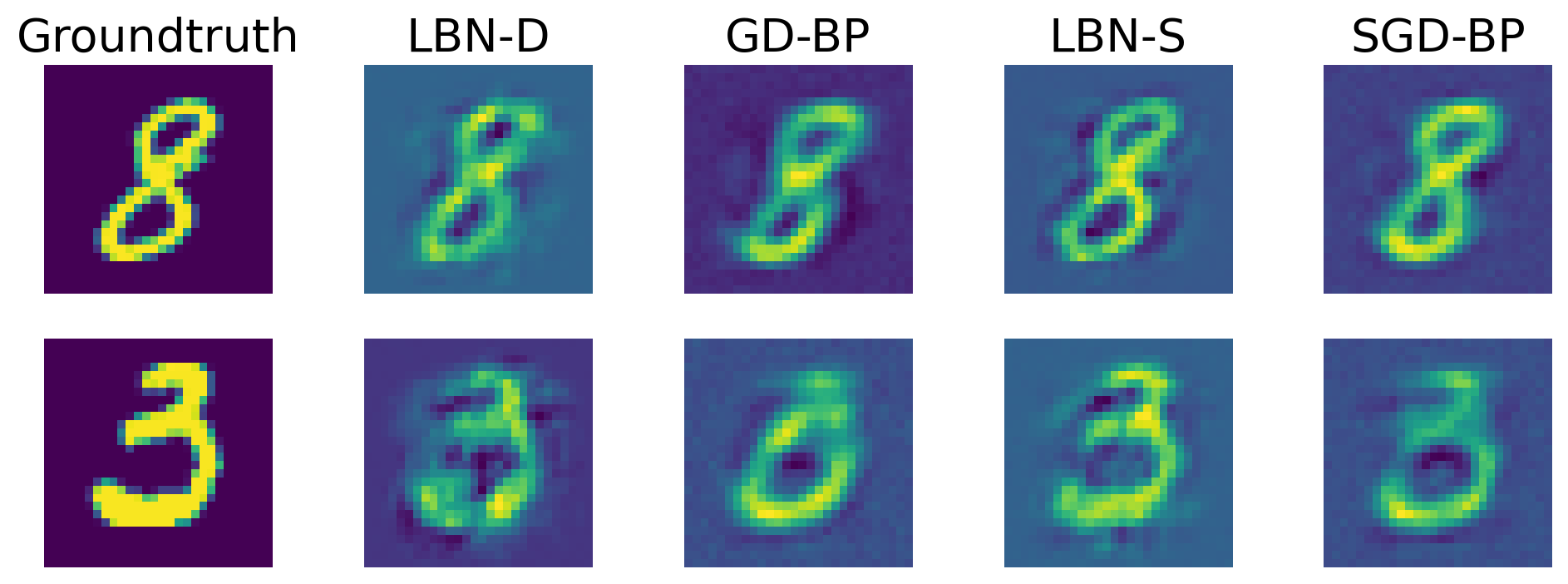}
        \caption{Reconstructed images from the MNIST validation dataset for the LBN-D, LBN-S, GD-BP and SGD-BP training strategies respectively, along with ground truth images.}
        \label{figure:sparse-autoencoder-mnist-val-reconstruction}
    \end{minipage}
\end{figure}

In Figure \ref{figure:sparse-autoencoder-mnist-train-reconstruction}, we visualise ground truth images and reconstructions of two randomly selected images from the training dataset for all four training approaches. In Figure \ref{figure:sparse-autoencoder-mnist-val-reconstruction} we visualise the same quantities but based on two randomly selected images from the validation dataset. All four approaches are capable of providing good quality reconstructions in light of the high sparsity levels, for both the training and testing datasets. However, the proposed LBN-D and LBN-S approaches provide sharper edges and finer details than the GD-BP and SGD-BP approaches. LBN-S slightly outperforms LBN-D and is capable of defining even clearer and sharper edges. More reconstructed images can be found in Appendix \ref{sec:additional-visualisation}. 
% \begin{figure}[!ht]
%     \centering
%     \includegraphics[scale=0.45]{figures/Numerical_Recons_SAE_MNIST_Recons_2_Train_Samples.png}
%     \caption{Images from the MNIST training dataset reconstructed with the sparse autoencoder described in this section, computed with different training strategies. Ground truth images are visualised in the left-most column.}
%     \label{figure:sparse-autoencoder-mnist-train-reconstruction}
% \end{figure}
% \begin{figure}[!ht]
%     \centering
%     \includegraphics[scale=0.45]{figures/Numerical_Recons_SAE_MNIST_Recons_2_Val_Samples.png}
%     \caption{Reconstructed images from the MNIST validation dataset for the LBN-D, LBN-S, GD-BP and SGD-BP training strategies respectively, along with ground truth images.}
%     \label{figure:sparse-autoencoder-mnist-val-reconstruction}
% \end{figure}

\subsection{Sparse Denoising Autoencoder}\label{sec:numerical-results-denoising-sparse-autoencoder}
%As described in Section \ref{sec:numerical-results-sparse-autoencoder}, 

The sparse denoising autoencoder adopts the same network architecture as the sparse autoencoder described in Section \ref{sec:numerical-results-sparse-autoencoder}, which sets the hidden dimension to 784 across all layers. To produce noisy images, instances of Gaussian random variables with mean zero and standard deviation $10^{-3}$ are added to each pixel. We use the MSE as the last layer loss function to measure the difference of reconstruction and noisy images. We consider two training scenarios to evaluate model performances: 1) In the first scenario, we train the network with limited data from the Fashion-MNIST dataset, similar to Section \ref{sec:numerical-results-sparse-autoencoder}, where we take 1,000 training images (referred to as Fashion-MNIST-1K) and validate on 10,000 images. 2) For the second scenario, the training dataset consists of 10,000 images (referred to as Fashion-MNIST-10K) and the validation dataset consists of 10,000 images.\\ %In Figure \ref{Fashion-MNIST-sample-noisy-images} we visualise five clean and noisy sample images from the training dataset.
% \begin{figure}[!ht]
%     \centering
%     \includegraphics[scale=0.5]{fmnist_dataset.png}
%     \includegraphics[scale=0.5]{fmnist_noisy_dataset.png}
%     \caption{\textbf{Top row:} Clean sample images from the Fashion-MNIST dataset. \textbf{Bottom row:} The same images corrupted by additive Gaussian noise.}
%     \label{Fashion-MNIST-sample-noisy-images}
% \end{figure}\\
\textbf{Fashion-MNIST-1K} 
In the first group of experiments, we train a sparse autoencoder with imposed sparsity regularisation during training on the Fashion-MNIST-1K dataset. For the lifted Bregman training approach (LBN) we apply Algorithm \ref{alg:implicit-stochastic-lifted-bregman-algorithm} to minimise the learning objective \eqref{eq:mlp-objective_sparse_denoising_autoencoder} with $\alpha= 0.09$. We set the batch size to $|B_k|=20$, the step-size parameter $\tau_k$ to $\tau_{k}= 0.5$, for all $k$, and the inner iterations $N$ to $N = 15$ for solving each mini-batch sub-problem \eqref{eq:implicit-sgd-for-many-data-points-objective} via \eqref{eq:bregman-proximal-network-training-algorithm}. For comparison, we consider a vanilla stochastic (sub-)gradient method (SGD-BP) as described in \eqref{eq:stochastic-gradient-descent} and the stochastic (sub-)gradient descent approach with implicit parameter update (ISGD-BP) that follows \eqref{eq:implicit-sgd} to train the network parameters. Both approaches use the learning rate $1 \times 10^{-3}$ and apply the back-propagation Algorithm \ref{alg:backprop-algorithm} for computing the (sub-)gradients. The ISGD-BP does so by solving the inner problem with GD-BP for $N$ iterations. 
% \begin{figure}[!ht]
%     \centering
%     \includegraphics[width=0.4\textwidth]{figures/Numerical_Denoising_Sparse_Autoencoder_1K_FMNIST_objective.png}\hspace{0.5cm}
%     \includegraphics[width=0.4\textwidth]{figures/Numerical_Denoising_Sparse_Autoencoder_1k_FMNIST_sparsity.png}
%     \caption{Sparse denoising autoencoder trained on Fashion-MNIST-1K images. \textbf{Left:} Objective values over epochs for each learning approach described in Section \ref{sec:numerical-results-denoising-sparse-autoencoder}. The training loss curves (solid lines) record MSE reconstruction error plus $\alpha$ times $\ell_1$-norm regularisation while the validation loss curves (dashed lines) report the MSE reconstruction error. \textbf{Right:} Sparsity level of the code per training epochs for each learning scheme.}
%     \label{figure:sparse-denoising-autoencoder-1k-fmnist-objective-sparsity-compare}
% \end{figure}

\begin{figure}[!t]
    \centering
    \begin{minipage}[t]{0.49\textwidth}
    \centering
    \includegraphics[width=0.49\textwidth]{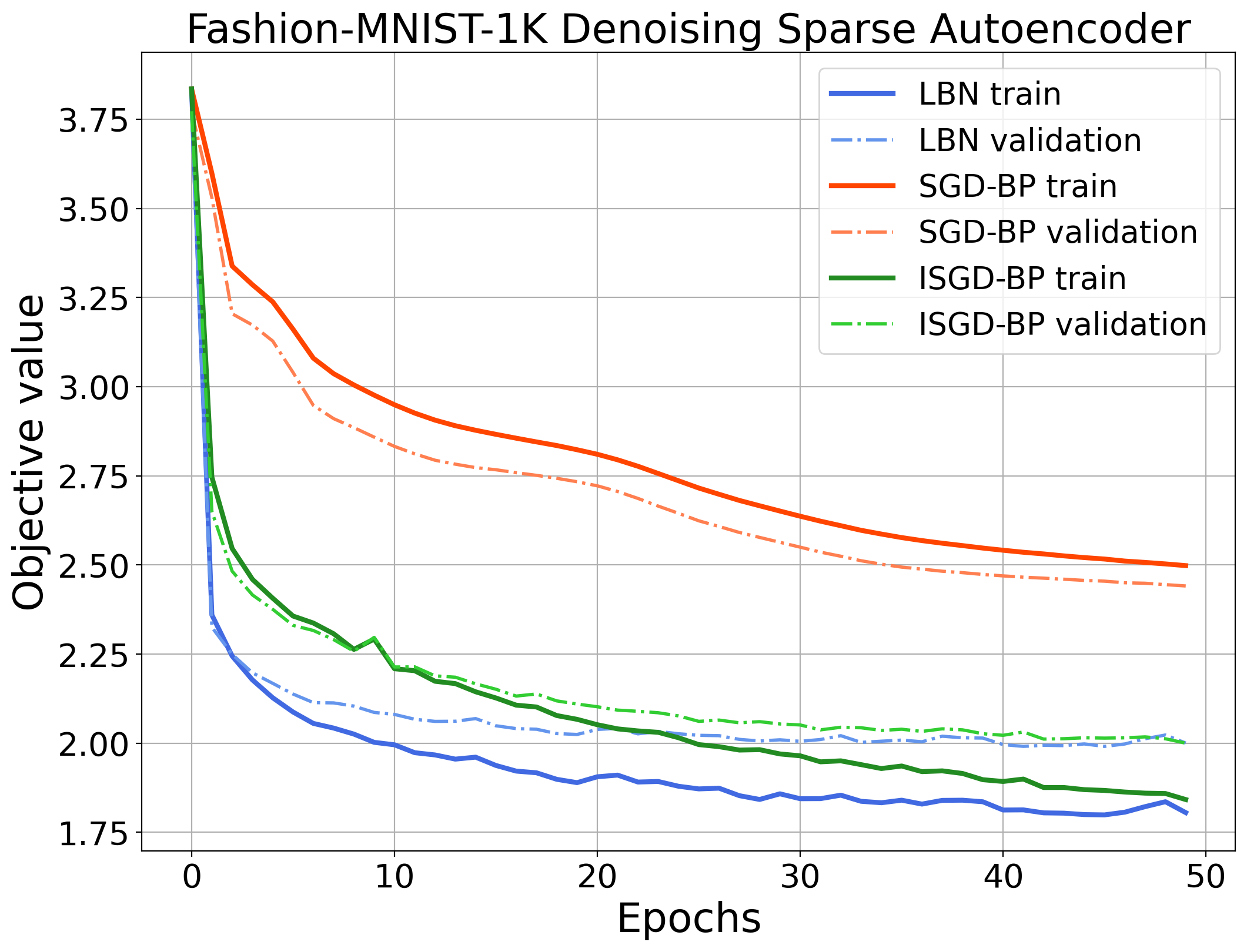}
    \includegraphics[width=0.49\textwidth]{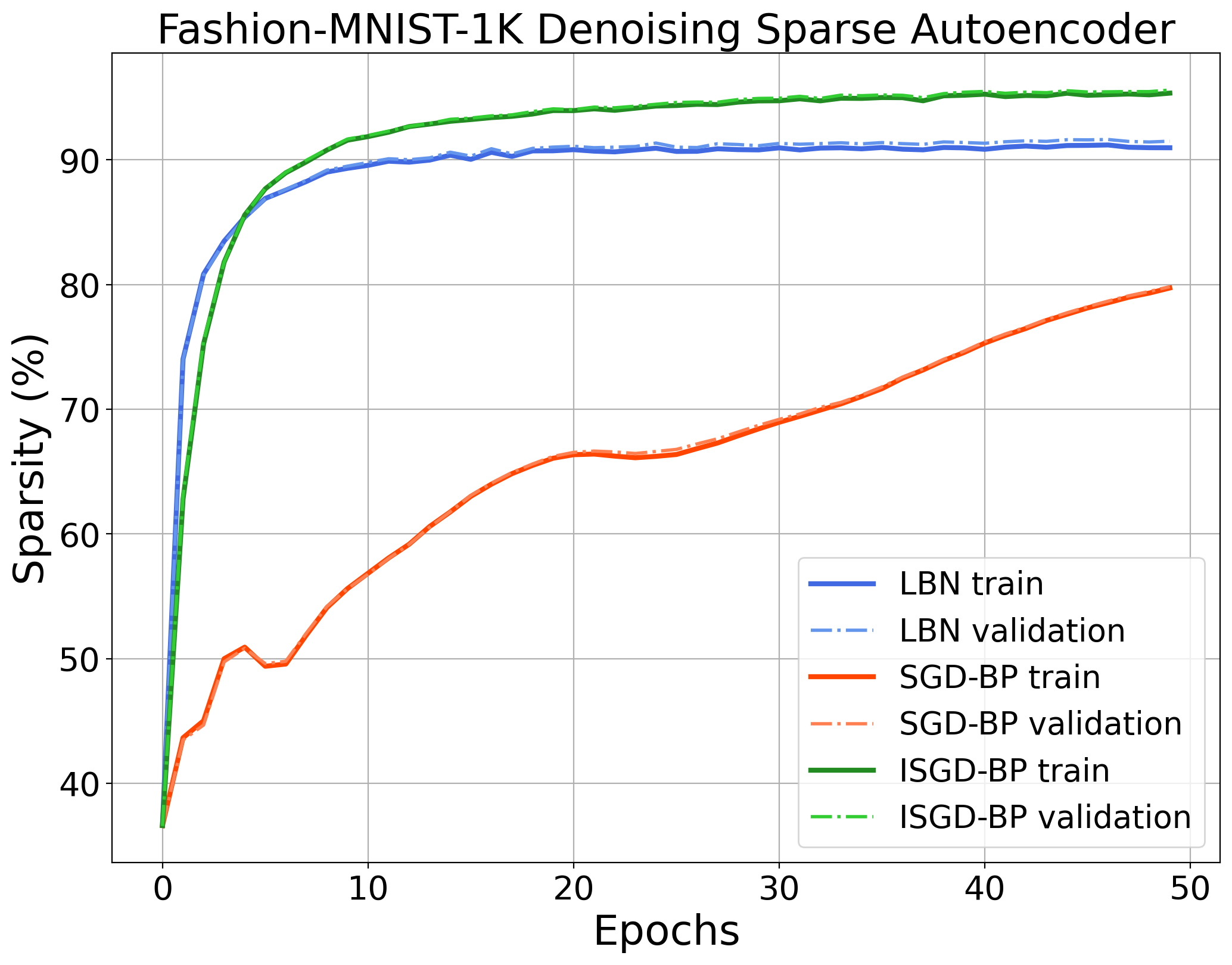}
    \caption{Sparse denoising autoencoder trained on Fashion-MNIST-1K images. \textbf{Left:} Objective values per epoch for each learning approach. The training loss record MSE reconstruction error plus $\alpha$ times $\ell_1$-norm regularisation while the validation loss record the MSE reconstruction error. \textbf{Right:} Sparsity level of the code per training epoch for each learning scheme.}
    \label{figure:sparse-denoising-autoencoder-1k-fmnist-objective-sparsity-compare}
    \end{minipage}\hfill
    \begin{minipage}[t]{0.49\textwidth}
    \centering
    \includegraphics[width=0.49\textwidth]{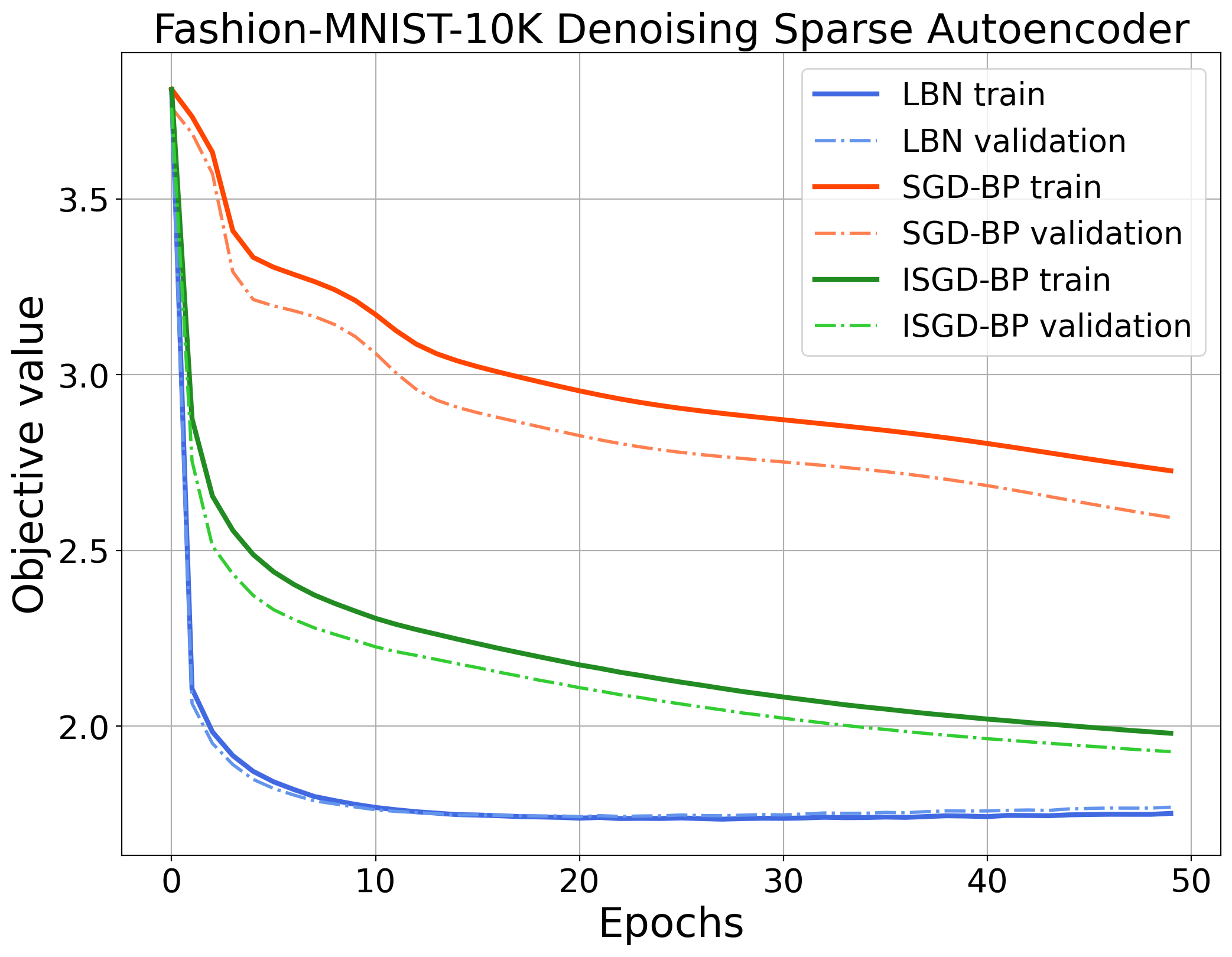}
    \includegraphics[width=0.49\textwidth]{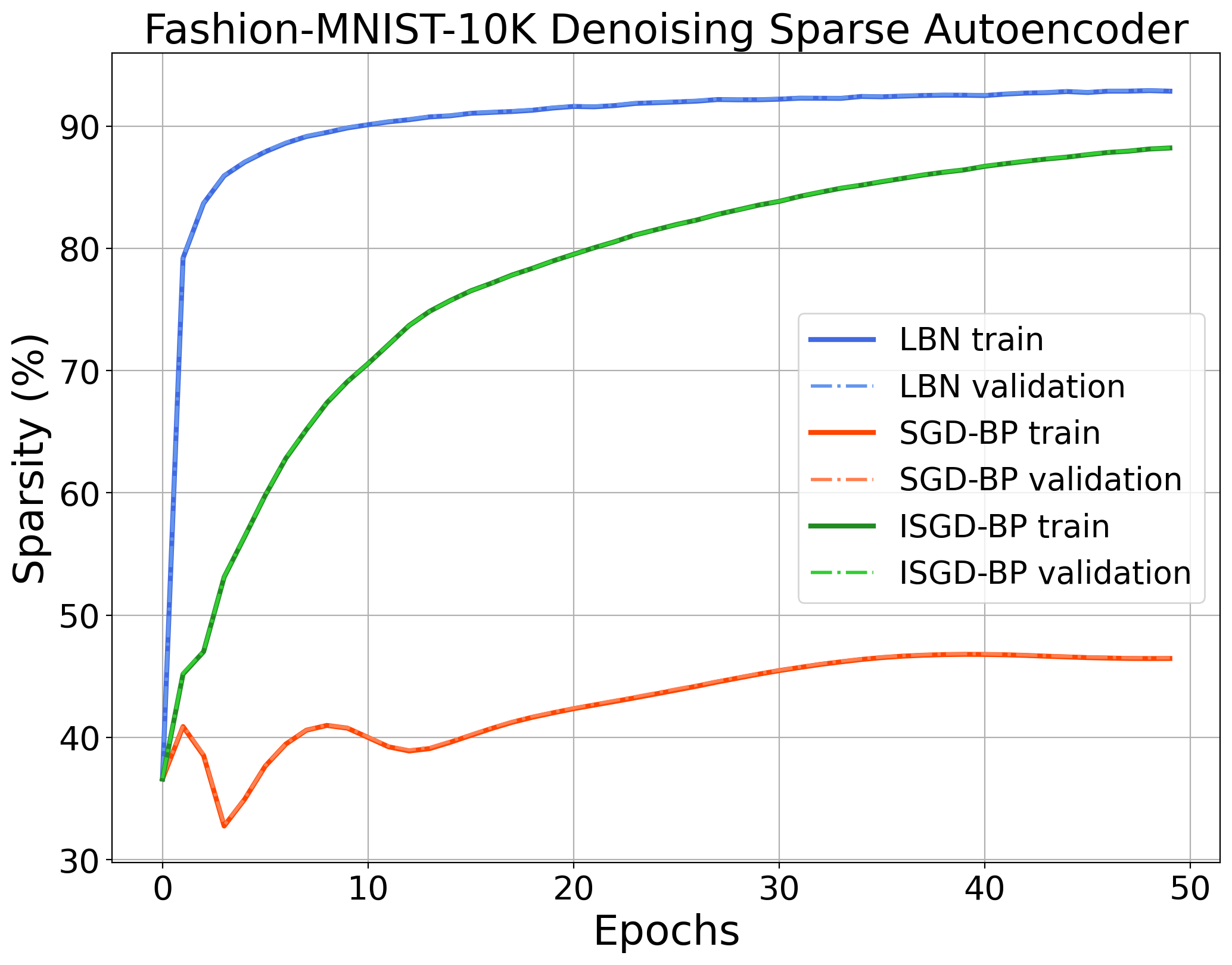}
    \caption{Sparse denoising autoencoder trained on the Fashion-MNIST-10K images. \textbf{Left:} Objective values per epoch for each learning approach. The training loss record MSE reconstruction error plus $\alpha$ times $\ell_1$-norm regularisation while validation loss curves record the MSE reconstruction error. \textbf{Right:} Sparsity level of the code over the 50 training epochs for each learning scheme.}
    \label{figure:sparse-denoising-autoencoder-10k-fmnist-objective-sparsity-compare}
    \end{minipage}
\end{figure}

All approaches train for 50 epochs. A log-scale plot of the objective value decay for all three training approaches is visualised in Figure \ref{figure:sparse-denoising-autoencoder-1k-fmnist-objective-sparsity-compare} (Left) and the tracked changes of the hidden node sparsity level are plotted in Figure \ref{figure:sparse-denoising-autoencoder-1k-fmnist-objective-sparsity-compare} (Right).

When trained on the Fashion-MNIST-1K dataset, we observe that LBN achieves faster convergence and outperforms ISGD-BP and SGD-BP especially at early training epochs. Towards later epochs, ISGD-BP starts to align and eventually achieves comparable training and validation errors. %Likely due to the higher dataset complexity and limitation on trainset size, both ISGD-BP and LBN struggle with further reducing validation loss values towards later epochs. 

The ISGD-BP approach and LBN approach are capable of achieving and maintaining higher sparsity levels from earlier epochs onwards, while the SGD-BP approach on the other hand seems to increase the sparsity level much more slowly compared to the other two approaches. It should, however, be emphasised that the SGD-BP algorithm is computationally less expensive per epoch.

\textbf{Fashion-MNIST-10K} 
We conduct the second set of experiments on the Fashion-MNIST-10K dataset. We choose a bigger batch size of $|B_k|=200$ and penalise the $\ell_1$-norm regularisation with $\alpha= 0.055$ during training. For both the LBN and the ISGD-BP approach, we set $\tau_{k} = 1$ and perform $N = 30$ iterations in each mini-batch sub-problem. We choose a learning rate of $4\times 10^{-4}$ for the ISGD-BP approach as well as for the SGD-BP approach.

% \begin{figure}[!ht]
%     \centering
%     \includegraphics[width=0.4\textwidth]{figures/Numerical_Denoising_Sparse_Autoencoder_10K_FMNIST_objective.png}\hspace{0.5cm}
%     \includegraphics[width=0.4\textwidth]{figures/Numerical_Denoising_Sparse_Autoencoder_10K_FMNIST_sparsity.png}
%     \caption{Sparse denoising autoencoder trained on the Fashion-MNIST-10K images. \textbf{Left:} Objective values per epoch for each learning scheme. The training loss curves record MSE reconstruction error plus $\alpha$ times $\ell_1$-norm regularisation while validation loss curves record the MSE reconstruction error. \textbf{Right:} Sparsity level of the code over the 50 training epochs for each learning scheme.}
%     \label{figure:sparse-denoising-autoencoder-10k-fmnist-objective-sparsity-compare}
% \end{figure}

In Figure \ref{figure:sparse-denoising-autoencoder-10k-fmnist-objective-sparsity-compare} we visualise the training and validation loss curves (Left) and the sparsity rate (Right) over the 50 training epochs when training the proposed algorithms on the Fashion-MNIST-10K dataset. As expected, the gap between training and validation loss curves in this experiment are closer due to larger amount of training samples while network complexity is unchanged. Note that the proposed LBN approach exhibits stronger performance and outperforms the other two (sub-)gradient based methods when a larger amount of training data is available. Not only is the LBN approach capable of achieving lower reconstruction loss values after 50 epochs but also achieves and maintains higher sparsity rates more quickly. 

\begin{figure}
    \centering
    \includegraphics[width=0.32\textwidth]{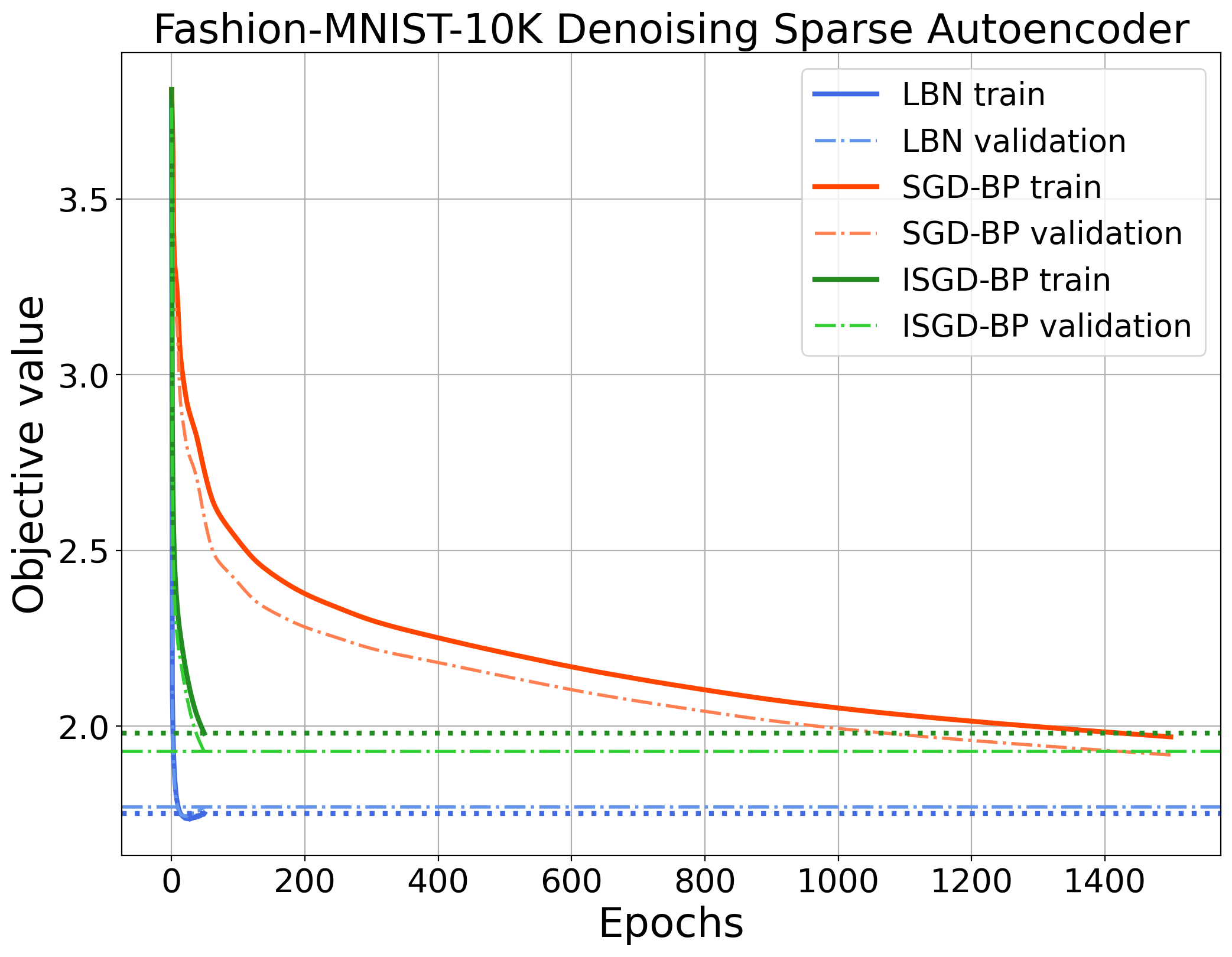}
    \includegraphics[width=0.32\textwidth]{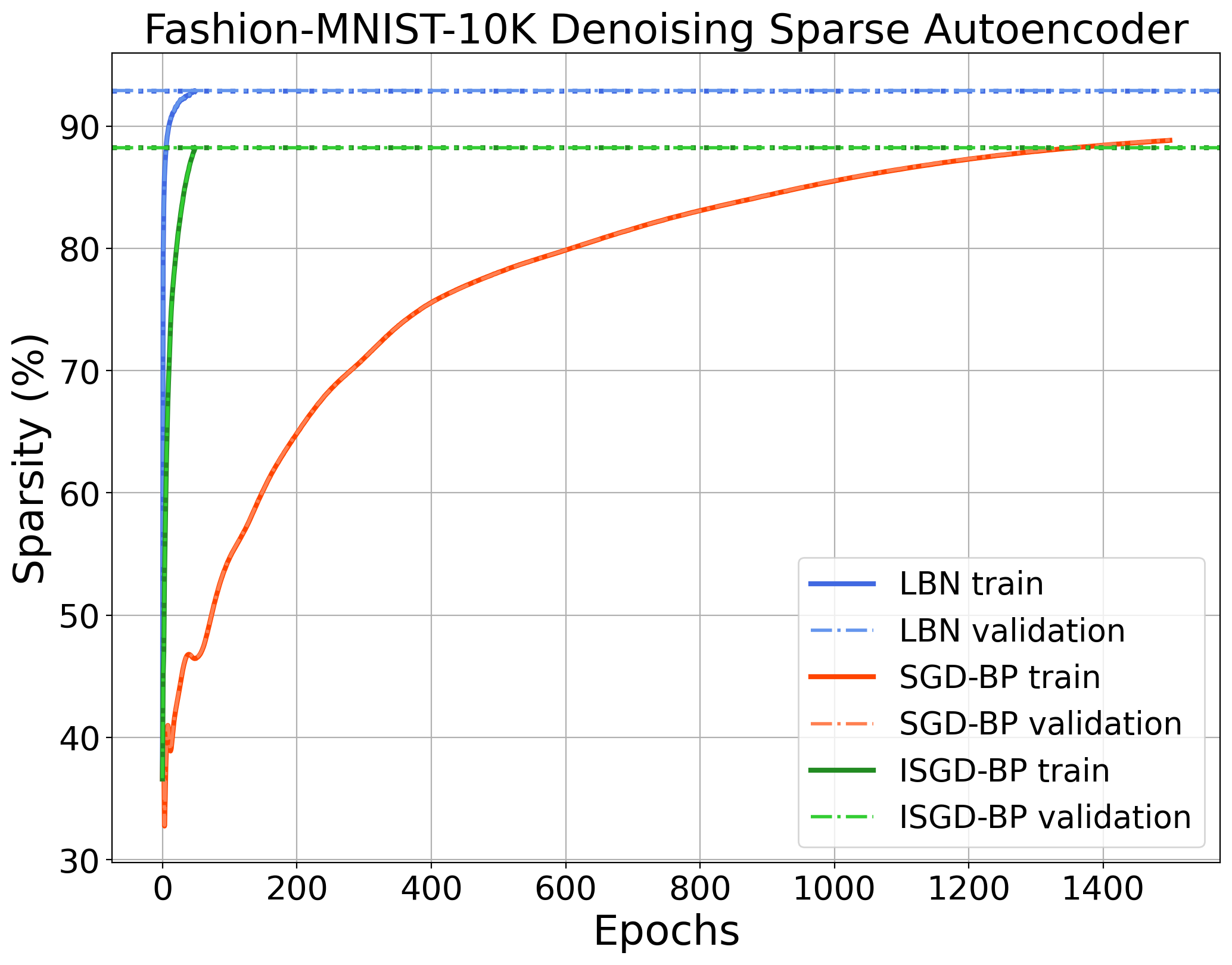}
    \includegraphics[width=0.32\textwidth]{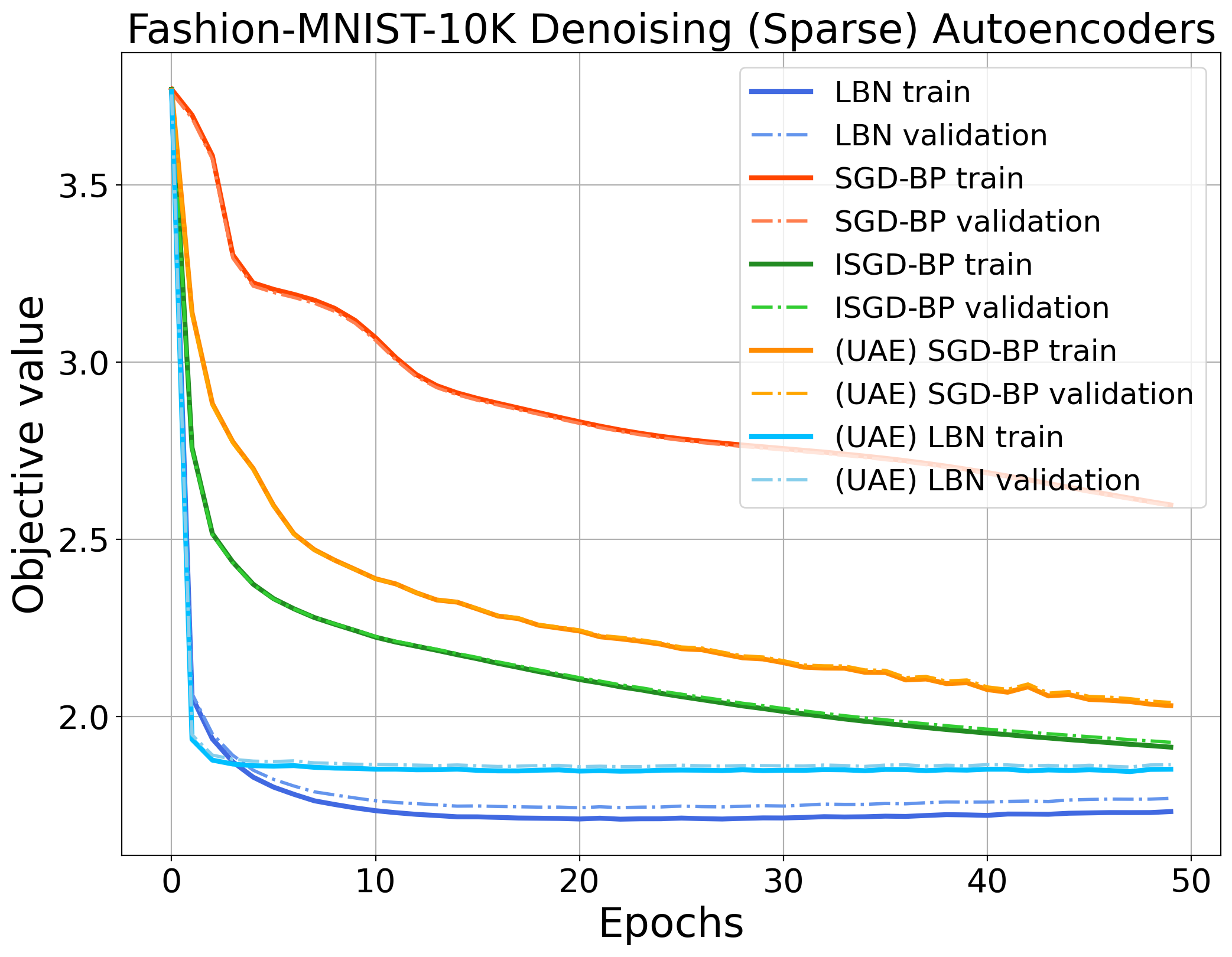}
    \caption{Comparison of sparse denoising autoencoders trained on the Fashion-MNIST-10K dataset with LBN, ISGD-BP and SGD-BP for 1500 epochs, and of sparse denoising autoencoders with undercomplete autoencoders. \textbf{Left: } Objective values per epoch for all learning approaches. Horizontal lines mark end-of-training and end-of-validation objective values of LBN and ISGD-BP. \textbf{Middle: } Sparsity level of the code over the 1500 training epochs. Horizontal lines show end-of-training sparsity rates of LBN and ISGD-BP. \textbf{Right: } Objective values per epoch for each learning approach. The training and validation loss curves report MSE reconstruction errors over 50 epochs.}
    \label{figure:sparse-denoising-autoencoder-fmnist-objective-sparsity-1500E-compare}
\end{figure}

As both the LBN and ISGD-BP approach require 30 inner iterations per batch sub-problem, we also compare both approaches to the SGD-BP approach for 1500 epochs to ensure that the total number of iterations for all three approaches are comparable. As visualised in Figure \ref{figure:sparse-denoising-autoencoder-fmnist-objective-sparsity-1500E-compare}, we verify with fairer comparison with regards to runtime that the SGD-BP approach matches training and validation errors of ISGD-BP, but falls short of achieving training and validation values as low as those obtained with LBN. 

%Enforcing sparsity instead of compressing the hidden dimension allows more flexible representation of signals, and more dynamically and adaptive representation of the signals. 

% \begin{figure}
%     \centering
%     \includegraphics[scale=0.285]{figures/Numerical_Denoising_Sparse_Autoencoder_10K_FMNIST_objective_1500_Epochs.png}
%     \includegraphics[scale=0.285]{figures/Numerical_Denoising_Sparse_Autoencoder_10K_FMNIST_sparsity_1500_Epochs.png}
%     \caption{A further comparison of training a sparse denoising autoencoder with LBN, ISGD-BO and SGD-BP for 1500 epochs. \textbf{Left: } Objective values per epoch for all learning approaches. Horizontal lines mark end-of-training and end-of-validation objective values of LBN and ISGD-BP. \textbf{Right: } Sparsity level of the code over the 1500 training epochs. Horizontal lines show end-of-training sparsity rates of LBN and ISGD-BP.}
%     \label{figure:sparse-denoising-autoencoder-fmnist-objective-sparsity-1500E-compare}
% \end{figure}

% \begin{figure}[!t]
%     \centering
%     \begin{minipage}[t]{0.59\textwidth}
%     \centering
%     \includegraphics[width=0.49\textwidth]{figures/Numerical_Denoising_Sparse_Autoencoder_10K_FMNIST_objective_1500_Epochs.png}
%     \includegraphics[width=0.49\textwidth]{figures/Numerical_Denoising_Sparse_Autoencoder_10K_FMNIST_sparsity_1500_Epochs.png}
%     \caption{A further comparison of training a sparse denoising autoencoder with LBN, ISGD-BO and SGD-BP for 1500 epochs. \textbf{Left: } Objective values per epoch for all learning approaches. Horizontal lines mark end-of-training and end-of-validation objective values of LBN and ISGD-BP. \textbf{Right: } Sparsity level of the code over the 1500 training epochs. Horizontal lines show end-of-training sparsity rates of LBN and ISGD-BP.}
%     \label{figure:sparse-denoising-autoencoder-fmnist-objective-sparsity-1500E-compare}
%     \end{minipage}\hfill
%     \begin{minipage}[t]{0.39\textwidth}
%     \centering
%     \includegraphics[width=0.75\textwidth]{figures/Numerical_Denoising_Autoencoders_10K_FMNIST_objective.png}
%     \caption{Comparison of sparse denoising autoencoder and undercomplete autoencoder trained on the Fashion-MNIST-10K dataset: Objective values per epoch for each learning approach. The training and validation loss curves report MSE reconstruction errors over 50 epochs. }
%     \label{figure:sparse-denoising-autoencoders-10k-fmnist-objective-compare}
%     \end{minipage}
% \end{figure}

\textbf{Comparison with undercomplete autoencoders.} 
So far, we haven't properly motivated the use of sparse autoencoders in comparison to traditional autoencoders that explicitly reduce the dimension of the code. A sparse code with $m$ non-zero entries requires basically the same amount of memory than a $m$-dimensional code but offers greater flexibility because the location of the non-zero entries can vary for different network inputs. This advantage should lead to better validation errors of sparse autoencoders compared to traditional, undercomplete autoencoders. In order to verify the advantages of training sparse autoencoders over training undercomplete autoencoders (UAE), we conduct two additional experiments where we train undercomplete autoencoders on the Fashion-MNIST-10K dataset without imposing any additional regularisation.

We observe that the sparse denoising autoencoder trained with $\alpha=0.055$ eventually reaches 91\% sparsity rate after 50 epchs, which suggests that on average 70 nodes in the code of the final model are activated. Hence, we train a 4-layer fully connected undercomplete autoencoder with ReLU activation function for every layer but set the number of nodes in the middle layer to 70, which implies $W_2 \in \mathbb{R}^{784 \times 70}$ and $W_3 \in \mathbb{R}^{70 \times 784}$. The undercomplete autoencoder is trained with both LBN and SGD-BP and we compare if a reconstruction quality comparable to the sparse autoencoder can be achieved. 

% \begin{figure}[!ht]
%     \centering
%     \includegraphics[scale=0.285]{figures/Numerical_Denoising_Autoencoders_10K_FMNIST_objective.png}
%     \caption{Comparison of sparse denoising autoencoder and undercomplete autoencoder trained on the Fashion-MNIST-10K dataset: Objective values per epoch for each learning approach. The training and validation loss curves report MSE reconstruction errors over 50 epochs. }
%     \label{figure:sparse-denoising-autoencoders-10k-fmnist-objective-compare}
% \end{figure}

In Figure \ref{figure:sparse-denoising-autoencoder-fmnist-objective-sparsity-1500E-compare}, we validate that the lifted Bregman approach helps to improve both undercomplete and sparse autoencoder to achieve faster convergence and lower objective values compared to the (sub-)gradient based methods. We also confirm that after 50 epochs the sparse denoising autoencoder achieves lower training and validation errors compared to the UAE approach.

%In Figure \ref{figure:sparse-denoising-autoencoders-10k-fmnist-objective-compare}, we validate that the lifted Bregman approach helps to improve both undercomplete and sparse autoencoder to achieve faster convergence and lower objective values compared to the (sub-)gradient based methods. We also confirm that after 50 epochs the sparse denoising autoencoder achieves lower training and validation errors compared to the UAE approach.

More specifically, the UAE model trained with the lifted Bregman approach sees bigger reconstruction loss decreases in the early epochs, but is outperformed by the sparse autoencoder at later epochs. This experiment seems to suggest that by leveraging the power of sparsity, sparse autoencoders trained via lifted Bregman approaches are capable of finding more flexible data representations compared to autoencoders with explicit dimension reduction.

\begin{figure}[!t]
    \centering
    \begin{minipage}[t]{0.49\textwidth}
    \centering
    \includegraphics[width=\textwidth]{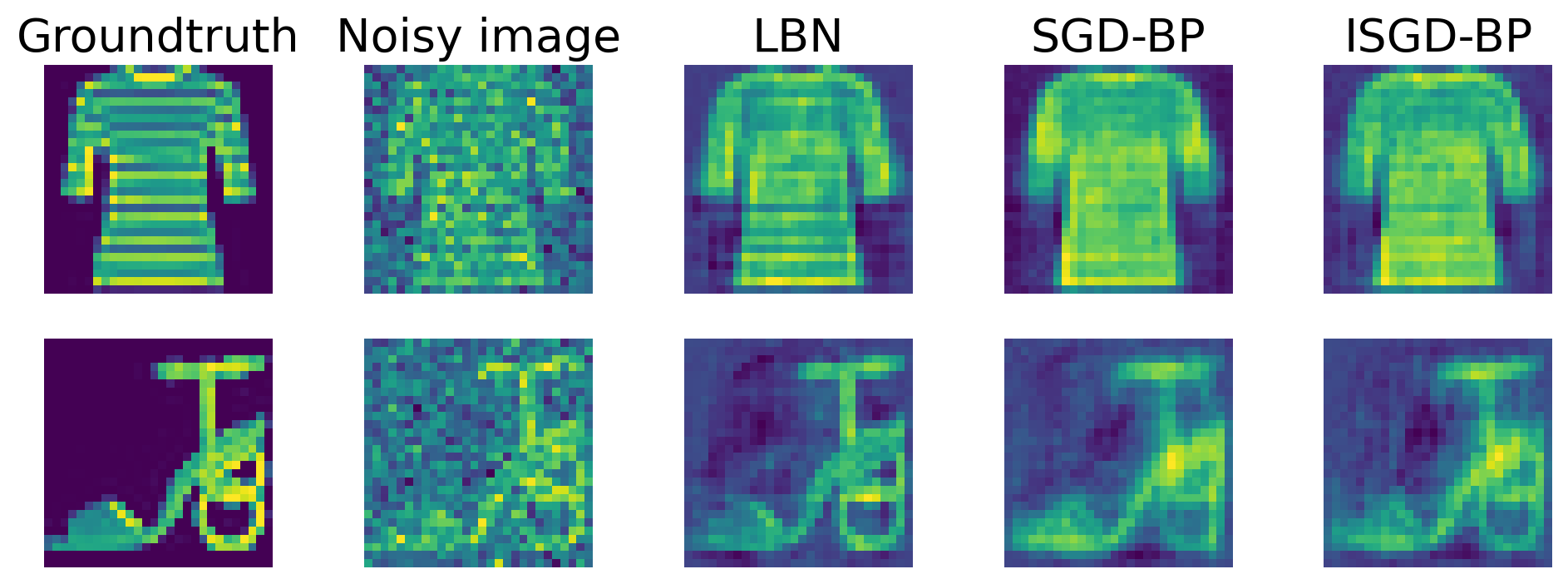}
    \caption{Denoised images from the Fashion-MNIST-1K training dataset, computed with the LBN, the SGD-BP and the ISGD-BP training strategies for 50 epochs, along with the noise-free and noisy images.}
    \label{figure:sparse-autoencoder-1k-fmnist-denoised-train-samples}
    \end{minipage}\hfill
    \begin{minipage}[t]{0.49\textwidth}
    \centering
    \includegraphics[width=\textwidth]{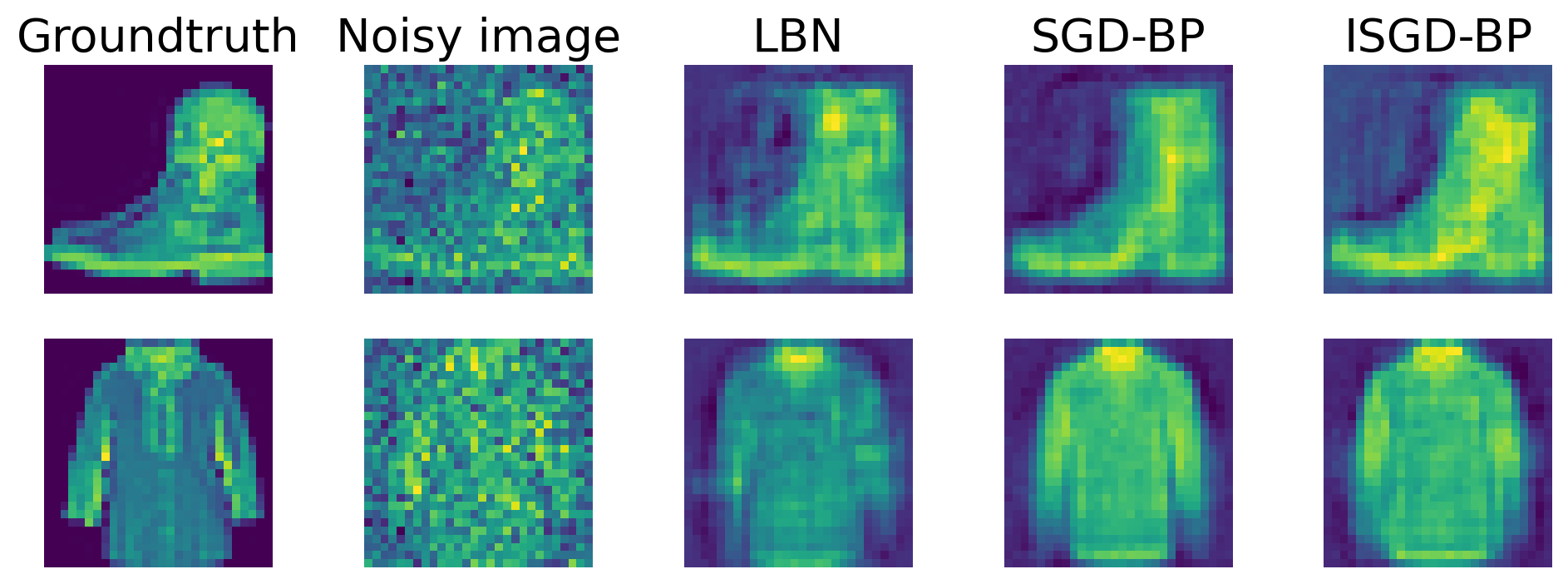}
    \caption{Denoised images from the Fashion-MNIST-1K validation dataset computed with the LBN, the SGD-BP and the ISGD-BP training strategies for 50 epochs, along with the noise-free and noisy images.}
    \label{figure:sparse-autoencoder-1k-fmnist-denoised-val-samples}
    \end{minipage}
\end{figure}

\begin{figure}[!t]
    \centering
    \begin{minipage}{0.49\textwidth}
    \includegraphics[width=\textwidth]{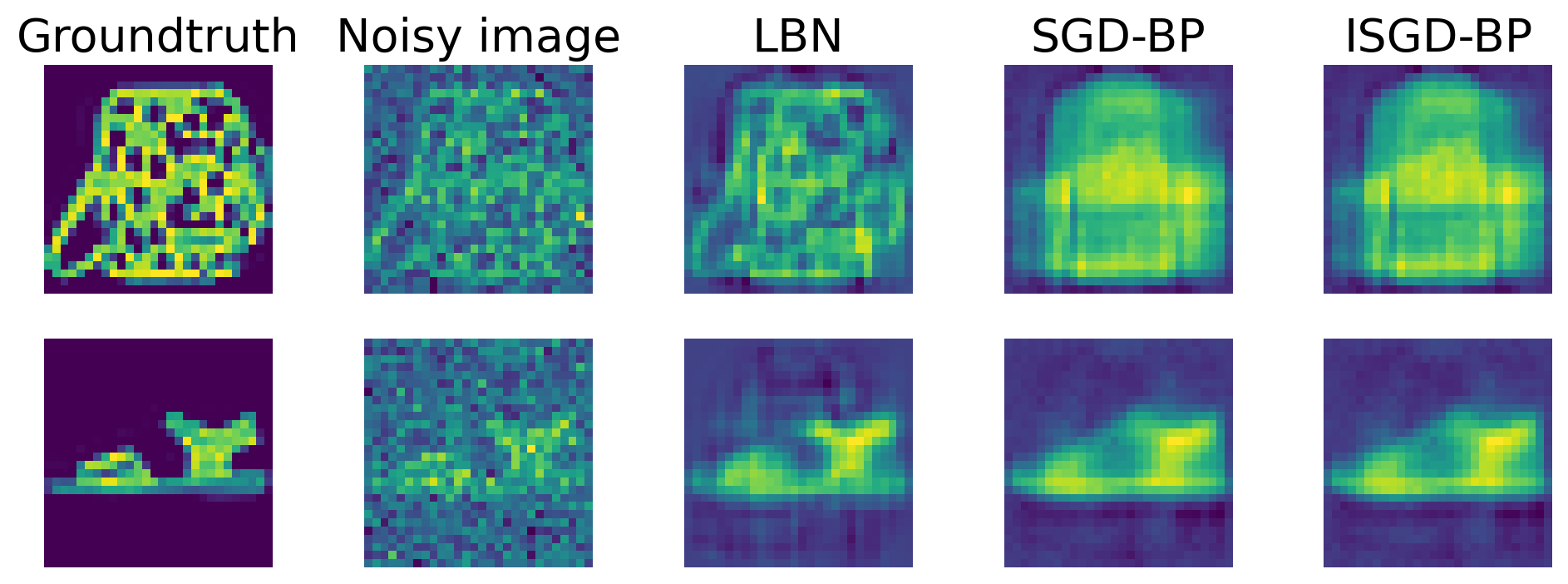}
    \caption{Denoised images from the Fashion-MNIST-10K training dataset computed with the LBN, the SGD-BP (trained for 1500 epochs) and the ISGD-BP training strategies respectively, along with the noise-free and noisy images.}
    \label{figure:sparse-autoencoder-fmnist-denoised-train-samples}
    \end{minipage}\hfill
    \begin{minipage}{0.49\textwidth}
    \centering
    \includegraphics[width=\textwidth]{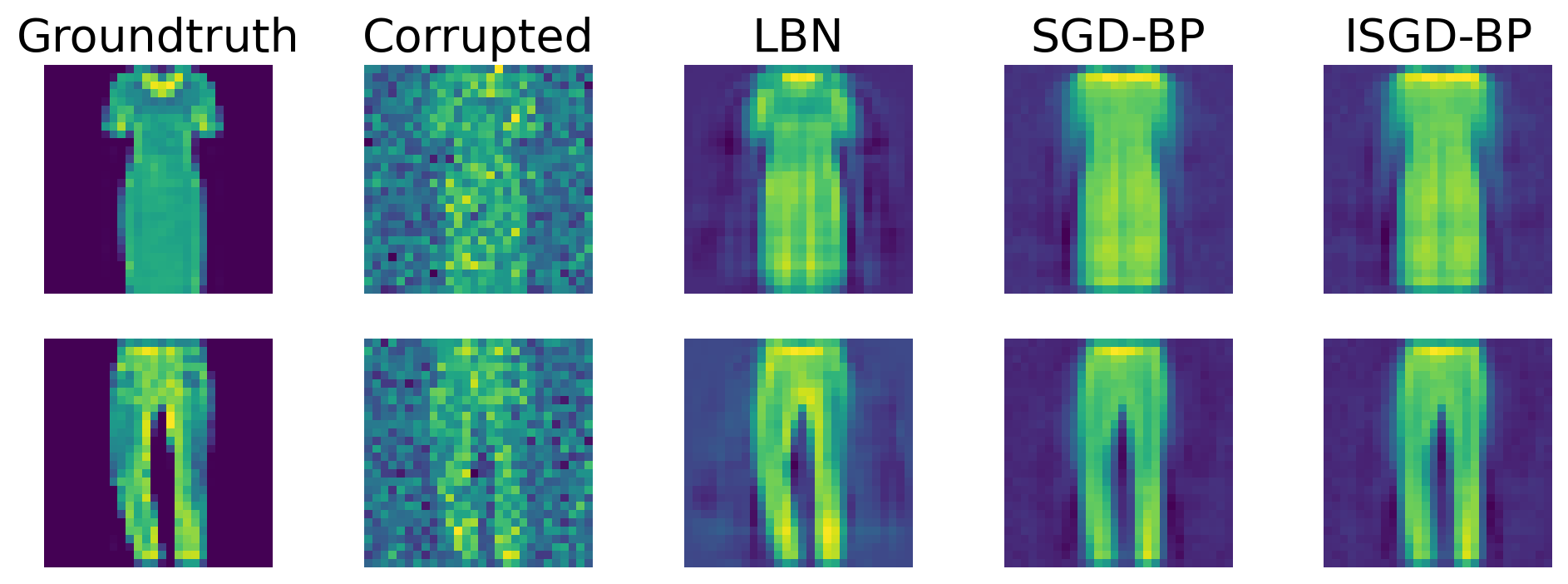}
    \caption{Denoised images from the Fashion-MNIST-10K validation dataset computed with the LBN, the SGD-BP (trained for 1500 epochs) and the ISGD-BP training strategies respectively, along with the noise-free and noisy images.}
    \label{figure:sparse-autoencoder-fmnist-denoised-val-samples}
    \end{minipage}
\end{figure}

We visualise a selection of denoised sample images from the Fashion-MNIST-1K training dataset in Figure \ref{figure:sparse-autoencoder-1k-fmnist-denoised-train-samples} and from the Fashion-MNIST-10K training dataset in Figure \ref{figure:sparse-autoencoder-fmnist-denoised-train-samples}. Although the loss values of LBN and ISGD-BP were comparable after 50 epochs, the visual comparison between the network outputs suggests that the network trained with LBN is able to restore finer details (for instance the stripes on the t-shirt).

The same network also exhibits better denoising performance for the validation images as seen in Figure \ref{figure:sparse-autoencoder-1k-fmnist-denoised-val-samples}. 

Similar observations can be made when we examine denoised sample images from the models trained on the Fashion-MNIST-10K dataset (Figure \ref{figure:sparse-autoencoder-fmnist-denoised-train-samples}). The LBN trained network is capable of recovering fine structures such as the dotted backpack object. Both SGD-BP and ISGD-BP struggle with fine details and are prone to producing coarse approximations. The observations still hold true for the validation dataset as can be seen from Figure \ref{figure:sparse-autoencoder-fmnist-denoised-val-samples}. For example, the LBN approach is able to recover the curvature of the trousers, but SGD-BP and ISGD-BP only recover straight-legged trousers. 
%We suspect this is due to the a restricted model capacity by the imposed sparsity constraints. 
%but maybe the safe guess is what makes sgd/gd/isgd based methods more generalisable?

% \begin{figure}[!ht]
%     \centering
%     \includegraphics[scale=0.45]{figures/Numerical_Denoise_SAE_10K_FMNIST_Denoised_2_Train_Samples.png}
%     \caption{Denoised images from the Fashion-MNIST-10K training dataset computed with the LBN, the SGD-BP (trained for 1500 epochs) and the ISGD-BP training strategies respectively, along with the noise-free and noisy images.}
%     \label{figure:sparse-autoencoder-fmnist-denoised-train-samples}
% \end{figure}

% \begin{figure}
%     \centering
%     \includegraphics[scale=0.45]{figures/Numerical_Denoise_SAE_10K_FMNIST_Denoised_2_Val_Samples.png}
%     \caption{Denoised images from the Fashion-MNIST-10K validation dataset computed with the LBN, the SGD-BP (trained for 1500 epochs) and the ISGD-BP training strategies respectively, along with the noise-free and noisy images.}
%     \label{figure:sparse-autoencoder-fmnist-denoised-val-samples}
% \end{figure}

\section{Conclusions \& outlook}\label{sec:conclusions-outlook}

In this work, we proposed a novel framework for learning parameters of feed-forward neural networks. To achieve that, we introduced a new loss function based on Bregman distances and we provided a detailed mathematical analysis of this new loss function. By replacing the quadratic penalties in the MAC-QP approach, we generalised MAC-QP and lifted network frameworks as introduced in \citep{carreira2014distributed,askari2018lifted,li2019lifted,gu2020fenchel} and established rigorous mathematical properties. 

The advantage of the proposed framework is that it formulates the training of a feed-forward network as a minimisation problem in which computing partial derivatives of the network parameters does not require differentiation of the activation functions. We have demonstrated that this framework can be realised computationally with a variety of different deterministic and stochastic minimisation algorithms. 

%This framework also circumvents potential drawbacks such as the vanishing gradient problem and overcomes the limitation of the non-parallelisability which the back-propagation algorithm in standard training approaches suffers. 

%This based on a special type of Bregman distances. One prominent advantage of the new formulation is that  we do not require the computations of derivatives with respect to the individual (non-smooth) activation functions when taking partial derivatives with respect to the network's parameters. Also by breaking down the learning problem into layer-wise sub-problems, our proposed approach allows for many parallelisable or distributed optimisation approaches. 

With numerical experiments in classification and compression, we compared the implementation of the proposed framework with proximal gradient descent and implicit stochastic gradient method to explicit and implicit stochastic gradient methods with back-propagation. 

For classification, we have shown that the trained network is fully non-linear in contrast to the lifted network approach proposed in \cite{askari2018lifted}, while also achieving similar classification accuracy as networks trained with first- order methods with back-propagation. 

For compression, we showed with our sparse autoencoder example that training parameters with proximal gradient descent or implicit stochastic gradient descent can be more efficient than first-order methods with back-propagation. We have also shown that it is straightforward to integrate regularisation acting on the activation variables. This enabled us to design an autoencoder with sparse codes that is shown to produce more adaptive compressions compared to conventional autoencoders.

There are many possible directions for future work. In this work, the proposed concept is limited to feed-forward neural network architectures, but it can easily be extended to other types of architectures such as residual neural networks \citep{he2016deep}. Residual networks of the form $x_{l + 1} = x_l - \sigma(f(x_l, \Theta_l))$ can for instance be trained with a slight modification of \eqref{eq:bregman-lifted} in the form of
\begin{align*}
    \min_{\Thetas,\Xs}  \sum_{i=1}^{s} \left[\ell\left(y^i, x_{\layer}^{i}\right) + \sum_{l=1}^{\layer-1} \,  \lambda_l \bregmanloss\left(x_{l}^{i} - x_{l - 1}^{i}, f(x_{l - 1}^{i}, \Theta_l) \right)\right] \, ,
\end{align*}
and many modifications with different linear transformations of the auxiliary variables.

Other notable research directions include feasibility studies on how to incorporate other popular layers, such as pooling layers, batch normalisation layers \citep{ioffe2015batch,ba2016layer}, or transformer layers \citep{vaswani2017attention}, the application of other optimisation methods such as stochastic coordinate descent, stochastic ADMM, or stochastic primal-dual hybrid gradient methods as well as tailored non-convex optimisation methods and other regularisation techniques like dropout \cite{hinton2012improving}. 

Similar to suggestions made in \cite{zach2019contrastive}, we could use the network energy for anomaly detection. As our proposed framework has opened more pathways to use other possible combinations of optimisation strategies, it would also be of particular interest to investigate more thoroughly if good generalisation properties of a deep neural network are the result of the chosen optimisation strategy or the chosen network architecture itself. 

\section*{Acknowledgements}
The authors would like to thank the Isaac Newton Institute for Mathematical Sciences, Cambridge, for support and hospitality during the programme \emph{Mathematics of Deep Learning} where work on this paper was undertaken. This work was supported by EPSRC grant no EP/R014604/1. XW acknowledges support from the Cantab Capital Institute for the Mathematics of Information and the Cambridge Centre for Analysis (CCA).

%\begin{itemize}
%    \item Conclusions here
%    \item Mention limitations
%    \item Mention open questions
%    \item Mention future work
%\end{itemize}

%\mb{TODO for Victor, include plots you showed me, and also}
%\begin{enumerate}
%    \item Include plots you showed me
%    \item Compare with SGD for 100 epochs to see if it catches up with LBN
%    \item Compare runtimes for different algorithms (total runtimes, runtimes per epoch etc.)
%    \item Compare validation errors between sparse Autoencoder and conventional Autoencoder with dimension of encoding similar to number of non-zero elements in encoding of sparse autoencoder.
%    \item Try same results for denoising autoencoder
%    \item Try different datasets (FashionMNIST, Cifar10, yale face dataset)
%    \item compare classification on fmnist and mnist dataset
%    \item compare mnist dataset on reconstruction task, small data samples check %validation performance
%    \item compare fmnist dataset on denosing sparse ae task, compare against sgd, lower loss of ours and higher sparse rate, 
%    \item compare against implicit sgd also ours lower loss, perhaps high sparsity for implicit bp sgd? 
%    \item compare against plain sgd and vanilla autoencoder, although low loss but reconstruction result not good, not detail
%\end{enumerate}

%%%%%%%%%%%%%%%%%%%%%%%%%%%%%%%%%%%%%%%%%%%%%%%%%%%%%%%%%%%%%%%%%

%\break 

%\pagebreak 

%\bibliographystyle{natbib}
\bibliography{references.bib}

\appendix 

\section{Back-propagation}\label{sec:backprop}

\begin{algorithm}[!t]
\caption{Back-Propagation Algorithm}   
\begin{algorithmic}
\label{alg:backprop-algorithm}
    \STATE Input $\Thetas$ \;
        \FOR{$i = 1, \ldots, s $}
            \FOR{$l = 1, \ldots, \layer$}
                \STATE Perform forward pass to compute $z_{l}$ and $x_{l}$:\\
                \STATE $z_{l}^{i} = f(x_{l-1}^{i},\Theta_{l})$ \\
                \STATE $x_{l}^{i} = \sigma_{l} (z_{l}^{i})$ \\
            \ENDFOR
        \ENDFOR
        \FOR{$i = 1, \ldots, s $}
            \FOR{$l = \layer, \ldots, 1$}
                \STATE Perform backward pass to compute $\delta_{l}^{i}$:\\
                \STATE  $\delta_{l}^{i} =
                \begin{cases}
                \sigma_{l}'(z_{\layer}^i) \odot \nabla_{x} \ell(y^{i}, x_{\layer}^{i}) \;\;\; \text{for } l = \layer\\
                \sigma_{l}'(z_{l}^{i}) \odot \frac{\partial f}{\partial x_{l}} \delta_{l+1}\;\;\; \text{for } l \in \{1, \ldots, \layer-1\}\\
                \end{cases}$ \\
        \ENDFOR
    \ENDFOR
%\STATE Partial derivatives: $\frac{\partial L}{\partial b_{l}^{j}} = \delta_{l}^{j} \;\;\; \text{for} j \in {1, \ldots n_{l}} $\\
\STATE Partial derivatives: $\frac{\partial \ell}{\partial \Theta_{l}^{j k}} = \delta_{l}^{j} \frac{\partial f}{\partial \Theta_{l}^{k}} \;\;\; \text{for } j \in \{1, \ldots, n_{l}\} \text{ and } k \in \{1,\ldots, n_{l-1} \}$\\
\end{algorithmic}
\end{algorithm}

\noindent We can write down the corresponding Lagrangian function for \eqref{eq:mlp-constrained} as
\begin{equation}\label{eq:mlp-lagrangian}
    \mathcal{L}(\Thetas,\Xs,\Zs,\mu,\delta) = \sum_{i=1}^{s} \left[ \ell(y, x_{\layer}^{i})+ \sum_{l=1}^{\layer} \langle \mu_{l}^{i}, x_{l}^{i} - \sigma(z_l^{i}) \rangle  + \sum_{l=1}^{\layer} \langle \delta_{l}^{i}, z_{l}^{i} - f(x_{l-1}^{i},\Theta_{l}) \rangle\right] \,,
\end{equation}%
%\begin{align}
%    \mathcal{L}(\{W\}_{l=1}^{\layer},\{x\}_{l=1}^{\layer-1},\lambda_{l=1}^{\layer}) = \sum_{i=1}^{s} \left[ \ell(y^i, W_L^\top x^i_{\layer-1})+ \sum_{l=1}^{\layer-1} \langle \lambda_{l}, (x^i_{l} - \sigma(W_l^\top x^i_{l-1})) \rangle \right]
%\end{align}
where $\mu=\{\mu_l^{i}\}_{l=1,\dots, \layer}^{i=1,\dots, s}$ and $\delta=\{\delta_l^{i}\}_{l=1,\dots, \layer}^{i=1,\dots,s}$ are Lagrangian multipliers. 

In the following deduction of the back-propagation algorithm, we fix for one training data sample and drop the dependency on $i$ for the ease of notation. The optimality conditions of \eqref{eq:mlp-lagrangian} can be split into the following sub-conditions:
\begin{align}%\label{eq:mlp-lagrangian-optimality}
   %& \frac{\partial \mathcal{L}}{\partial \Theta_{\layer}} = 0 \implies  \nabla_{\Theta_{\layer}}\ell(y^i, x_{\layer}^{i}) = 0 \label{eq:mlp-lagrangian-optimality-W_L}\\
   & \frac{\partial \mathcal{L}}{\partial x_{\layer}} = 0 \implies  \nabla_{x_{\layer}} \ell(y,x_{\layer}) + \mu_{\layer} = 0  \;\; \label{eq:mlp-lagrangian-optimality-x_L-1}\\
   & \frac{\partial \mathcal{L}}{\partial \Theta_l} = 0 \implies \delta_{l} = \frac{\partial f}{\partial \Theta_{l}} \mu_{l}^{\top} \;\; \text{for } l = 1, \ldots \layer-1 \label{eq:mlp-lagrangian-optimality-w_l}\\
   & \frac{\partial \mathcal{L}}{\partial x_l} = 0 \implies \mu_{l} = \frac{\partial f}{\partial x_{l}} \delta_{l+1} \;\; \text{for } l = 1, \ldots \layer-1 \label{eq:mlp-lagrangian-optimality-x_l}\\
   & \frac{\partial \mathcal{L}}{\partial z_l} = 0 \implies  \delta_{l} = \sigma_{l}'(z_{l})\mu_{l}  \;\; \text{for } l = 1 \ldots \layer \label{eq:mlp-lagrangian-optimality-z_l}\\
   & \frac{\partial \mathcal{L}}{\partial \mu_l} = 0 \implies x_{l} = \sigma_{l}(f(x_{l-1},\Theta_{l})) \;\; \text{for } l = 1, \ldots \layer \label{eq:mlp-lagrangian-optimality-mu} \\
   & \frac{\partial \mathcal{L}}{\partial \delta_l} = 0 \implies z_{l} = f(x_{l-1},\Theta_{l}) \;\; \text{for } l = 1, \ldots \layer \label{eq:mlp-lagrangian-optimality-delta}
\end{align}
Notice that sub-conditions \eqref{eq:mlp-lagrangian-optimality-mu} and \eqref{eq:mlp-lagrangian-optimality-delta} carry out exactly the forward-pass of the Algorithm \ref{alg:backprop-algorithm}. Merging sub-conditions \eqref{eq:mlp-lagrangian-optimality-x_l} and \eqref{eq:mlp-lagrangian-optimality-z_l} gives
\begin{equation*}%\label{eq:mlp-backprop-relation}
    \delta_{l} = \sigma_{l}'(z_{l}) \frac{\partial f}{\partial x_{l}} \delta_{l+1}  \;\; \text{for } l = 1 \ldots \layer-1 \,\, ,
\end{equation*}
which is the backward-pass in Algorithm \ref{alg:backprop-algorithm} for computing the partial derivatives. Sub-conditions \eqref{eq:mlp-lagrangian-optimality-x_L-1} and \eqref{eq:mlp-lagrangian-optimality-w_l} on the other hand are computational steps for taking the (sub-)gradient update step (cf. \cite{higham2019deep}). 

% \section{Subdifferential properties}\label{sec:subdiff-properties}
% We verify condition \eqref{eq:subdifferential-equivalence} directly with the following lemma and proof.
% \begin{lemma}[Subdifferential characterisation]
% Suppose the function $\Psi:\mathbb{R}^n \rightarrow \mathbb{R} \cup \{ \infty \}$ is a proper, lower semi-continuous and convex function. Then the inclusion $z \in \partial \left( \frac12 \| \cdot \|^2 + \Psi\right)(x)$ is equivalent to \eqref{eq:subdifferential-equivalence}.
% \end{lemma}
% %\begin{proof}
% %\end{proof}

\subsection{Related works}\label{sec:related-works}
To tackle the major limitations of the classical lifted training approach, several other recent works also formulate training of lifted networks differently. Fenchel lifted neural network \citep{gu2020fenchel} and Lifted Proximal Operator Machine (LPOM) \citep{li2019lifted} both extend the observation in \cite{zhang2017convergence} to a broader class of activation functions and demonstrate improved learning performance. In this section, we briefly present the approaches proposed in these works and point out the their connections and differences to the proposed lifted Bregman framework.

Fenchel lifted neural networks \citep{gu2020fenchel} extend the classical lifted network approach \citep{askari2018lifted} as discussed in Section \ref{sec:lifted} by reformulating the learning problem \eqref{eq:mlp-objective} of a neural network with affine-linear function $f(x_{l-1},\Theta_{l}) = W_{l}^\top x_{l-1} + b_l$ as
\begin{align*}
    & \min_{\Xs, \{W_l, b_l\}_{l=1}^{\layer}} \;\; \sum_{i=1}^{s} \ell(y^i, W_{\layer}^\top x_{\layer-1}^{i} + b_{\layer}) \\
    %& \text{\;\;\;\;\;\;\;s.t. } x^i_{0} = x \, \nonumber \\
    & \text{\;\;\;\;\;\;\;s.t. } B_{l}(x_{l}, W_{l}^{\top}x_{l-1}^{i} + b_l) \leq 0 \, , \; \text{for} \; l=0 \ldots \layer-1 \,. \nonumber
\end{align*}%
where the equality constraints in Problem \eqref{eq:mlp-constrained} are converted into constraints on a collection of pre-defined biconvex functions $B_{l}$. 
For one fixed data pair $(x,y)$, the constraint $B_{l}(x_{l}, W_{l}^\top x_{l-1}) \leq 0$ provides the characterisation of an activation function via
\begin{equation*}
x_{l}^{i} =\sigma(W_{l}^\top x_{l-1}^{i} + b_l)  \qquad \Leftrightarrow \qquad B_{l}(x_{l}^{i},W_{l}^\top x_{l-1}^{i} + b_l) \leq 0 \, .
\end{equation*}
\noindent The Lagrange formulation of the constrained minimisation problem yields
\begin{equation*}%\label{eq:fenchel-lifted}
    \min_{\Xs, \{W_l,b_l\}_{l=1}^{\layer}} \;\; \sum_{i=1}^{s} \left [ \ell(y^i, W_{\layer}^\top x^i_{\layer-1}) + \sum_{l=1}^{\layer} \lambda_l B_{l}(x_{l+1}, W_{l+1}^\top x_{l}) \right] \,\, ,
    %& \text{\;\;\;\;\;\;\;s.t. } x^i_{0} = x \,, \nonumber
\end{equation*}%
where $\lambda_l$ are the Lagrange multipliers. 
The bi-convex functions $B_{l}$ essentially take on the same role as the Bregman loss function \eqref{eq:bregman-loss}, but without making the connection to Bregman distances and without many of the theoretical results provided in this paper. In the case of the ReLU activation function $\sigma(z) = \max(z, 0)$, the functions $B_{l}$ are defined as
\begin{align*}
    B_{l}(v, u) &= \begin{cases} \frac12 v^2 + \frac12 u_{+}^2 -uv  &\text{if } v \geq 0 \\
    \infty & \text{otherwise} \end{cases} \, .
\end{align*}
Using Item 4 in Theorem \ref{thm:bregman-loss}, it is not difficult to see that this is the scalar-valued equivalent of  the Bregman loss function
\begin{align*}
    \bregmanloss(x, z) &= \begin{cases} \frac12 \| x \|^2 + \frac12 \| \max(z, 0) \|^2 - \langle z,x \rangle &\text{if } x \geq 0  \\ \infty & \text{otherwise} \end{cases} \, .
\end{align*}
However, without using generalised Bregman distances and proximal activation functions, the derivation of $B_{l}$ in \cite{gu2020fenchel} requires deducing the energy from the optimality system. Further, derivations for activation functions other than ReLU seem possible but not as straightforward. %Indeed the paper briefly noted the possibility of re-writing sigmoid or tanh activation functions with constraints on $B_{l}$, but the derivation and discussion is limited only to the case of using ReLU activation function.
Moreover, the analysis in \cite{gu2020fenchel} requires the activation function $\sigma$ to be strictly monotone, while such strong restrictions are not required in this work. 
%Alternatively, the $B_l$ characterisation can be equivalently described from a Bregman perspective. As an example, we consider the case of a neural network with ReLU activation function. Using property 2 from Theorem \ref{thm:bregman-loss}, we have the following relation:
%\begin{align*}
%    \bregmanloss(y,z) &= \frac12 \| y - \sigma(z) \|^2 + D_\Psi^{z - \sigma(z)}(y, \sigma(z)) \\
%    &= \frac12 \| y - \sigma(z) \|^2 + \Psi(y) + \Psi^*(z-\sigma(z)) - \langle z-\sigma(z),y \rangle \\
%    &= \frac12 \| y \|^2 + \frac12 \| \sigma(z) \|^2 - \langle y, \sigma(z) \rangle + \Psi(y) + \Psi^*(z-\sigma(z)) - \langle z,y \rangle + \langle  \sigma(z),y \rangle \\
%    &= \frac12 \| y \|^2 + \frac12 \| \sigma(z) \|^2 + \Psi(y) + \Psi^*(z-\sigma(z)) - \langle z,y \rangle \\
%    &=  \frac12 \| y \|^2 + \frac12 \| \sigma(z) \|^2 - \langle z,y \rangle \text{ if } v \geq 0 
%\end{align*}
%and $\infty$ otherwise. This is characterised in the work of \cite{gu2020fenchel} by:

LPOM on the other hand takes a similar approach to \cite{gu2020fenchel} through the use of proximal operators. The authors define %$f(x_l^{i})$ as: 
\[f(x) := \int_0^x (\sigma^{-1}(y) -y) dy \qquad \text{and} \qquad g(x) := \int_0^x (\sigma(y) - y) dy \]
such that the proximal map of $f$ as defined in Definition \ref{def:proximal-map} is $\sigma$, and the proximal map of $g$ is $\sigma^{-1}$.
%\begin{equation*}
%    \text{prox}_f(y) = \argmin_{x} f(x) + \frac{1}{2}\|x-y\|^2 = \sigma(y)\,\,.
%\end{equation*}%
The authors of \cite{li2019lifted} propose to approximate \eqref{eq:mlp-constrained} with
\begin{equation*}
 \min_{\Xs,\{W_l,b_l\}_{l=1}^{\layer}} \;\; \sum_{i=1}^{s} \ell(y^i, W_{\layer}^\top x^i_{\layer-1}+b_{\layer}) + \sum_{l=1}^{\layer-1} \mu_l ( f(x_{l}^{i}) + g(W_{l}^\top x_{l-1}^{i}) + \frac12 \|x_{l}^{i} - (W_{l}^\top x_{l-1}^{i} + b_l)\|^2 ) \, , \\
\end{equation*}
such that the optimality condition of the objective with respect to $x_l^i$ reads
\begin{equation*}
    0 = \mu_{l}(\sigma^{-1}(x_{l}^{i}) - (W_{l}^\top x_{l-1}^{i} + b_l)) + \mu_{l+1}(W_{l+1})(\sigma(W_{l+1}^\top x_{l}^{i} +b_l) - x_{l+1}^{i}) \,\,\, .
\end{equation*}
%This modification approximates the original constrained Problem \eqref{eq:mlp-constrained} and is essential to ensure the nonlinear equality constrains $x_{l}^{i} = \sigma(W_{l}^\top x_{l-1}^{i} + b_l)$ are fulfilled.
Similar to the proposed approach, this approach ensures that the network constraints $x_{l}^{i} = \sigma(W_{l}^\top x_{l-1}^{i} + b_l)$ are satisfied. However, like the Fenchel lifted neural network network approach, the LPOM work requires more restrictive assumptions such as invertibility of the activation function $\sigma$, and provides fewer mathematical insights into the newly defined penalisation functions.

\section{Numerical realisation}\label{app:numerical-realisation}
In this section we present two other suitable deterministic strategies for the computational minimisation of \eqref{eq:bregman-lifted}.

\subsection{Coordinate descent}
The function $E$ in Section \ref{sec:numerical-realisation} is convex in each individual variable if all other variables are kept fixed. It therefore makes sense to consider alternating minimisation approaches, also known as \emph{coordinate descent}, for the minimisation of $E$, i.e.
\begin{subequations}
\begin{align}
    \Theta_l^{k + 1} &= \argmin_{\Theta_l} E\left(\Theta_1^{k + 1}, \ldots, \Theta_{l - 1}^{k + 1}, \Theta_l, \Theta_{l + 1}^k, \ldots, \Theta_{\layer}, x_1^k, \ldots, x_{\layer - 1}^k\right) \, ,\label{eq:coordinate-descent-1}
    \intertext{respectively}
    x_j^{k + 1} &= \argmin_{x_j} E\left(\Theta_1^{k + 1}, \ldots, \Theta_\layer^{k + 1}, x_1^{k + 1}, \ldots, x_{j - 1}^{k + 1}, x_j, x_{j + 1}^k, \ldots, x_{\layer - 1}^k \right) \, ,
\end{align}\label{eq:coordinate-descent-all}%
\end{subequations}
for $l \in \{1, \ldots, \layer\}$ and $j \in \{1, \ldots, \layer - 1\}$. Obviously, the order of updating the individual parameters is completely arbitrary and can be replaced with permutations. However, we want to point out that it does make sense to alternate with respect to the $\mathbf{\Theta}$ and $\mathbf{X}$ block, and to further alternate between the block of even and odd indices within the $\mathbf{X}$ block. Suppose $\layer$ is even, then this alternating optimisation between blocks can be written as
\begin{subequations}
\begin{align}
    \mathbf{\Theta}^{k + 1} &= \argmin_{\Theta_1, \ldots, \Theta_\layer} E\left(\Theta_1, \ldots, \Theta_{\layer}, x_1^k, \ldots, x_{\layer - 1}^k\right) \, , \\
    \{ x_{2l - 1}^{k + 1} \}_{l = 1}^{\frac{\layer}{2}} &= \argmin_{x_1, x_3, \ldots, x_{\layer - 1}} E\left(\Theta_1^{k + 1}, \ldots, \Theta_{\layer}^{k + 1}, x_1, x_2^k, x_3, \ldots, x_{\layer - 2}^k, x_{\layer - 1}\right) \, ,\\
    \{ x_{2l}^{k + 1} \}_{l = 1}^{\frac{\layer}{2}} &= \argmin_{x_2, x_4, \ldots, x_{\layer - 2}} E\left(\Theta_1^{k + 1}, \ldots, \Theta_{\layer}^{k + 1}, x_1^{k + 1}, x_2, x_3^{k + 1}, \ldots, x_{\layer - 2}, x_{\layer - 1}^{k + 1}\right) \, .
\end{align}\label{eq:coordinate-descent-blocks1}%
\end{subequations}
This alternating scheme has the advantage that the individual optimisation problems are still convex, but whilst the three blocks have to be updated sequentially, every variable in each of the three blocks can be updated in parallel. This becomes evident if we write down the individual optimisation problems, i.e.
\begin{subequations}
\begin{align}
    \Theta_l^{k + 1} {} = {} &\argmin_{\Theta_l} \left\{ \left( \frac12 \| \cdot \|^2 + \Psi \right)^\ast\left(f(x_{l - 1}^k, \Theta_l) \right) - \left\langle x_l^k, f(x_{l - 1}^k, \Theta_l) \right\rangle \right\} \, , \label{eq:coordinate-descent-first-block}\\
    x_i^{k + 1} {} = {} &\argmin_{x_i} \left\{ \left( \frac12 \| \cdot \|^2 + \Psi\right)(x_i) - \left\langle x_i, f(x_{i - 1}^k, \Theta_i^{k + 1}) \right\rangle \right. \label{eq:coordinate-descent-second-block}\\
    &\qquad \qquad \left. + \left( \frac12 \| \cdot \|^2 + \Psi\right)^\ast\left(f(x_i, \Theta_{i + 1}^{k + 1}) \right) - \left\langle f(x_i, \Theta_{i + 1}^{k + 1}), x_{i + 1}^k \right\rangle \right\} \, , \nonumber\\
    x_j^{k + 1} {} = {} &\argmin_{x_j} \left\{ \left( \frac12 \| \cdot \|^2 + \Psi\right)(x_j) - \left\langle x_j, f(x_{j - 1}^{k + 1}, \Theta_j^{k + 1}) \right\rangle \right. \label{eq:coordinate-descent-third-block}\\
    &\qquad \qquad \left. + \left( \frac12 \| \cdot \|^2 + \Psi\right)^\ast\left(f(x_j, \Theta_{j + 1}^{k + 1}) \right) - \left\langle f(x_j, \Theta_{j + 1}^{k + 1}), x_{j + 1}^{k + 1} \right\rangle \right\} \, ,\nonumber
\end{align}\label{eq:coordinate-descent-blocks2}%
\end{subequations}
for $l \in \{1, \ldots, \layer\}$, $i \in \{1, 3, \ldots, \layer - 1\}$ and $j \in \{2, 4, \ldots, \layer - 2\}$.

Note that the individual updates in \eqref{eq:coordinate-descent-all} and in \eqref{eq:coordinate-descent-blocks1}, respectively \eqref{eq:coordinate-descent-blocks2}, are still implicit. If we, for instance, want to make the updates with respect to $\Theta$ explicit, we can modify \eqref{eq:coordinate-descent-1} to
\begin{equation}
\begin{aligned}
    &\Theta_l^{k + 1} = \argmin_{\Theta_l} \left\{ E\left(\Theta_1^{k + 1}, \ldots, \Theta_{l - 1}^{k + 1}, \Theta_l, \Theta_{l + 1}^k, \ldots, \Theta_{\layer}, x_1^k, \ldots, x_{\layer - 1}^k\right) \right. \\
    &\left. + D_J\left( \left( \Theta_1^{k + 1}, \ldots, \Theta_{l - 1}^{k + 1}, \Theta_l, \Theta_{l + 1}^k, \ldots, \Theta_{\layer}^k \right), \left( \Theta_1^{k + 1}, \ldots, \Theta_{l - 1}^{k + 1}, \Theta_l^k, \Theta_{l + 1}^k, \ldots, \Theta_{\layer}^k \right) \right) \right\} \, ,
\end{aligned}\label{eq:bregman-coordinate-descent}
\end{equation}
where $D_J$ denotes the Bregman distance w.r.t. the function $J(\mathbf{\Theta}) = \sum_{l = 1}^\layer \frac{1}{2\tau_{\Theta_l}} \| \Theta_l \|^2 - E(\mathbf{\Theta}, \mathbf{X}^k)$. Similar modifications can be deployed for the block-coordinate descent scheme \eqref{eq:coordinate-descent-blocks1}, respectively \eqref{eq:coordinate-descent-blocks2}, and the auxiliary parameters $\mathbf{X}$. In Appendix \ref{sec:admm}, we briefly address how to computationally solve the individual sub-problems \eqref{eq:coordinate-descent-first-block}, \eqref{eq:coordinate-descent-second-block} and \eqref{eq:coordinate-descent-third-block}, with a special case of alternating minimisation approach: the alternating direction method of multipliers (ADMM), see \cite{gabay1983chapter}.

\subsection{Alternating direction method of multipliers}\label{sec:admm}

In the previous section we have explored how to break the optimisation problem into different blocks via coordinate descent. Each of these blocks could be solved with proximal gradient descent as described in Section \ref{sec:prox-grad}. However, we can alternatively solve these blocks with ADMM instead. Beginning with sub-problem \eqref{eq:coordinate-descent-first-block}, we can introduce auxiliary variables $z_l = f(x_{l - 1}^k, \Theta_l)$ to replace \eqref{eq:coordinate-descent-first-block} with the constrained formulation
\begin{align*}
    (\Theta_l^{k + 1}, z_l) {} = {} &\argmin_{\Theta, z} \left\{ \left( \frac12 \| \cdot \|^2 + \Psi \right)^\ast\left(z \right) - \left\langle x_l^k, z \right\rangle \right\} \quad \text{subject to} \quad z = f(x_{l - 1}^k, \Theta) \, ,
\end{align*}
for all $l \in \{1, \ldots, \layer\}$. By introducing Lagrange multipliers $\{ \mu_l \}_{l = 1}^\layer$, we can transform these problems into saddle-point problems of the form
\begin{align*}
    (\Theta_l^{k + 1}, z_l, \mu_l) {} = {} &\argmin_{\Theta, z}\max_{\mu} \left\{ \left( \frac12 \| \cdot \|^2 + \Psi \right)^\ast\left(z \right) - \left\langle x_l^k, z \right\rangle + \left\langle \mu,  z - f(x_{l - 1}^k, \Theta) \right\rangle \right\} \, .
\end{align*}
ADMM computationally solves this saddle-point problem by alternatingly minimising and maximising the associated augmented Lagrangian
\begin{align*}
    \mathcal{L}_{l}^k(\Theta, z; \mu) = \left( \frac12 \| \cdot \|^2 + \Psi \right)^\ast\left(z \right) - \left\langle x_l^k, z \right\rangle + \left\langle \mu,  z - f(x_{l - 1}^k, \Theta) \right\rangle + \frac{\delta_l}{2} \left\| z - f(x_{l - 1}^k, \Theta) \right\|^2 \, ,
\end{align*}
for a positive parameter $\delta_l$, i.e.
\begin{align*}
    \Theta_l^{j + 1} &= \argmin_{\Theta} \mathcal{L}_{l}^k(\Theta, z_l^j; \mu_l^j) \, , \\
    z_l^{j + 1} &= \argmin_{z} \mathcal{L}_{l}^k(\Theta_l^{j + 1}, z; \mu_l^j) \, , \\
    \mu_l^{j + 1} &= \argmax_{\mu} \mathcal{L}_{l}^k(\Theta_l^{j + 1}, z_l^{j + 1}; \mu) - \frac{1}{2\delta_l} \| \mu - \mu_l^j \|^2 \, .
    %\mu_l^j + \delta_l \left( z_l^{j + 1} - f(x_{l - 1}^k, \Theta_l^{j + 1}) \right) \, .
\end{align*}
In the context of \eqref{eq:coordinate-descent-first-block}, these updates read
\begin{align*}
    \Theta_l^{j + 1} &= \argmin_{\Theta} \left\{ \frac{\delta_l}{2} \left\| f(x_{l - 1}^k, \Theta) - z_l^j \right\|^2 - \langle \mu_l^j , f(x_{l - 1}^k, \Theta) \rangle \right\} \, , \\
    z_l^{j + 1} &= \text{prox}_{\frac{1}{\delta_l} (\frac12 \| \cdot \|^2 + \Psi)^\ast} \left( \frac{1}{\delta_l} \left( x_l^k - \mu_l^j + \delta_l \, f(x_{l - 1}^k, \Theta_l^{j + 1}) \right) \right) \, , \\
    \mu_l^{j + 1} &= \mu_l^j + \delta_l \left( z_l^{j + 1} - f(x_{l - 1}^k, \Theta_l^{j + 1}) \right) \, .
\end{align*}
Please note that the proximal map $\text{prox}_{\frac{1}{\delta_l} (\frac12 \| \cdot \|^2 + \Psi)^\ast}$ can easily be expressed via the extended Moreau decomposition \cite[Theorem 6.45]{beck2017first}, i.e. $\delta_l \, \text{prox}_{\frac{1}{\delta_l} (\frac12 \| \cdot \|^2 + \Psi)^\ast}(x/\delta_l) = x - \text{prox}_{\delta_l (\frac12 \| \cdot \|^2 + \Psi)}(x)$. We further know $\text{prox}_{\delta_l (\frac12 \| \cdot \|^2 + \Psi)}(x) = \text{prox}_{\frac{\delta_l}{1 + \delta_l} \Psi}(x / (1 + \delta_l) )$; hence, we observe 
\begin{align*}
    z_l^{j + 1} &= \text{prox}_{\frac{1}{\delta_l} (\frac12 \| \cdot \|^2 + \Psi)^\ast} \left( \frac{1}{\delta_l} \left( x_l^k - \mu_l^j + \delta_l \, f(x_{l - 1}^k, \Theta_l^{j + 1}) \right) \right) \, , \\
    &= \frac{1}{\delta_l} \left(  x_l^k - \mu_l^j + \delta_l \, f(x_{l - 1}^k, \Theta_l^{j + 1}) - \text{prox}_{\frac{\delta_l}{1 + \delta_l} \Psi}\left( \frac{1}{1 + \delta_l} \left(  x_l^k - \mu_l^j + \delta_l \, f(x_{l - 1}^k, \Theta_l^{j + 1}) \right) \right) \right) \, ,
\end{align*}
which is why we can compute the update solely based on a scaled version of the proximal map of $\Psi$, which is the activation function of the neural network. 
\begin{example}
In analogy to Example \ref{exm:proximal-gd} we design a feed-forward network architecture with $\Psi = \chi_{\geq 0}$, implying $\text{prox}_\Psi(z) = \max(z, 0)$, and $f(x, \Theta_l) = W^\top_l x + b_l$, for $\Theta_l = (W_l, b_l)$. Computing a solution to \eqref{eq:coordinate-descent-first-block} via ADMM yields the algorithm
\begin{align}
    W_l^{j + 1} &= \left(x_{l - 1}^k (x_{l - 1}^k)^\top \right)^{-1} x_{l - 1}^k \left( z_l^j + \frac{1}{\delta_l} \mu_l^j - b_l^j \right)^\top \, , \label{eq:admm-weight-update}\\
    b_l^{j + 1} &= z_l^j - \left( \frac{1}{\delta_l} \mu_l^j - (W_l^{j + 1})^\top x_{l - 1}^k \right) \, , \nonumber\\
    r_l^{j + 1} &= x_l^k - \mu_l^j + \delta_l \left( (W_l^{j + 1})^\top x_{l - 1}^k + b_l^{j + 1} \right) \, , \nonumber\\
    z_l^{j + 1} &= \frac{1}{\delta_l} \left( r_l^{j + 1} - \, \text{prox}_{\frac{\delta_l}{1 + \delta_l} \Psi}\left( \frac{r_l^{j + 1}}{1 + \delta_l}  \right) \right) \, , \nonumber\\
    \mu_l^{j + 1} &= \mu_l^j + \delta_l \left( z_l^{j + 1} - \left( (W_l^{j + 1})^\top x_{l - 1}^k + b_l^{j + 1} \right) \right) \, ,\nonumber
\end{align}
for $l \in \{1, \ldots, \layer\}$, $j \in \mathbb{N}$ and suitable initial values $b_l^0$, $z_l^0$ and $\mu_l^0$. It is important to point out that \eqref{eq:admm-weight-update} is not well defined for only one sample because $x_{l - 1}^k (x_{l - 1}^k)^\top$ is not invertible. However, we can overcome this issue by extending $E$ to include more samples, so that instead of requiring invertibility of $x_{l - 1}^k (x_{l - 1}^k)^\top$, we require invertibility of $\sum_{i = 1}^s x_{i, l-1}^k (x_{i, l - 1}^k)^\top$, which can be achieved if one has sufficiently distinctive samples. Or we can add the term $\frac{\tau}{2} \| W_l - W_l^j \|^2_{\text{Fro}}$ to the augmented Lagrangian, with a positive multiple $\tau$, which doesn't affect the minimiser but guarantees that \eqref{eq:admm-weight-update} is well-defined by adding a positive multiple of the identity matrix to $x_{l - 1}^k (x_{l - 1}^k)^\top$, respectively $\sum_{i = 1}^s x_{i, l-1}^k (x_{i, l - 1}^k)^\top$.
\end{example}
In similar fashion to solving \eqref{eq:coordinate-descent-first-block}, we can also solve \eqref{eq:coordinate-descent-second-block} and \eqref{eq:coordinate-descent-third-block} via ADMM. Because \eqref{eq:coordinate-descent-second-block} and \eqref{eq:coordinate-descent-third-block} have identical structure but only different variables, we focus on \eqref{eq:coordinate-descent-second-block} without loss of generality. If we follow the same line-of-reasoning as before, the associated augmented Lagrangian reads
\begin{align*}
    \mathcal{L}_l^k(x, z; \mu) {} = {} &\left( \frac12 \| \cdot \|^2 + \Psi\right)(x) - \left\langle x, f(x_{l - 1}^k, \Theta_l^{k + 1}) \right\rangle + \left( \frac12 \| \cdot \|^2 + \Psi\right)^\ast(z) - \langle z, x_{l + 1}^k \rangle \\
    &+ \left\langle \mu, z - f(x, \Theta_{l + 1}^{k + 1}) \right\rangle + \frac{\delta_l}{2} \left\| z - f(x, \Theta_{l + 1}^{k + 1}) \right\|^2 \, .
\end{align*}
However, the challenge lies in the optimality condition for $x$, which even for fixed variables $z$ and $\mu$ doesn't lead to a closed-form solution for $x$ (unless $\Psi$ is constant). Alternatively, we can introduce another variable $v$ with constraint $v = x$ and formulate
\begin{align*}
    \mathcal{L}_l^k(x, v, z; \mu) {} = {} &\left( \frac12 \| \cdot \|^2 + \Psi\right)(x) - \left\langle x, f(x_{l - 1}^k, \Theta_l^{k + 1}) \right\rangle + \left( \frac12 \| \cdot \|^2 + \Psi\right)^\ast(z) - \langle z, x_{l + 1}^k \rangle \\
    &+ \left\langle \mu, \left( \begin{matrix} z \\ v \end{matrix} \right) - \left( \begin{matrix} f(v, \Theta_{l + 1}^{k + 1}) \\ x \end{matrix} \right) \right\rangle + \frac{\delta_l}{2} \left\| \left( \begin{matrix} z \\ v \end{matrix} \right) - \left( \begin{matrix} f(v, \Theta_{l + 1}^{k + 1}) \\ x \end{matrix} \right)  \right\|^2 \, .
\end{align*}
Alternatingly minimising and maximising this augmented Lagrangian yields the updates 
\begin{align*}
    x_l^{j + 1} &= \text{prox}_{\frac{1}{1 + \delta_l} \Psi}\left( \frac{\delta_i v_i^j + (\mu_l^j)_2 + f(x_{l - 1}^k, \Theta_l^{k + 1})}{1 + \delta_l} \right) \, , \\
    v_l^{j + 1} &= \argmin_v \left\{ \frac{\delta_l}{2} \left\| f(v, \Theta_{l + 1}^{k + 1}) - z_l^j \right\|^2 - \left\langle (\mu_l^j)_1, f(v, \Theta_{l + 1}^{k + 1}) \right\rangle \right\} \, , \\
    z_l^{j + 1} &= \text{prox}_{\frac{1}{\delta_l} ( \frac12 \| \cdot \|^2 + \Psi)^\ast}\left( f(v_l^{j + 1}, \Theta_{l + 1}^{k + 1}) + \frac{1}{\delta_l} \left( x_{l + 1}^k - (\mu_l^j)_1 \right) \right)\, , \\
    \mu_l^{j + 1} &= \left( \begin{matrix} z_l^{j + 1} \\ v_l^{j + 1} \end{matrix} \right) - \left( \begin{matrix} f(v_l^{j + 1}, \Theta_{l + 1}^{k + 1}) \\ x_l^{j + 1} \end{matrix} \right) \, ,
\end{align*}
for all $l \in \{ 1, 3, \ldots, \layer - 1 \}$. %We conclude this section on ADMM with another example in the spirit of Example \ref{exm:proximal-gd}. 

\subsection{Constrained optimisation}
Note that it is not necessarily known a-priori how to choose the scalar parameters $\{ \lambda_l \}_{l = 1}^\layer$ in \eqref{eq:bregman-lifted}. Alternatively, we can also re-formulate the constrained problem formulation \eqref{eq:mlp-constrained} to
\begin{align}
\begin{split}
&\min_{\Thetas,\Xs} \sum_{i=1}^{s} \ell(y^i, x_{\layer}^{i})  \\
\text{subject to} \qquad &\bregmanloss(x_l^i, f(x_{l - 1}^i, \Theta_l)) \leq 0 \quad \text{for all } \, l \in \{1, \ldots, \layer - 1\} \, ,
\end{split}\label{eq:mlp-constrained-new}
\end{align}
similar to \cite{gu2020fenchel}. Note that we can write \eqref{eq:mlp-constrained-new} as
\begin{align*}
    \min_{\Thetas,\Xs} \sum_{i=1}^{s} \ell(y^i, x_{\layer}^{i}) + \sum_{l = 1}^{\layer - 1} \chi_{\leq 0}\left( \bregmanloss(x_l^i, f(x_{l - 1}^i, \Theta_l)) \right) \, ,
\end{align*}
for the characteristic function $\chi_{\leq 0}$ over the non-positive orthant, which is equivalent to the Lagrange formulation
\begin{align}
    \min_{\Thetas,\Xs} \sup_{\{\lambda_l \}_{l = 1}^{\layer - 1}} \sum_{i=1}^{s} \ell(y^i, x_{\layer}^{i}) + \sum_{l = 1}^{\layer - 1} \left[ \lambda_l \bregmanloss(x_l^i, f(x_{l - 1}^i, \Theta_l)) - \chi_{\geq 0}(\lambda_l) \right] \, ,\label{eq:mlp-lagrange-new}
\end{align}
if we express $\chi_{\leq 0}$ in terms of its conjugate $\chi_{\geq 0}$. Albeit being concave in $\{\lambda_l \}_{l = 1}^{\layer - 1}$, the objective function is not convex in $\Thetas,\Xs$ in general, which is why swapping $\min$ and $\sup$ will not necessarily lead to the same optimisation problem. Nevertheless, we could employ coordinate descent or alternating direction method of minimisation techniques. We could for example solve \eqref{eq:mlp-lagrange-new} with an alternating minimisation and proximal maximisation algorithm of the form
\begin{align*}
    (\Thetas^{k + 1},\Xs^{k + 1}) &\in \argmin_{\Thetas, \Xs} \left\{ \sum_{i=1}^{s} \ell(y^i, x_{\layer}^{i}) + \sum_{l = 1}^{\layer - 1} \lambda_l^k \bregmanloss(x_l^i, f(x_{l - 1}^i, \Theta_l))  \right\} \, , \\
    \lambda_l^{k + 1} &= \max\left(0, \lambda_l^k + \tau \bregmanloss((x_l^i)^{k + 1}, f((x_{l - 1}^i)^{k + 1}, \Theta_l^k)) \right) \, ,
\end{align*}
where the minimisation with respect to $\Thetas,\Xs$ could be solved with any of the algorithms described in Section \ref{sec:numerical-realisation} or in this section.

\begin{figure}[!t]
    \centering
    \includegraphics[scale=0.25]{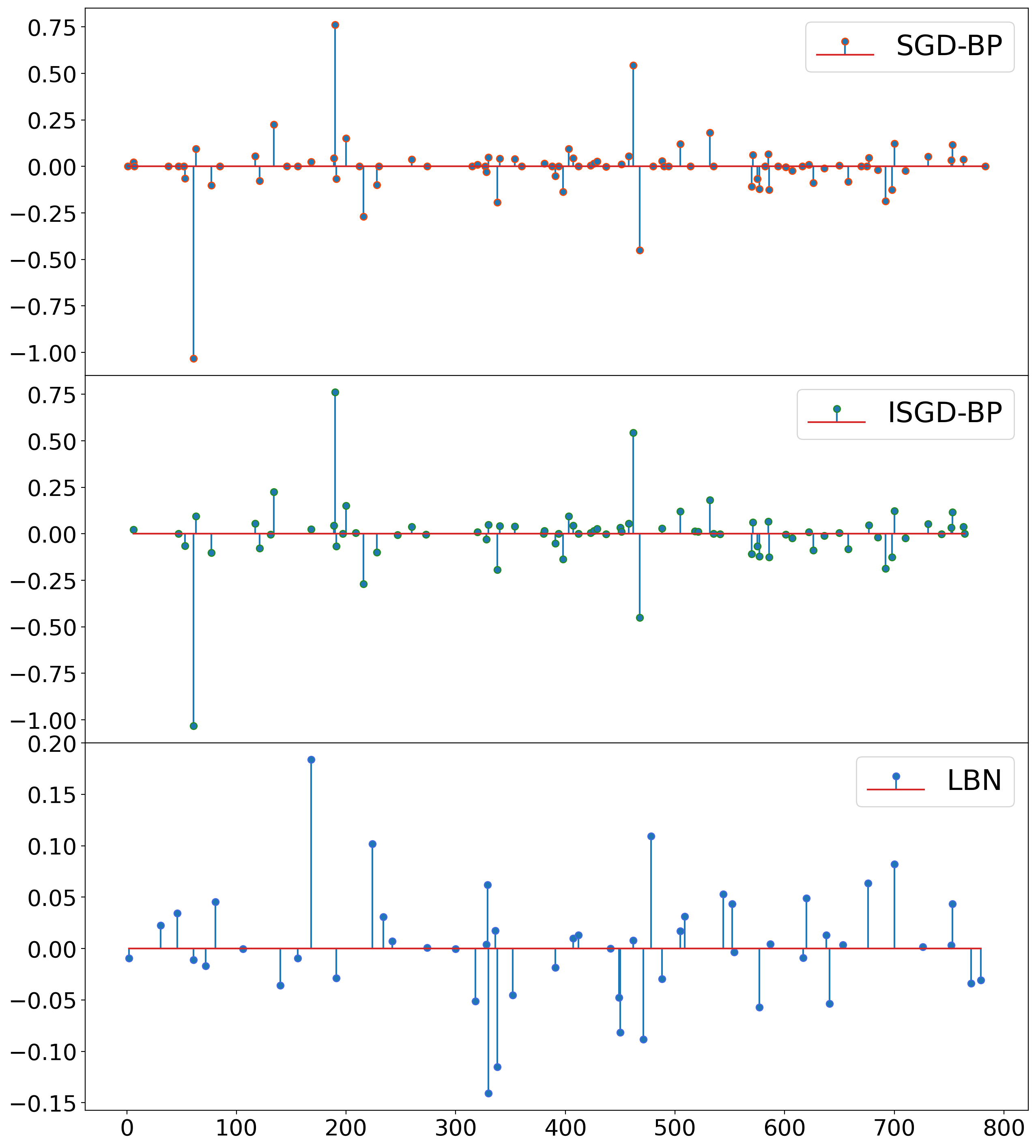}
    \caption{Visualisation of the sparse coding of the sandal image object from Figure \ref{figure:sparse-autoencoder-fmnist-denoised-train-samples} computed with the SGD-BP, ISGD-BP, and LBN, respectively. Note that the scaling of the sparse coefficients is slightly different for LBN in comparison to SGD-BP and ISGD-BP.} \label{figure:sparse-autoencoder-fmnist-denoised-sparse-code}
\end{figure}

\section{Numerical Results}
In the section we provide numerical results in addition to the ones provided in Section \ref{sec:numerical-realisation}. 

\subsection{Visualisations} \label{sec:additional-visualisation}
This section includes additional visualisations of the sparse autoencoder and sparse denoising autoencoder outputs. The compression of the input data is achieved by ensuring that only relatively few coefficients of the code are non-zero. As mentioned earlier, the advantage over conventional autoencoders is that the position of the non-zero coefficients can vary for every input. We also provide visualisations of codes from a sample input image computed with the SGD-BP, ISGD-BP and LBN approach in Figure \ref{figure:sparse-autoencoder-fmnist-denoised-sparse-code}. 

\begin{sidewaysfigure}[ht]
    \includegraphics[scale=0.35]{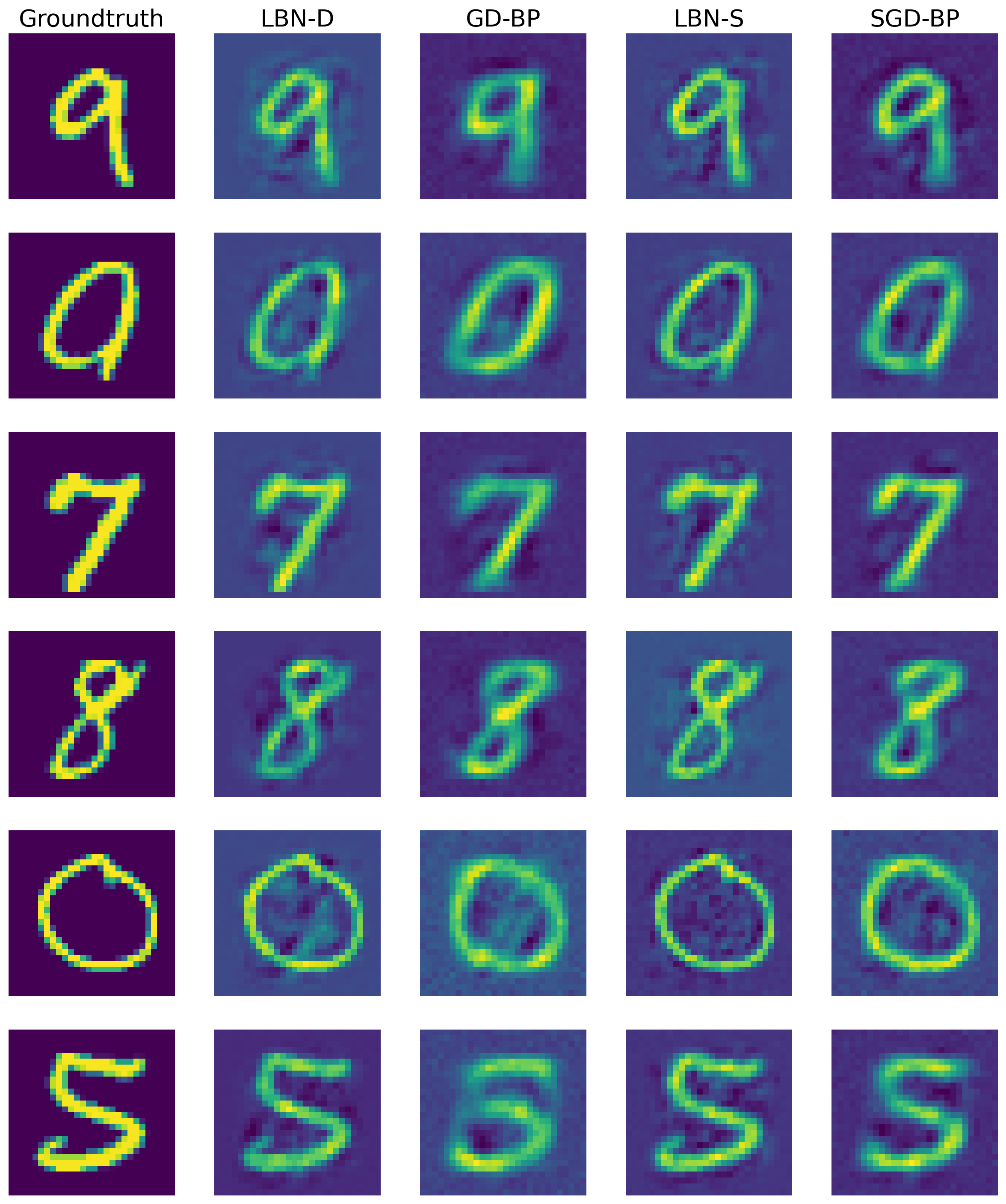}
    \hspace{10mm}
    \includegraphics[scale=0.35]{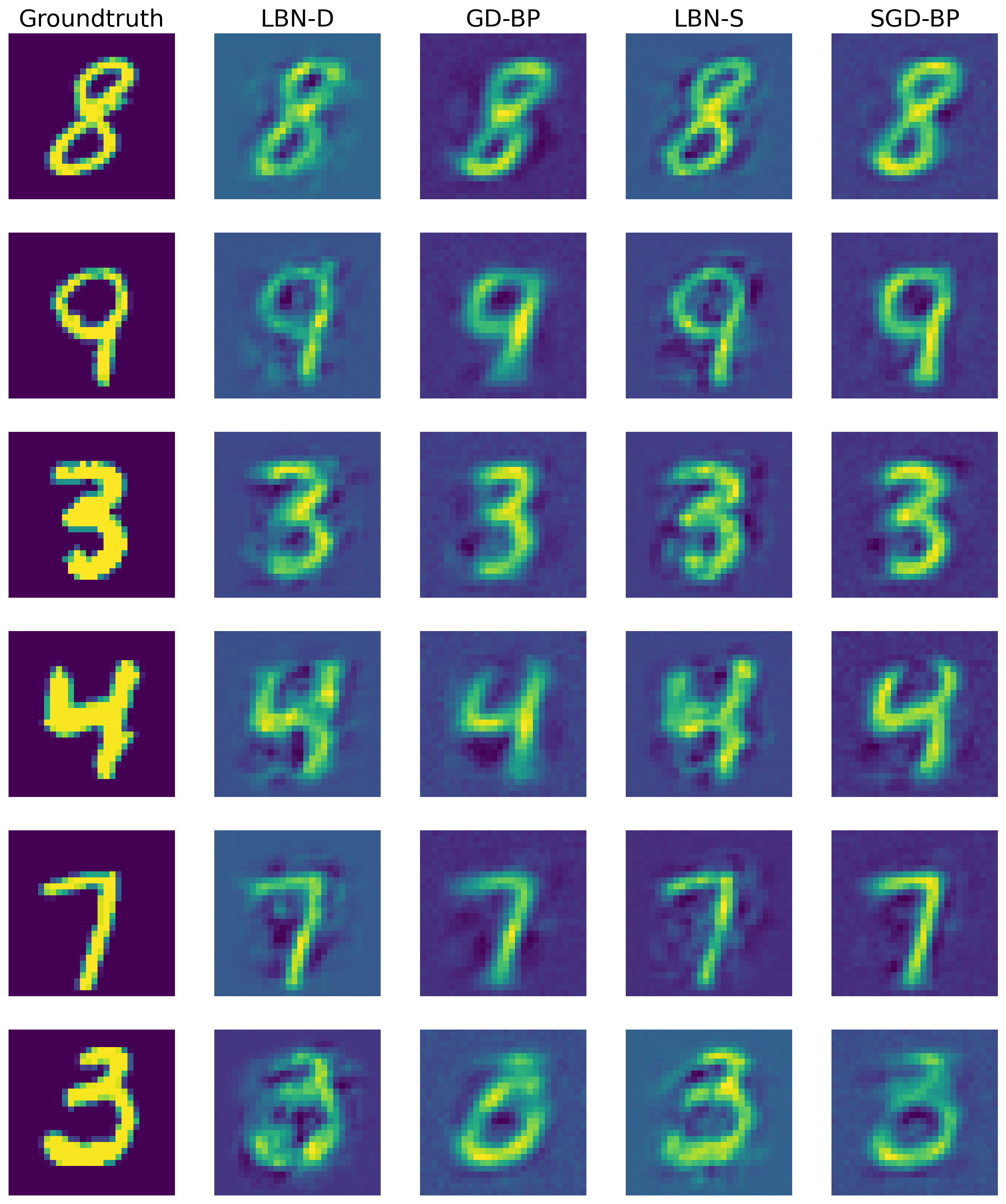}
    \caption{\textbf{Left}:This plot shows selected denoised MNIST training images computed with LBN-D, GD-BP, LBN-S and SGD-BP. \textbf{Right}: This plot shows selected denoising results on the MNIST validation dataset.}
    \label{fig:sparse_autoencoder_MNIST_recons_samples}
\end{sidewaysfigure}

\begin{sidewaysfigure}[ht]
    \includegraphics[scale=0.35]{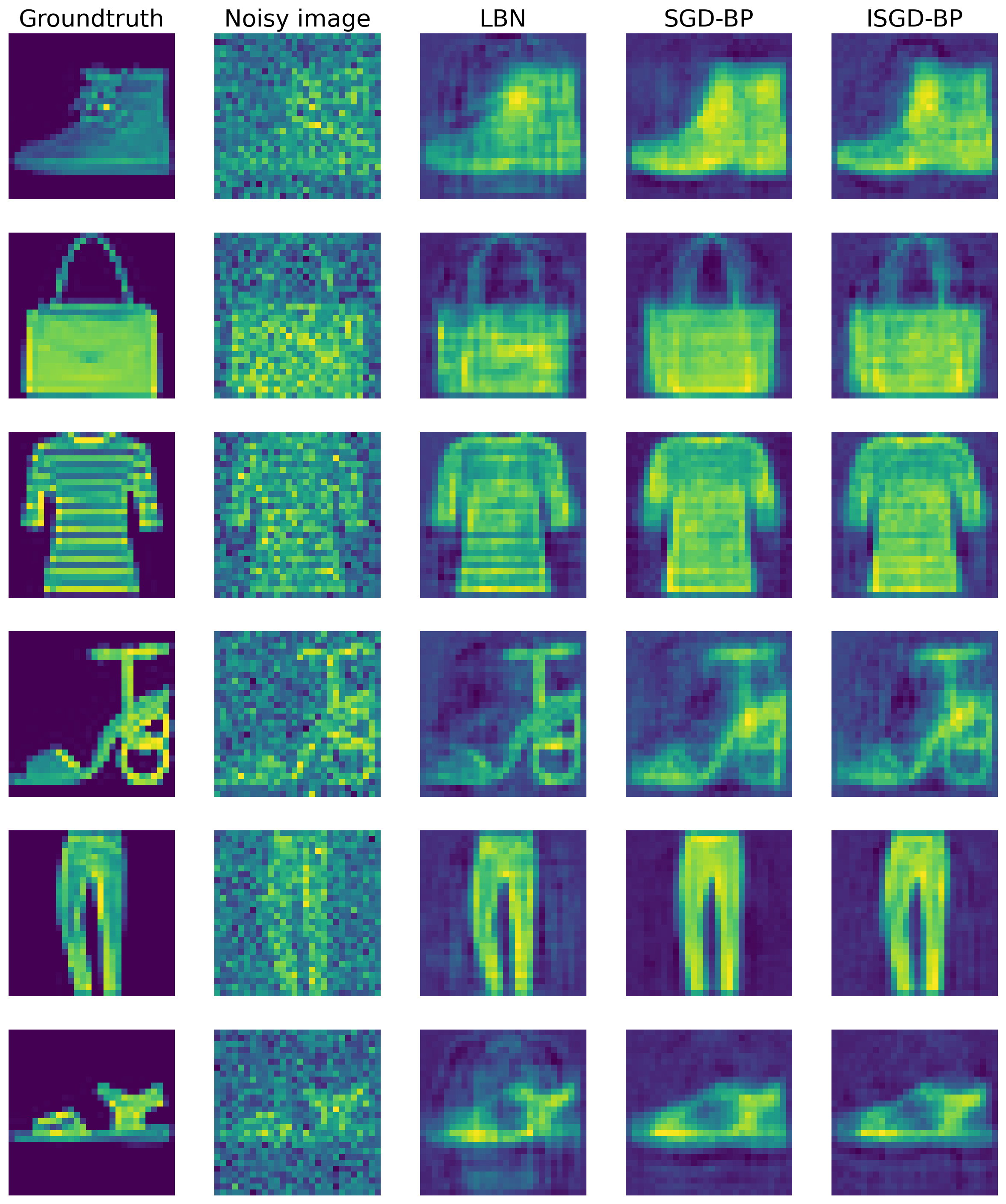}
    \hspace{1cm}
    \includegraphics[scale=0.35]{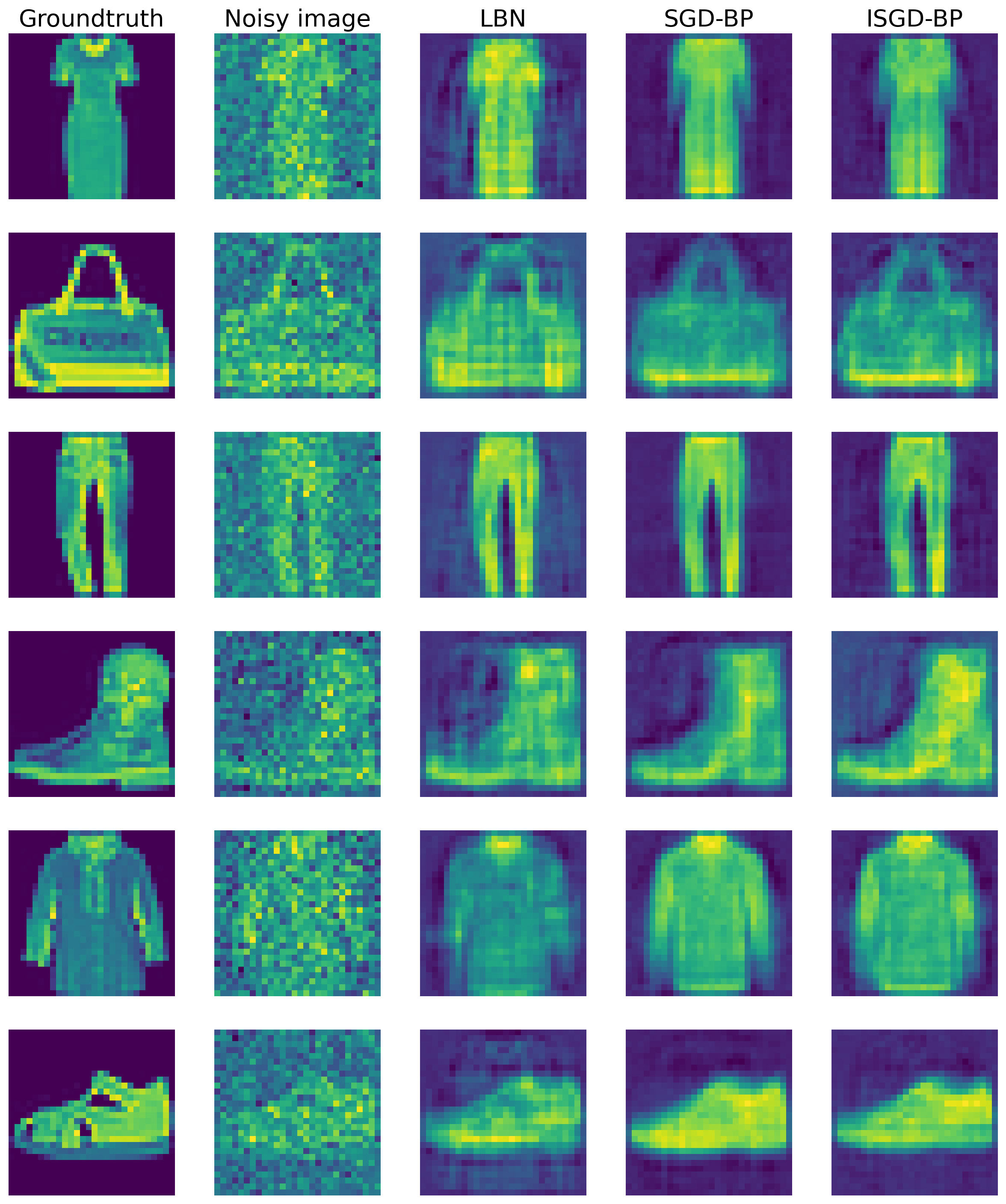}
    \caption{\textbf{Left}: Denoised images from the Fashion-MNIST-1K training dataset computed with LBN, SGD-BP and ISGD-BP for 50 epochs, along with the noise-free and noisy images. \textbf{Right}: Denoised images from the Fashion-MNIST-1K validation dataset computed with LBN, SGD-BP and ISGD-BP for 50 epochs, along with the noise-free and noisy images.}
    \label{fig:denoising_sparse_autoencoder_1K_FMNIST_denoised_samples}
\end{sidewaysfigure}

%\begin{figure}[!ht]
%    \centering
%    \includegraphics[scale=0.35]{figures/Numerical_Denoising_Sparse_Autoencoder_1K_FMNIST_Denoised_Train_Samples.png}
%    \caption{Results of denoised images from the Fashion-MNIST-1K training dataset under the LBN, the SGD-BP and the ISGD-BP training strategies after 50 epochs of training respectively, along with the clear and noisy images.}
%\end{figure}

%\begin{figure}[!ht]
%    \centering
%    \includegraphics[scale=0.35]{figures/Numerical_Denoising_Sparse_Autoencoder_1K_FMNIST_Denoised_Val_Samples.png}
%    \caption{Results of denoised images from the Fashion-MNIST-1K validation dataset under the LBN, the SGD-BP and the ISGD-BP training strategies after 50 epochs of training respectively, along with the clear and noisy images.}
%\end{figure}

\begin{sidewaysfigure}[ht]
    \includegraphics[scale=0.35]{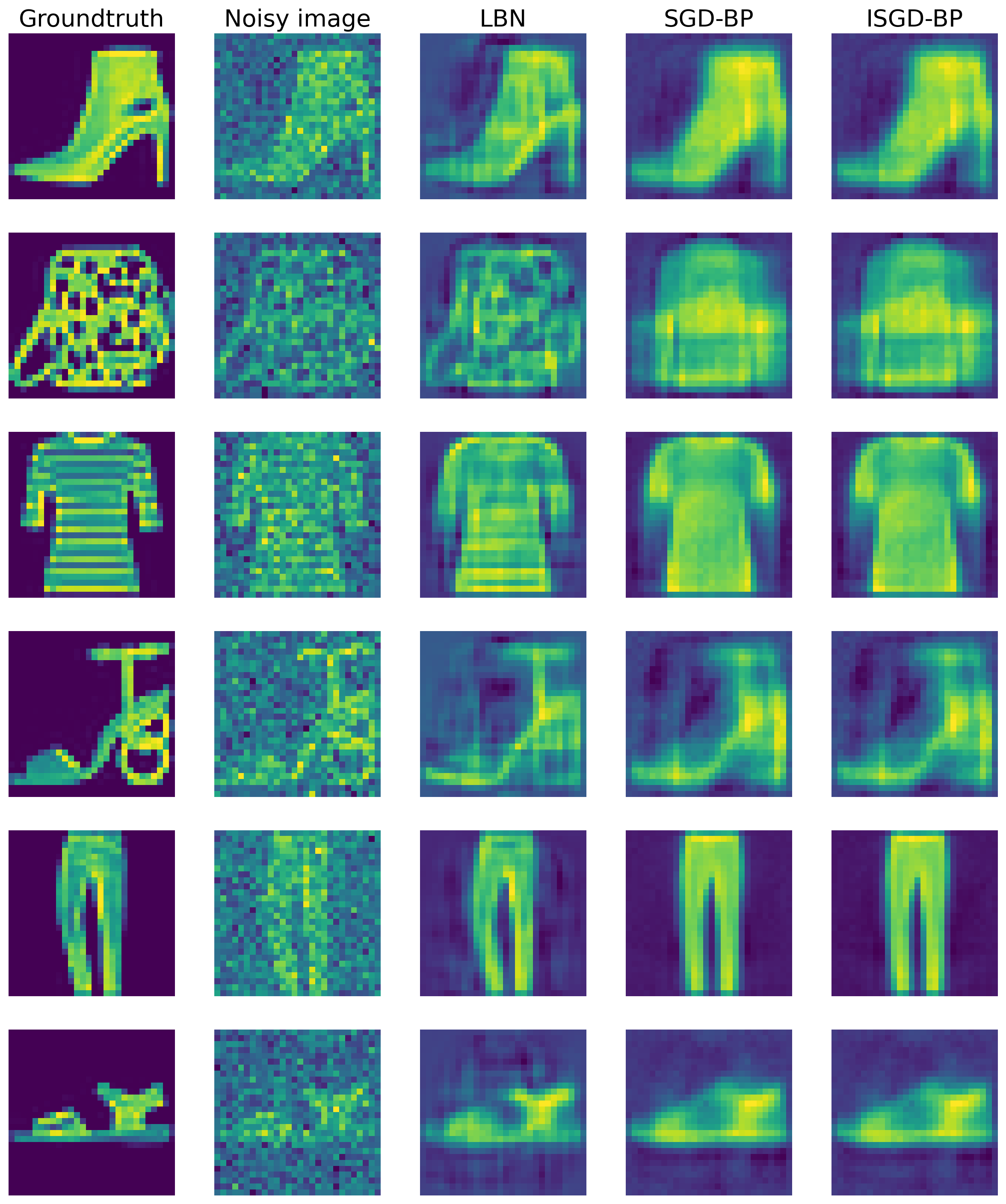}
    \hspace{10mm}
    \includegraphics[scale=0.35]{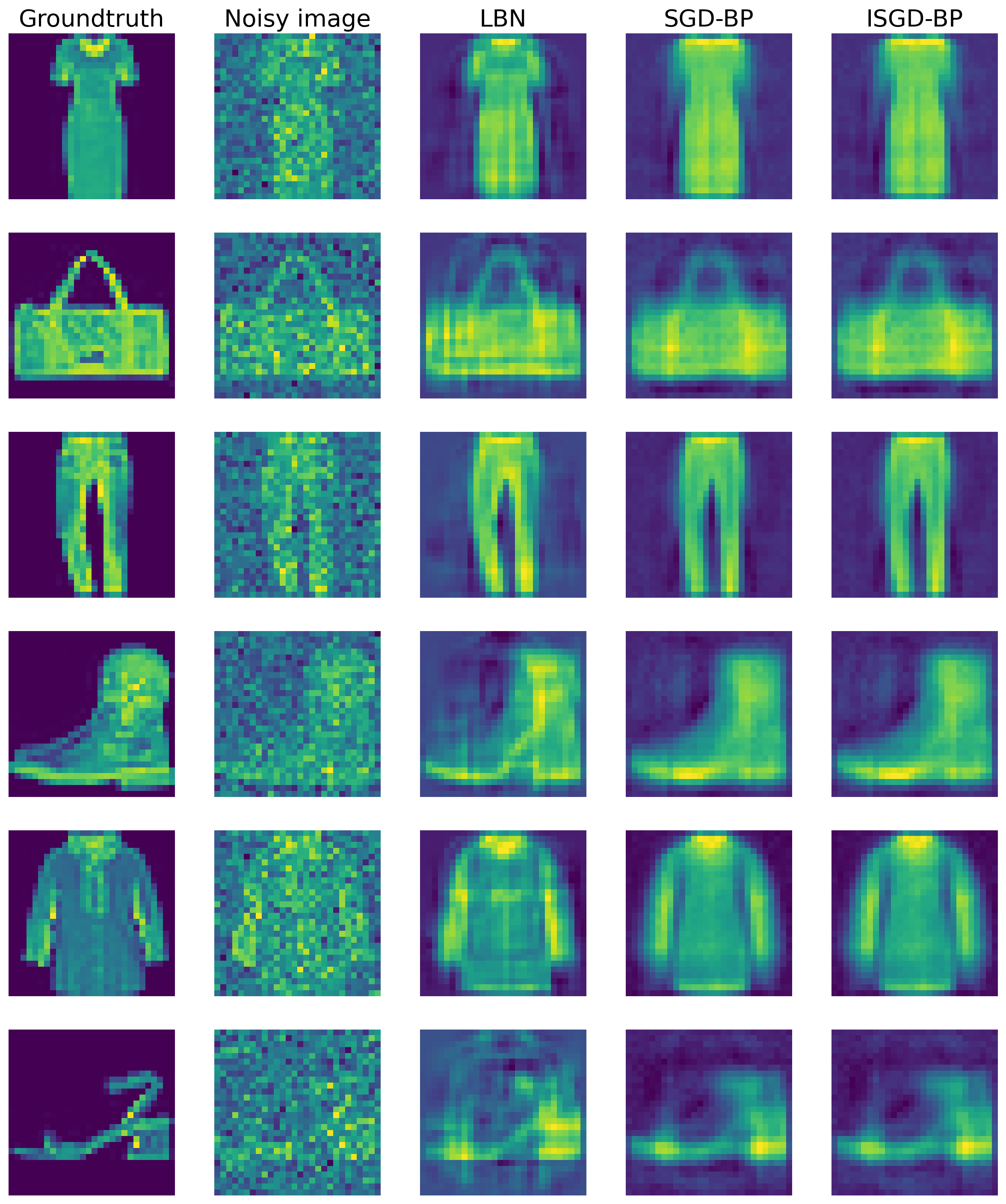}
    \caption{\textbf{Left}: Denoised images from the Fashion-MNIST-10K training dataset computed with LBN, SGD-BP (trained for 1500 epochs) and ISGD-BP, along with the noise-free and noisy images. \textbf{Right}: Denoised images from the Fashion-MNIST-10K validation dataset computed with LBN, SGD-BP (trained for 1500 epochs) and ISGD-BP, along with the noise-free and noisy images.}
    \label{fig:denoising_sparse_autoencoder_10K_FMNIST_denoised_samples}
\end{sidewaysfigure}

\end{document}